\documentclass[10pt]{amsart}

\usepackage{eucal,url,amssymb,stmaryrd,enumerate,amscd,verbatim}
\usepackage{amsmath}
\usepackage{xcolor}
\usepackage[pagebackref,colorlinks=true,linkcolor=cyan,citecolor=olive]{hyperref}
\usepackage{amsfonts}
\usepackage{amsthm}
\usepackage[margin=1in]{geometry}
\usepackage{eurosym,lipsum,latexsym,esint}
\usepackage{esint}
\usepackage[obeyFinal,textsize=footnotesize,textwidth=.5in]{todonotes}

\numberwithin{equation}{section}

\newtheorem{theorem}{Theorem}[section]

\newtheorem{lemma}[theorem]{Lemma}
\newtheorem{proposition}[theorem]{Proposition}
\newtheorem{prop}[theorem]{Proposition}
\newtheorem{corollary}[theorem]{Corollary}
\newtheorem{cor}[theorem]{Corollary}

\newtheorem{definition}[theorem]{Definition}

\numberwithin{equation}{section}

\newtheorem{innercustomthm}{}
\newenvironment{customthm}[1]{\renewcommand\theinnercustomthm{#1}\innercustomthm}
  {\endinnercustomthm}

\newcommand{\COG}{\mathcal{C}^\infty_0(\bG)}

\newcommand{\eps}{\varepsilon}

\newcommand{\Wsploc}{W^{s,p}_{loc}(\bG)}
\newcommand{\WsplocOm}{W^{s,p}_{loc}(\Om)}
\newcommand{\WspOm}{W^{s,p}(\Om)}
\newcommand{\WspG}{W^{s,p}(\bG)}
\newcommand{\tailspace}{\mathfrak{L}_s^{p-1}(\bG)}
\newcommand{\Dsp}{\mathcal{D}^{s,p}}
\newcommand{\DspOm}{\mathcal{D}^{s,p}(\Om)}
\newcommand{\DspG}{\mathcal{D}^{s,p}(\bG)}

\newcommand{\Rn}{\mathbb{R}^n}

\newcommand{\Om}{\Omega}

\newcommand{\Hn}{\mathbb{H}^n}

\newcommand{\bG}{\mathbf{G}}
\newcommand{\algg}{\mathfrak g}
\newcommand{\lap}{\mathcal{L}_s}
\newcommand{\frlap}{\mathcal{L}_{s,p}}

\newcommand{\norm}[1]{\lVert#1\rVert}

\newcommand{\volg}{\, d{g}}

\newcommand{\ddt}{\frac {d}{dt}}

\DeclareMathSymbol{\R}{\mathbin}{AMSb}{"52}
\DeclareMathSymbol{\N}{\mathbin}{AMSb}{"4E}
\newcommand{\G}[0]{{\bf{G}}}
\newcommand{\D}[0]{\mathcal{L}}

\newcommand{\paren}[1]{\left[ \, {#1}\, \right]}

\newcommand{\pl}{\left(}
\newcommand{\pr}{\right)}
\newcommand{\abs}[1]{\left| {#1}_{{}}\,\negthickspace\right|}

\newcommand{\tpwr}[1]{\left( #1\right)^{\langle p-1 \rangle}}
\newcommand{\fe}[1]{{\langle #1\rangle}}
\newcommand{\fep}{{\langle p-1\rangle}}

\newcommand{\feq}{{\langle q\rangle}}

\DeclareMathOperator*{\essinf}{ess\,inf}

\DeclareMathOperator*{\Tail}{T}

\DeclareMathOperator*{\osc}{osc}
\DeclareMathOperator*{\supp}{supp}
\DeclareMathOperator*{\loc}{loc}

\allowdisplaybreaks

\begin{document}

\begin{abstract}
    We  prove a Harnack inequality and H\"older continuity of solutions to non-linear nonlocal equations on homogeneous groups.  In addition, we show basic properties of the involved Sobolev spaces, which are crucial for the analysis of such operators.
\end{abstract}

\title[Regularity of Solutions to Non-Local Equations on Homogeneous Groups.]{Regularity of Solutions to Non-Local Equations on Homogeneous Groups.}

\author{Galen Novello }
\address[Galen Novello]{ Department of Mathematics and Statistics\\
University of New Mexico\\
Albuquerque, New Mexico, 87131-0001\\
}
\email{gnovello@unm.edu}

\author{Dimiter Vassilev}
\address[Dimiter Vassilev]{ Department of Mathematics and Statistics\\
University of New Mexico\\
Albuquerque, New Mexico, 87131-0001\\
}
\email{vassilev@unm.edu}

\date{\today}
\maketitle

\tableofcontents

\section{Introduction}
In this paper we  prove the local boundedness,  Harnack inequality and H\"older continuity of solutions to non-linear nonlocal equations on a homogeneous group $\bG$.  In addition, we show foundational properties of the involved spaces, which are crucial for the analysis of the considered non-local operators. Throughout the paper   a homogeneous group will be understood in the sense of \cite{Stein68}, \cite{KS69}, \cite{KV71}, and \cite{FS82}, i.e., $\G$ is a connected and simply connected Lie group endowed with a non-isotropic family of dilations.  The relevant definitions  are recalled  in section \ref{s:background}.

Nonlocal operators have been studied in the setting of Euclidean space relying on the various definitions of the fractional Laplacian and its generalizations. Following the Caffarelli-Silvestre extension method see \cite{CS} and \cite{Silvestre07},  a notable non-trivial generalization to H-type groups followed after the  CR conformal fractional sub-Laplacian $\D_s$ was introduced in \cite{BFM}  and \cite{FGMT} on the Heisenberg group. 
Motivated by the fractional Yamabe problem in Euclidean space, \cite{GonzalesQing13} and \cite{ChangGonzales-Adv11}, the papers  \cite{BFM}  and \cite{FGMT} considered the fractional CR Yamabe problem on the Heisenberg group. In particular, \cite{FGMT} showed that the CR conformality   leads to an extension problem involving the Kummer operator rather than the Bessel operator in a half-plane. On the other hand, on groups of H-type,  a fundamental fact proved in
\cite[proposition 4.1]{RTaim} and \cite[Theorem 1.2]{RT} for $0<s<1/2$, is the following remarkable M. Riesz type representation of the operator $\D_s$,
\begin{equation}\label{riesz}
   \D_s u(g)=\int_{\bG} \frac {u(g)-u(h)}{|h|^{Q+2s}}dh,
\end{equation}
up to a multiplicative constant involving the dimensions of the two layers of the H-type group. Here $|h|$ stands for the Kor\'anyi-Folland gauge. Using the heat equation approach in \cite{GT20Adv}, \cite{GTinter}, see also \cite{GT23PA}, formula \eqref{riesz} was extended to  the whole range $0<s<1$. Motivated by these  formulas \cite{Gar-Loi-Vas-2024-OP} and \cite{Gar-Loi-Vas-2024-FR} considered a version of  $\D_s$ by replacing the Kor\'anyi-Folland gauge in \eqref{riesz}  with a general homogeneous norm on groups of homogeneous type, which also constitute the ambient space of the current paper. 

In the linear case $p=2$, the fractional operator on homogeneous groups was considered in \cite{Gar-Loi-Vas-2024-OP} and \cite{Gar-Loi-Vas-2024-FR}.  We stress that the Lie groups considered here are not assumed to be Carnot groups.  Thus, we do not impose the H\"ormander condition \cite{Ho} on the invariant vector fields of homogeneity one, i.e., these vector fields might not generate the entire Lie algebra.  

Here we consider  more general non-local operators $\D_{s,p}$, which should be considered as a generalization of  the fractional $p$-Laplacian to the setting of homogeneous groups. Our results contain as special cases  the corresponding results for weak solutions to nonlocal equations defined as in \eqref{e: Main equation} in both Euclidean and Carnot group settings by taking appropriate homogeneous norms. 

Throughout the paper, we denote with $Q$  the homogeneous dimension of the group defined in \eqref{e:homogeneous dimension}. Furthermore,  the exponents $s$  and $p$ will be fixed with the assumption that  
\begin{equation}\label{e:range for p}
0<s<1 \qquad \text{and}\quad 1 < p < Q/s,
\end{equation}
unless explicitly said otherwise.  
We denote by $p^{*}(s)$ the fractional Sobolev exponent and with $p^{*}(s)'$ its H\"{o}lder conjugate,  defined by
\begin{equation}
 p^{*}(s) = \frac{pQ}{Q-ps}\quad \text{and}\quad p^*(s)'=\frac {pQ}{(p-1)Q+ps}\, .   
\end{equation} 
This definition is motivated by the following Sobolev type embedding which can be proven as in  \cite[Theorem 1.3]{Gar-Loi-Vas-2024-FR}, this time using H\"older's inequality for $p$.
\begin{customthm}{Global Sobolev Embedding} 
Let  $0 < s < 1$, $1 < p < Q/s$, and $p^{*}(s) =  \frac{pQ}{Q-ps}$.  There is a constant $S_{p,s,Q}$ such that for any $u \in {\Dsp}(\G)$ we have
\begin{equation}\label{t: Global Sobolev}
\left( \int_{\G} |u(g)|^{p^{*}(s)} dg \right)^{\frac{1}{p^{*}(s)}}
\leq S_{p,s,Q}\left( \int_{\G} \int_{\G} \frac{\abs{u(g) - u(h)}^{p}}{|h^{-1}g|^{Q+ps}} dh dg \right)^{1/p}. 
\end{equation}
\end{customthm}
The precise definitions of the relevant spaces  can be found in sub-section \ref{ss:homog groups}, and definitions \ref{d:Wsploc} and \ref{d:tail space def}. In addition, sections \ref{ss:sobolev spaces} - \ref{ss:localizations}  include foundational properties of the involved spaces, which are crucial for the analysis of the considered non-local operators.

For $F\in L^{p^*(s)'}_{\loc}(\bG)$ a function $u  \in \Wsploc \cap \tailspace$ is called a weak solution to the equation $\D_{s,p}u= F$
 if the following identity holds true for all $\phi \in C_{0}^{\infty}(\G)$,
\begin{equation}\label{e: Main equation}
D(u,\phi) =D_{s,p}(u,\phi) \overset{def}{=}\int_{\G} \int_{\G} \frac{ \tpwr{u(g) - u(h)} \paren{\phi(g) - \phi(h)}}{|h^{-1}g|^{Q+ps}} dgdh = \int\limits_{\G}F(g) \phi(g)\, dg,
\end{equation}
 where $|h^{-1}g|$ denotes the left-invariant distance arising from a fixed homogeneous norm on the group, $C_{0}^{\infty}(\G)$ denotes the space of smooth functions with compact support and $t^\fep$ denotes the function on the real line defined by 
 \begin{equation}\label{e:signed power}
  t\mapsto t^\fep=|t|^{p-2}t.  
 \end{equation}
 In {Section \ref{a:A}} we include some useful inequalities satisfied by the signed power function.  A sub-solution (resp. super-solution)  is defined in the usual manner by requiring
\begin{equation}\label{e: Main equation sub-sup-solution}
D(u,\phi) \leq\  (\ \text{resp.} \geq\ )\  \int\limits_{\G}F(g) \phi(g) dg
\end{equation}
for all non-negative $\phi \in C_{0}^{\infty}(\G)$, see also sub-section \ref{ss:weak solutions defs} for further details.  

In \cite{Gar-Loi-Vas-2024-FR},  it was shown that for $u,\, \phi\in C_0^\infty(\bG)$
we have
\begin{equation}\label{e:pv integral}
{D_{s,2}}(u,\phi)=\int_{\bG}\phi(g)\, \lap u (g) dg.
\end{equation}
Here the linear fractional operator $\lap$ is the infinitesimal generator of the Dirichlet form \eqref{e: Main equation} in the case $p=2$, i.e., 
\begin{equation}\label{e:frac L via 2nd order}
\lap u(g)=\int_{\bG} \frac {2u(g)-u(gh)-u(gh^{-1})}{|h|^{Q+2s}}dh = 2 \lim_{\eps\rightarrow
0}\int_{\bG\setminus B(g,\eps)} \frac {u(g)-u(h)}{|
h^{-1}\, g|^{Q+2s}}dh.
\end{equation}
Furthermore, by \cite[Theorem 1.2]{Gar-Loi-Vas-2024-FR} we have
\begin{equation}\label{e:limits frac operator}
\lim_{s\rightarrow 0+} \frac {s}{\sigma_Q}\lap u(g)=u(g) \quad\text{and}\quad \lim_{s\rightarrow 1-} \frac {2m(1-s)}{\tau_m}\lap u(g)=-\sum_{i=1}^m X_i^2u(g),
\end{equation}
where $\sigma_Q$ is the measure of the unit sphere (defined by the fixed norm) and, cf. \eqref{e:homog coord} for the definition of the homogeneous coordinates $x_i$, 
\begin{equation}\label{e:tau_m def}
\tau_m=\int_{|h|=1} {\sum_{i=1}^m |x_i(h)|^2}\, d\sigma(h),
\end{equation}
with $x_1, \dots, x_m$ denoting the coordinates of homogeneity one, assuming that  the fixed  homogeneous norm $|g|$ is a radial function of the coordinates $(x_1,\dots,x_m)$, i.e., its dependence on these coordinates is only as a smooth function of $x_1^2+\dots+x_m^2$.

As usual, the non-local setting requires the consideration of the so called tail term
defined for a measurable function $u$ 
and $r > 0$ via the formula
\begin{equation}\label{e: Tail}
    \Tail(u;g_{0}, r) \overset{def}{=} \left( r^{sp} \int_{\G \setminus B_{r}(g_0)} \frac{|u(g)|^{p-1}}{|g^{-1}g_{0}|^{Q+sp}} dg \right)^{\frac{1}{p-1}}.
\end{equation}
For ease of notation, we may sometimes write
$$K(g,h)=\frac{1}{|h^{-1}g|^{Q+ps}} $$ for the considered symmetric kernel.

Following are the main results of this paper. We note that  $L^{q}_{loc}(\bG)\subset L^{p^*(s)'}_{loc}(\bG)$ since $$\frac{Q}{ps}-p^*(s)'= \frac{(p-1)Q(Q-ps)}{ps\left((p-1)Q+ps \right)}. $$ 
\begin{theorem}[Local boundedness]
\label{t: Local Boundedness Interpolation Inequality}
Let  $q>Q/(sp)$ and $F\in L^{q}_{loc}(\bG)$.
 If $u $ is a sub-solution to $\D_{s,p}u= F$ then $u_{+}\in L^\infty_{\loc}(\bG)$. Furthermore there are constants  $C=C(Q,p,s)$ such that for $g_0\in\bG$  and $\delta \in (0,1]$  we have the inequality
\begin{equation}\label{e: Local Boundedness Interpolation Inequality}
\sup_{B_{R/2}} u_{+} \leq  
 \delta \Tail \left(u; g_{0}, {R}/{2}\right)
+ \delta R^{\frac {1}{(p-1)} \left(sp-\frac Qq\right) }
||F_{+}||_{L^{q}(B_{R})}^{\frac{1}{p-1}}
+  C {\delta^{-\frac{1}{\gamma}} }\left( \fint_{B_{R}} u_{+}^{p} \right)^{\frac{1}{p}},
\end{equation}
where $\gamma=\gamma(Q,s,p,q)\in (0,sp/Q)$ is defined in \eqref{e:def of gamma}.
\end{theorem}
The proof of Theorem \ref{t: Local Boundedness Interpolation Inequality} is presented in section \ref{s:proof of loc bdd}.  Versions of theorem \ref{t: Local Boundedness Interpolation Inequality} when $q=\infty$ were proven previously in \cite[Theorem 1.1]{DiC-Kuu-Pal-2016} and \cite[Theorem 1.1]{Man-Pal-Pic-2023} where similar equations on bounded domains  were considered. In the first paper the result is proven for solutions $u\in W^{s,p}(\Rn)$ on a bounded domain in the Euclidean case for vanishing $F$, see also \cite{Cozzi17} for approach relying on fractional De Giorgi classes and \cite{Behn2025DeGiorgi} for case of homogeneous systems. On the other hand,  \cite[Theorem 1.1]{Man-Pal-Pic-2023} states a result for solutions in the space $W^{s,p}(\Hn)$ on a bounded domain in the Heisenberg group $\Hn$ with a non-explicit constant  $\textbf{c}$ in \cite[(1.6)]{Man-Pal-Pic-2023}  depending on $||F||_{L^{\infty}(B_{R})}$.  The explicit form \eqref{e: Local Boundedness Interpolation Inequality} allows applications to  equations with right-hand side depending on $u$. 

The above theorem is a first step towards the next two results of the paper. In the linear case, an  approach to Harnack's inequality and H\"older continuity of solutions based on probability methods for ``jump'' processes was developed in \cite{BassLevin02}, \cite{ChenKumagai2003, ChenKumagai2008}, \cite{Bogdan_et_al2003}, \cite{Kassmann09}. In the seminal works \cite{CS} and \cite{Silvestre06} using the ``extension technique"  the non-local fractional equation in Euclidean space was transformed into a local degenerate elliptic equation in one higher dimension. This allowed the application of the classical De Giorgi-Nash-Moser regularity theory using  the regularity results of Fabes,  Kenig and Serapioni \cite{FabesKenSer82} concerning divergence form equations with $A_2$ weights. 

In the case of the fractional p-Laplacian in Euclidean space \cite{DiC-Kuu-Pal-2014} and \cite{DiC-Kuu-Pal-2016} showed the Harnack’s inequality and the local H\"older regularity of solutions  for solutions of the homogeneous fractional p-Laplacian equation in a bounded domain  generalizing the ideas of De Giorgi. In both papers  the solutions were in the space $W^{s,p}(\Rn)$. An alternative to  \cite{DiC-Kuu-Pal-2014} and \cite{DiC-Kuu-Pal-2016} approach for local solution was developed in \cite{Cozzi17} using fractional De Giorgi classes, see also \cite{cassanello_et_al25holder} for another method.  

We note that due to the generality of  our setting all known arguments involving  explicit solutions/ barriers and comparison, are not at our disposal. 

\begin{theorem}[Harnack inequality]
\label{t: Harnack}
Let  $q>Q/(sp)$, $F\in L^{q}_{loc}(\bG)$ and
  $u $ be a solution to $\D_{s,p}u= F$. If $u \geq 0 $ on $ B_R(g_0)$,
then there is a constant $C=C(s,p,Q)$
such that for all  $r<R/2^9$ we have
\begin{equation}\label{e: Harnack}
\sup_{B_{r}(g_{0})}  u \leq C \left[ \inf_{B_{r}(g_{0})} u
+  \left(\frac{r}{R}\right)^{\frac{sp}{p-1}}\Tail(u_{-};g_{0}, R)\ 
+   r^{\frac{{sp-\frac Qq}}{p-1}} \norm{F}_{L^{q}(B_{r})}^{1/(p-1)} \right].
\end{equation}
\end{theorem}
 The proof of theorem \ref{t: Harnack} will be given in section \ref{s:harnack}.  As in \cite{DiC-Kuu-Pal-2014}  we rely on the the expansion of positivity estimates. However, we consider local solutions to non-homogeneous equations and show uniform growth estimates, cf. corollary \ref{c:  intro of inf}, that are needed in the proof of the weak Harnack inequality lemma \ref{l: mean bounded by inf lemma} which relies on the Krylov-Safonov lemma.

The last result in the paper is the H\"older continuity of solutions to the equation $\frlap u =F$, cf. section \ref{ss:weak solutions defs}. A direct implementation of the  De Giorgi-Nash-Moser theory to the case of the fractional Laplacian in Euclidean space was developed in \cite{Kassmann03, Kassmann07, Kassmann09} in the linear $p=2$ case which was extended to the $p-$Laplacian case in  \cite{DiC-Kuu-Pal-2014}, see \cite{Ros-OtonSerra16} and \cite{IannizzottoMosconiSquassina2016} for the the global H\"older continuity. Regularity for viscosity solutions on bounded domains in Euclidean space was proven in \cite{Lindgren16} while  \cite{KuusiMingioneSire18} showed regularity for solutions with measure data.  

\begin{theorem}
[H\"older Continuity]
\label{t: Holder Continuity} Let  $q>Q/(sp)$.
If $u$ is a solution to the equation $\frlap u =F$ on $\bG$ with  $F\in L^{q}_{loc}(\bG)$, then $u$ is locally H\"older continuous. In fact, for any $0<R_0<1$ there  are  constants  $C= C(p,Q,s)$ 
and $\beta\in(0,1)$ such that  on any ball $B_r=B_r(g_0)$ with $r < R_{0}/2^9$ 
we have 
$$\osc_{B_{r}} u
\leq C \left(\frac {r}{R_0}\right)^\beta \left[\sup_{B_r} u + k_0\right],\quad k_0=  \left(\frac{1}{R_{0}}\right)^{\frac{sp-\frac Qq}{p-1}} \left[ C\Tail(u;g_{0}, R_{0})
    + C u(g_{0}) +  {R_0}^{\frac{sp-\frac Qq}{p-1}}\norm{F}_{L^{q}(B_{R_0})}^{1/(p-1)}\right].
$$
\end{theorem}
The proof of theorem \ref{t: Holder Continuity}, given in section \ref{s:holder cont}, follows  the classical idea of Moser where the H\"older continuity is derived from the Harnack inequality.

The above results can be applied to solutions of the fractional Yamabe equation studied in \cite{Gar-Loi-Vas-2024-OP} and \cite{Gar-Loi-Vas-2024-FR}. In fact, these two papers show that after  rescaling, there are (nonnegative) functions achieving equality in the $L^2$ Sobolev type inequality \eqref{t: Global Sobolev}, they are in the homogeneous space ${\mathcal{D}^{s,2}(\bG)}$, and  satisfy the
equation
\begin{equation}\label{e:fracYamabe1}
\frlap u=u^{2^*(s)-1},
\end{equation}
where $\frlap$ is as in \eqref{e:frac L via 2nd order}.  Furthermore, we have $u\in L^{{2^*(s)}/{2},\infty}(\bG)\cap L^\infty(\bG)$ and $|\cdot|^{{Q-2s}}\,u\in L^\infty(\bG)$. From theorem \ref{t: Harnack} and \ref{t: Holder Continuity} it follows that these solutions are positive and H\"older continuous.

\section{Definitions and background}\label{s:background}

\subsection{Homogenous groups}\label{ss:homog groups}
Throughout the paper, $\bG$ will be  a homogeneous group as defined in  \cite[Chapter 1]{FS82}. Thus, $\bG$ is a connected simply connected nilpotent Lie group. Furthermore, the exponential map $\exp:\algg\rightarrow \bG$ is a diffeomorphism between  the Lie algebra $\algg$ and the group $\bG$ and $\algg$ is endowed with a family of non-isotropic group dilations $\delta_\lambda$ for $\lambda>0$. The latter means that there is a basis $X_j$, $j=1,\dots,n$ of the Lie algebra $\algg$ and positive real numbers $d_j$, such that,
\begin{equation}\label{e:dilations on algebra}
1=d_1\leq d_2\leq\dots\leq d_n \quad \text{and}\quad
\delta_\lambda X_j=\lambda^{d_j} X_j.
\end{equation}
Using the exponential map, we define 1-parameter family of automorphisms of the group $\bG$ given by $\exp\circ\delta_\lambda\circ\exp^{-1}$. We will use the same notation $\delta_\lambda$ for the group automorphisms. As customary, $Q$ will denote the homogeneous dimension,
\begin{equation}\label{e:homogeneous dimension}
Q=d_1+\dots+d_n.
\end{equation}
We will denote with $\volg$ or $dH$ a fixed Haar measure. Since $\bG$ is a nilpotent group, the fixed Haar measure is bi-invariant. Furthermore, the homogeneous dimension and the Haar measure are related by the identity  $d(\delta_\lambda g)=\lambda^Q\, dg$.

Let $X_1,\dots X_n$ be a basis of the Lie algebra $\algg$ consisting of eigenvectors of the dilations $\delta_\lambda$ with corresponding eigenvalues $\lambda^{d_1}, \dots, \lambda^{d_n}$. We will denote with $X_j$, resp. $Y_j$, $j=1,\dots, n$ the corresponding left-invariant, resp. right-invariant, vector fields defined as usual by their action on a smooth function $f$ on $\bG$ by the formulas
\begin{equation}\label{e:inv vec fields}
(X_jf)(g)=\ddt f(g\exp(tX_j))\vert_{t=0} \quad \text{and}\quad (Y_jf)(g)=\ddt f(\exp(tY_j)g)\vert_{t=0}.
\end{equation}
The translation invariance implies that each of the vector fields $X_j$ and  $Y_j$ is formally skew-symmetric \cite[p. 21]{FS82}, i.e., for $u,\, v\in \mathcal{C}^\infty_0(\bG)$ and $j=1,\dots,n$ we have
\begin{equation}\label{e:skew-adjoint v.f.}
\int_{\bG} u(g)X_jv(g)\, dg=-\int_{\bG} v(g)X_j u (g)\, dg\quad\text{and}\quad \int_{\bG} u(g)Y_jv(g)\, dg=-\int_{\bG} v(g)Y_j u (g)\, dg.
\end{equation}
The canonical coordinate functions of the first kind on $\algg$ with respect to the fixed basis $X_1,\dots X_n$ will be denoted by $x_1,\dots,x_n$, thus for $x\in\algg$ the equation $x=(x_1,\dots,x_n)$ means
\begin{equation}\label{e:homog coord}
x=x_1X_1+\dots + x_nX_n.
\end{equation}
 In the coordinates $(x_1,\dots,x_n)$ on $\algg$ we have $X_j\vert_e=Y_j\vert_e=\partial_{x_j}\vert_0$ and the dilations take the form, recalling the normalization $d_1=1$,
\[
\delta_\lambda(x)=(\lambda x_1,\lambda^{d_2}x_2\dots,\lambda^{d_n}x_n).
\]
We will refer to these coordinates also as the homogeneous coordinates. They will be used throughout the paper.

We will denote by $| \cdot |$  a fixed homogeneous norm, i.e., $g\mapsto |g|$ is a continuous function on $\bG$ which is $\mathcal{C}^\infty$ smooth on $\bG\setminus \{e\}$, where $e$ is the group identity, $ |g|=0$ if and only if $g=e$,  and  for all $g\in\bG$ we have
\begin{equation}\label{e:homog norm}
(i) \ |g^{-1}|=|g|; \quad (ii)\ |\delta_\lambda g|=\lambda|g|.
\end{equation}
For the analysis in this paper, it is enough to assume that the norm is continuous rather than $\mathcal{C}^\infty$ outside the identity of the group.

Finally, we shall assume that the fixed norm satisfies the triangle inequality
\begin{equation}\label{e:triangle ineqs}
 |g\cdot h|\leq |g|+|h|, \qquad g,h\in \bG.
\end{equation}
According to \cite{HS90}, any homogeneous group allows a homogeneous norm which satisfies the triangle inequality.

We shall denote with $B_R(g)$ or $B(g,R)$ the open balls with center $g$ and radius $R$ defined using the fixed  homogeneous norm,
\begin{equation}\label{e:right-invarinat balls}
B_R(g)\equiv B(g,R)=\{h\, : |g^{-1}\, h|<R\}.
\end{equation}

For the balls centered at the identity or a fixed point $g_0$, when it has been specified explicitly, we will denote by $B_R$ the ball $B_R(e)$ or $B_R(g_0)$.
Notice that the distance $d(g,h)=|g^{-1}\cdot h|$ is invariant under left translations $L_{g_0}$, $L_{g_0}g=g_0\, g$, due to the equality
\begin{equation}\label{e:right inv dusatnce}
d(L_{g_0}g,L_{g_0}h)=|( g_0\, g)^{-1}\, ( g_0\, h)|=|g^{-1} g_0^{-1}\, g_0\,  h|=|g^{-1}\, h|=d(g,h),
\end{equation}
but the distance $d(g,h)$ is not necessarily invariant under right translations. In particular, the balls \eqref{e:right-invarinat balls} are obtained by left-translations from the balls centered at the identity.

From \cite[Prop. (1.15)]{FS82} the polar coordinates formula for the Haar measure gives the existence of a unique Radon measure $d\sigma(g)$ such that for $u\in L^1(\bG)$  we have the  identities
\begin{equation}\label{e:polar coordinates integral}
\int_{\bG} u(g)dg=\int_0^\infty\int_{\{|g|=1\}} u(\delta_r g) r^{Q-1}d\sigma(g)dr
\end{equation}
and,  see also \cite{FR66},
\begin{equation}\label{e:polar coord}
\int_{\{r<|g|<R\}}|g|^{-\gamma}\volg=\left\{
                  \begin{array}{ll}
                    \frac {\sigma_Q}{Q-\gamma}\left(R^{Q-\gamma} - r^{Q-\gamma} \right), & \gamma\not=Q\\
                    &\\
                    \sigma_Q\ln (R/r), & \gamma=Q,
                  \end{array}
                \right.
\end{equation}
where  $\sigma_Q = Q\omega_Q$ denotes  the measure of the unit sphere with   $\omega_Q = \int_{B_1} dg$ denoting the volume of the unit ball.

From \cite[Proposition (1.26)]{FS82} any derivative with respect to  the exponential coordinates $x_j$, $j=1,\dots,n$, can be expressed as a linear combination with polynomial coefficients of the left- (or right-) invariant vector fields. Conversely, any  left- (or right-) invariant vector field can be written as a linear combination with polynomial coefficients of the  derivative $\frac {\partial}{\partial_{x_j}}$,  $j=1,\dots,n$,
\begin{equation}\label{e:fields are polynomial}
    X_j=\sum_{jk}^n P_{jk}(x)\partial_{x_k}, \quad Y_j=\sum_{jk}^n Q_{jk}(x)\partial_{x_k},
\end{equation}
where $P_{jj}=Q_{jj}=1$, $P_{jk}=Q_{jk}=0$ if $d_k\leq d_j$ with $j\not= k $, and $P_{jk}$ and $Q_{jk}$ are polynomials of degree $d_k-d_j$ if $d_k>d_j$. Thus for any bounded domain $\Om$ and $C^1$ function $u$  the following pseudo-norms are equivalent
\begin{equation}\label{e:equiv of C and Gamma norms}
   \norm {\nabla u}_{C^0(\Om)} \overset{def}{=} \sum_{j=1}^n  \sup_{h\in \Om} |\partial_{x_j} u(h)|\quad
   \text{and}\quad  \norm {X u}_{\Gamma^0(\Om)}\overset{def}{=}  \sum_{j=1}^n  \sup_{h\in \Om } |X_ju(h)|.
\end{equation}
For $k\in \mathbb{N}$ and a multi-index $I=(i_1,i_2,\dots,i_k)$ with $i_1, i_2, \dots,i_k\in \{1,2,\dots,n\}$ we let
\[
X_I=X_{i_1}\cdot\cdot\cdot X_{i_k} \quad\text{and}\quad x_I(g)=x_{i_1}(g)\cdots x_{i_k}(g).
\]
We denote with $|I|=k$ the order of the differential operator $X_I$, while
\[
d(I)=d_{i_1}+\dots+d_{i_k}
\]
will stand for the homogeneous degree or the degree of homogeneity of the operator with respect to the fixed dilations. 

The finiteness of many of the hyper-singular integrals we shall encounter relies on  the mean value theorem \cite[(1.33)]{FS82}. In fact, the proof of \cite[(1.33)]{FS82} shows the existence of constants $T_0>0$ and $\beta\geq 1$ such that for $u\in {C}^1(B_R(e))$ and $g\in B_R(e)$ with $ B_{\beta \rho}(g)\Subset B_R(e)$ we have for $|h|<\rho$ the inequalities
\begin{equation}\label{e:MVT}
    \vert u(gh)-u(g)\vert  \leq T_0\sum_{j=1}^n |h|^{d_j}\sup_{h'\in B_{\beta \rho} (g)}|X_ju(h')|,
\end{equation}
noting explicitly that the $\sup$ involves derivatives along \emph{all} vector fields, where $X_j$ denote the left-invariant homogeneous basis of vector fields, and $\tilde B_{\beta \rho}(g) =R_g B_{\beta \rho}(e) $ denotes the right-invariant ball. {{In particular, if ${u}\in C_{0}^{1}(B_{R})$, then from \eqref{e:MVT}  we have for $g, \, h\in \bG$ the inequality
\begin{equation}\label{e:MVT cpct supp}
|{u}(gh)-{u}(g)| \leq T_0 \norm {X {u}}_{\Gamma^0(B_R)}(R^{d_n - 1}+1)|h|,
\end{equation}
using $|h|^{d_j}\leq |h|$ for $|h|\leq 1$, and  $|h|^{d_j}\leq R^{d_n -1}|h|$ for $1<|h|<R$. On the other hand, if $u\in C^1(B_R)$ with  $\norm {X {u}}_{\Gamma^0(B_R)}<\infty$, for example if ${u}\in C^{1}(\bar B_R)$, then
taking into account \eqref{e:equiv of C and Gamma norms}, it follows that \eqref{e:MVT cpct supp} will hold for $g, \, h\in B_R$.  }}

We also have the following  Taylor formula centered at $g$  of order two \cite[(1.37)]{FS82}, 
\[
u(gh)=P_{2}(h;g)+ R_{2}(h;g),
\]
where
\begin{equation}\label{e:Taylor pol ord 2}
P_{2}(h;g)=u(g)+\sum_{i:\, d_i\leq 2} (X_iu)(g)\,x_i(h)\ +\  \frac 12\sum_{i,j:\, d_i=d_j=1} (X_iX_ju)(g)\,x_i(h)x_j(h)
\end{equation}
and the remainder $R_{2}(h;g)$ satisfies the estimate
\begin{equation}\label{e:Taylor rem ord 2}
|R_{2}(h;g)|\leq \sum_{k=1}^{3}\frac {c^k}{k!}\sum_{I:\, |I|=k, d(I)> 2} |h|^{d(I)} \, \sup_{|h'|\leq \beta^3|h|} |X_Iu(gh')|, \quad g,\, h\in \bG.
\end{equation}
Thus,  using that $|h^{-1}|=|h| $ and $x_i(h^{-1})=-x_i(h)$, from \eqref{e:Taylor pol ord 2} and \eqref{e:Taylor rem ord 2} we have 
\begin{equation}\label{e:2nd order difference}
u(gh)+u(gh^{-1})-2u(g)=  \sum_{i,j:\, d_i=d_j=1} (X_iX_ju)(g)\,x_i(h)x_j(h) + \mathcal{R}(h;g),
\end{equation}
where  for some $d>2$, equal to the smallest $d(I)>2$ appearing in \eqref{e:Taylor rem ord 2}, and  for  all for $u\in C^\infty_0(B_R)$ we have the following  estimate 
\begin{equation}\label{e:2nd order difference remainder}
\vert \mathcal{R}(h;g) \vert \leq 2 T_2 M_2\,   (R^{d-1}+1)|h|^{d}, \quad h\in \bG,
\end{equation}
where  
\[
M_2=\sup \{|X_I u (h')| \,:\,  h'\in\bG, 1\leq |I|\leq 3, d(I)> 2 \}.
\]

\subsection{Sobolev spaces}\label{ss:sobolev spaces}
For an open set $\Om$ and a measurable function $u$ we define the semi-norm
\begin{equation}\label{e:sp seminorm def}
 [u]_{s,p,\Om}=\left(\int_\Om\int_\Om \frac{\abs{u(g) - u(h)}^{p}}{|h^{-1}g|^{Q+ps}} dgdh \right)^{1/p}.
\end{equation}
In the case $\Om=\bG$ we will often write simply $ [u]_{s,p}$ instead of $ [u]_{s,p,\bG}$.  As usual, the notation $U\Subset \Om$ will indicate two open sets with the closure $\bar U$ of $U$  a compact contained in $\Om$. 
\begin{definition}\label{d:Wsploc}
Let $\Om$ be an open set in $\bG$, not necessarily bounded.

a) The space  $W^{s,p}(\Om)$ is defined by
$$W^{s,p}(\Om) =\{u\in L^p(\Om):  [u]_{s,p,\Om}<\infty \}$$
equipped with the norm
\[
\norm{u}_{W^{s,p}(\Om)}=\left[ \norm{u}^p_{L^p(\Om)} +  [u]^p_{s,p,\Om}\right]^{1/p}.
\]

b) The space  $\WsplocOm$ is defined by
$$\WsplocOm=\{u : \text{ for all } U\Subset \Om, u\in L^p(U) \text{ and }  [u]_{s,p,U}<\infty   \}.$$

c) The space $W_0^{s,p}(\Om)$ is the closure of $\mathcal{C}^\infty_0(\Om)$ with respect to the $\norm{\ \cdot\  }_{W^{s,p}(\Om)}$ norm.

d) The space $\tilde{W}_{0}^{s,p}(\Om) = \{u\in \WspOm : \ \supp u \Subset \Om \}$.

\end{definition}
For an open set $\Om\subset\bG$ we define the following quantity for measurable functions $u$ and $\phi$, cf. \eqref{e:french power} for the notation $t^\fep$,
\begin{equation}\label{e:local forms}
D_{\Om}(u,\phi)=\int_{\Om}\int_{\Om}\frac {\tpwr {{u ( g ) - u(h)}}\paren{ \phi (g) - \phi (h)}}{|h^{-1}g|^{Q+ps}}  \,dg dh.
\end{equation}
In the case $\Om=\bG$ we will write  $D(u,\phi)$ for  $D_\bG(u,\phi)$. From H\"older's inequality we have
\begin{equation}\label{e:Holder for p-functional}
   | D_{\Om}(u,\phi)|\leq [u]_{s,p,\Om}^{p-1} \ [\phi]_{s,p,\Om}.
\end{equation}
When $\Om$ is a ball, $\Om=B_R$, instead of \eqref{e:local forms} we will use the notation
\begin{equation}\label{e:local forms balls}
\begin{aligned}
  D_R(u,\phi) & \overset{def}{=}D_{B_R}(u,\phi),\\    I_R(u,\phi)  & \overset{def}{=} \int_{G\setminus B_R}\int_{ B_R}\frac {\tpwr {{u ( g ) - u(h)}}\paren{ \phi (g) - \phi (h)}}{|h^{-1}g|^{Q+ps}}  \,dg dh,\\
  D(u,\phi)&  =  D_R(u,\phi)+ 2I_R(u,\phi).
\end{aligned}
\end{equation}
Before we continue, let us show that the semi-norms, hence the integrals appearing above, are finite for smooth functions with compact support. From \eqref{e:Holder for p-functional} it suffices to show that $[u]_{s,p,\Om}<\infty$ when $u \in C_{0}^{\infty}(\Om)$. First we consider the case of the local spaces. We note that in the next proposition we do not assume that $sp<Q$.

\begin{proposition}\label{p:C1 in Wspsloc}
Let $1\leq p<\infty$ and $0<s<1$. Let $\Om$ be an open set in $\bG$, not necessarily bounded.  If  $u\in  \mathcal{C}^1_{loc}(\Om)$, then $u\in \WsplocOm$,  and we have the continuous embedding $\mathcal{C}^1_{loc}(\Om)\hookrightarrow \WsplocOm$.
\end{proposition}

\begin{proof}
a) We begin with the model case  of the ball of radius $R$ centered at the identity, $\Om=B_R$, and $U=B_{r}\Subset B_{cr}\Subset B_R$, $0<r<cr<R$.
To estimate $[u]_{s,p,U}$ when $u\in\mathcal{C}^1(B_R)$, for $g\in B_{r}$ we define
 $$\delta= \frac 1\beta\min \left\{\frac {(c-1)r}{2},  1 \right\}$$
 where $\beta\geq 1$ is the constant appearing in the mean-value inequality \eqref{e:MVT}.
 From the triangle inequality  we have $$B_{\delta}(g)\Subset B_{\beta\delta}(g)\Subset B_{cr}\Subset B_R.$$ Thus,  noting that $\delta\leq 1$, for $|h|<\delta$ we have from \eqref{e:MVT} the inequalities
 \[
 \vert u(gh)-u(g)\vert  \leq T_0 \norm {X u}_{\Gamma^0(B_{\beta\delta} (g))} |h|\leq T_0 \norm {X u}_{\Gamma^0(B_r)} |h|.
 \]
Therefore, after performing a translation    we obtain from \eqref{e:MVT}, using the notation from \eqref{e:equiv of C and Gamma norms}, the inequalities:
\begin{multline}\label{e:loc frac norm C1 sing part}
  \int_{B_{r}}\int_{B_{r}\cap B_{\delta}(g)} \frac {|u(g)-u(h)|^p}{|g^{-1}h|^{Q+ps}}\,dh dg \leq \int_{B_{r}}\int_{B_{\delta}(e)} \frac {|u(gh)-u(g)|^p}{|h|^{Q+ps}}\,dh dg  \\
  \leq T_0^p  \norm {X u}_{\Gamma^0(B_{r})}^p
  \int_{ B_{r}}\int_{B_{\delta}(e)} \frac {|h|^p}{|h|^{Q+ps}}\,dh dg
  \leq \frac{\sigma_Q}{(1-s)p} T_0^p \norm {X u}_{\Gamma^0(B_{cr})}^p |B_r|\delta^{(1-s)p}\\
  = \frac{\sigma_Q^2\delta^{(1-s)p}}{Q(1-s)p}T_0^p  \norm {X u}_{\Gamma^0(B_{cr})}^p
  r^Q,
\end{multline}
recalling \eqref{e:polar coord} for the computation of the integrals of radial functions, see lemma \ref{l: phi^p-kernel-lemma}.
On the other hand, we have
\begin{multline}\label{e:loc frac norm C1 non-sing part}
 \int_{B_{r}}\int_{B_{r}\setminus B_{\delta}(g)} \frac {|u(g)-u(h)|^p}{|g^{-1}h|^{Q+ps}}\,dh dg \leq   {\delta^{-(Q+ps)}} \int_{B_{r}}\int_{B_{r}} {|u(g)-u(h)|^p}\,dh dg\\
 \leq  2^{p-1}{\delta^{-(Q+ps)}} \int_{B_{r}}\int_{B_{r}} {|u(g)|^p+|u(h)|^p}\,dh dg \leq 2^{p}\delta^{-(Q+ps)}|B_{r}|\,\norm{u}^{p}_{L^p(B_{r})}.
\end{multline}
From \eqref{e:loc frac norm C1 sing part} and \eqref{e:loc frac norm C1 non-sing part} it follows
\begin{equation}\label{e:loc frac norm for C1}
    [u]_{s,p,B_{r}}\leq C_{\bG,p,s} r^{Q}\left[ \norm {X u}_{\Gamma^0(B_{cr})}^p  \delta^{(1-s)p}
    + \delta^{-(Q+ps)}\norm{u}^{p}_{L^p(B_{r})}\right]<\infty.
\end{equation}
The general case of a pre-compact set $U\Subset \Om$ follows by a finite cover of  $U$ by balls of radius $cr$ with any $c$ such that $0<r<cr<R$ where $R=\frac 12 d(U,\bG\setminus U)$, which leads to
\begin{equation}\label{e:loc frac norm for C1 general case}
    [u]_{s,p,U}\leq C_{\bG,p, s,U, \Om, R}\left[ \norm {X u}_{\Gamma^0(U_{R})}^p+\norm{u}^{p}_{L^p(U)}\right], \quad U_R=\{g\in\Om: d(g,U)\leq R\}.
\end{equation}

\end{proof}

Following the classical setting, we define the ``tail-space'' as follows.
\begin{definition}\label{d:tail space def}
Let $1< p<\infty$ and $0<s<1$. For $\Om\subset\bG$ we let
 \[
 \mathcal{L}_s^{p-1}(\Om)\overset{def}{=}\left\{ u\in L^1_{loc}(\Om) : \int_{\Om} \frac {|u(h)|^{p-1}}{1+|h|^{Q+ps}} \, dh  < \infty  \right\}
 \]
 with a corresponding norm $$||u||_{\mathcal{L}_{s}^{p-1}(\Omega)} = \left(  \int_{\Om} \frac {|u(h)|^{p-1}}{1+|h|^{Q+ps}} \, dh \right)^{\frac{1}{p-1}}. $$
\end{definition}

From \eqref{e:polar coord} the space of continuous functions with compact support $C_{0}(\bG)$  is embedded in $ \tailspace$. Since for $|g|<R$ we have
$$\frac {1}{1+|g|^{Q+sp}}\leq 1\leq \frac {1+R^{Q+sp}}{1+|g|^{Q+sp}},$$
 for a bounded domain $\Om$ the space $\mathcal{L}_s^{p-1}(\Om)$ is the same as the Lebesgue space
$ L^{p-1}(\Om)$. Furthermore, every constant function is in the tail space. In fact an application of lemma \ref{l:Om complement kernel estimate} gives that the tail of a constant equals the constant, cf.  \eqref{e: Tail}, 
\begin{equation}\label{e:tail of constant}
   \Tail(m; g_{0}, R_{0}) = \left(R_{0}^{sp} \int_{\G \setminus B_{R_{0}}} \frac{|m|^{p-1}}{|g^{-1}g_{0}|^{Q+sp}} \right)^{\frac{1}{p-1}} = m \left( R_{0}^{sp} R_{0}^{-sp} \right)^{1/(p-1)} = m.
\end{equation}

\begin{prop}\label{p:defs of frac operator}
Let $1< p<\infty$ and $0<s<1$. Let $\Om$ be an open set in $\bG$, not necessarily bounded. If $u\in\WsplocOm \cap \tailspace$, then $D(u,\phi)$ is finite for $\phi\in C_{0}^{1}(\Om)$.
Furthermore, we have
\begin{equation}\label{e:estimate in defs of frac operator general case}
|D(u,\phi)|\leq  [u]_{s,p,K_{R}}^{p-1}\, [\phi]_{s,p,K_{R}}\\
+ C\left[\norm{u}_{\mathcal{L}_s^{p-1}(\bG)}^{p-1} + \norm{u}^{p-1}_{L^p(K_{R})}\right]\norm{\phi}_{L^p(\Om)}
\end{equation}
with a constant $C=C(Q,s,p,K,R)$, where $K = \supp(\phi)$, $R=d(K,\bG\setminus \Om)/2$ and $K_{R} = \{g \in K : d(g,\partial \Omega) < R \}$. In particular, $D(u,\phi)$ is well defined continuous functional for $u\in\WsplocOm \cap \tailspace$ and $\phi\in\WsplocOm \cap L^p(\Om) $.
\end{prop}

\begin{proof}
 We will consider first the model case  of the ball of radius $R$ centered at the identity, $\Om=B_R$, and $\phi\in C_{0}^{\infty}(B_{r})$ with $0<r<R$. Let $c>0$ be such that  $0<r<cr<R$.  Since $\phi$ vanishes outside the ball $B_{cr}$ we have the decomposition, cf. \eqref{e:local forms balls},
\begin{equation}\label{e:global form finite 1}
  D(u,\phi) = D_R(u,\phi)  + 2I_R(u,\phi).
\end{equation}
We shall estimate each of the integrals $D_R(u,\phi)$ and $ I_R(u,\phi)$ separately.

First, taking into account that $\supp \phi \subset B_{r}$
and H\"older's inequality we have
\begin{multline}\label{e:local forms finite 1}
|D_R(u,\phi)|\leq |D_{cr}(u,\phi)| +2J_2 \leq [u]_{s,p,B_{cr}}^{p-1}\, [\phi]_{s,p,B_{cr}} + 2J_2, \\\text{ where }\\
J_2=\int_{B_R\setminus B_{cr}}\int_{B_{cr}}   \frac{|u \pl g \pr - u(h)|^{p-1} \abs{\phi (h)}}{|g^{-1}h|^{Q+sp}} \,dh dg=\int_{B_R\setminus B_{cr}}\int_{B_{r}}\frac {|u(g)-u(h)|^{p-1} \abs{\phi (h)} }{|g^{-1}h|^{Q+sp}}\,dh dg.
\end{multline}
 Note that the semi-norm $ [u]_{s,p,B_{cr}}$ is finite since we have assumed $u\in\WsplocOm$ while $[\phi]_{s,p,B_{cr}}$ is finite by proposition \ref{p:C1 in Wspsloc}.
 For $g\in{B_R\setminus B_{cr}}$ and $h\in {B_{r}}$ as in the domains of the integral in $J_2$ we have
\[
 |g^{-1}h|\geq |g|-|h|\geq \frac {c-1}{c}|g|.
\]
This, together with $|a-b|^{p-1}\leq 2^{p-1}\left( |a|^{p-1}+ |b|^{p-1}\right)$ and $\supp \phi\, \subset B_r$,
give  
\begin{multline}\label{e:J2}
J_2\leq  2^{p-1}\left( \frac {c}{c-1}  \right)^{Q+sp}\left[\int_{B_R\setminus B_{r}}\int_{B_{r}} \abs{\phi (h)}  \frac{|u(g)|^{p-1}+|u(h)|^{p-1}}{|g|^{Q+sp}}  \,dh dg \right]\\
\leq 2^{p-1}\left( \frac {c}{c-1}  \right)^{Q+sp}\left[\int_{B_R\setminus B_{r}}\int_{B_{r}} \abs{\phi (h)}  \frac{|u(g)|^{p-1}}{|g|^{Q+sp}}  \,dh dg + \int_{B_R\setminus B_{r}}\int_{B_{r}} \abs{\phi (h)}  \frac{|u(h)|^{p-1}}{|g|^{Q+sp}}  \,dh dg \right].
\end{multline}
From H\"older's inequality we obtain the estimates
\begin{equation}
    \int_{B_R\setminus B_{r}}\int_{B_{r}} \abs{\phi (h)}  \frac{|u(g)|^{p-1}}{|g|^{Q+sp}}  \,dh dg \leq |B_{r}|^{1/p'}\norm{\phi}_{L^p(B_{r})}\ \int_{B_R\setminus B_{r}}   \frac{|u(g)|^{p-1}}{|g|^{Q+sp}}\, dg
\end{equation}
and
\begin{multline}
 \int_{B_R\setminus B_{r}}\int_{B_{r}} \abs{\phi (h)}  \frac{|u(h)|^{p-1}}{|g|^{Q+sp}}  \,dh dg \leq \int_{B_R\setminus B_{r}} \frac{1}{|g|^{Q+sp}}  \,dg\, \int_{B_{r}} \abs{\phi (h)} |u(h)|^{p-1}\,dh \\
 \leq \frac{\sigma_Q}{sp} r^{-sp}\norm{\phi}_{L^p(B_{r})}\norm{u}^{p-1}_{L^p(B_{r})}.
\end{multline}
With these estimates we have
\begin{multline}\label{e:J2alternative}
J_2\leq  2^{p-1}\left( \frac {c}{c-1}  \right)^{Q+sp}\left[ |B_{r}|^{1/p'} \int_{B_R\setminus B_{r}}   \frac{|u(g)|^{p-1}}{|g|^{Q+sp}}\, dg + \frac{\sigma_Q}{sp} r^{-sp}\norm{u}^{p-1}_{L^p(B_{r})}\right]\norm{\phi}_{L^p(B_{r})}\\
\leq 2^{p-1}\left( \frac {c}{c-1}  \right)^{Q+sp}\Bigg[ |B_{r}|^{1/p'}\frac {1+r^{Q+sp}}{r^{Q+sp}}  \norm{u}_{\mathcal{L}_s^{p-1}(B_R\setminus B_{r})}^{p-1}
 + \ \frac{\sigma_Q}{sp}r^{-sp}\norm{u}^{p-1}_{L^p(B_{r})}\Bigg]\norm{\phi}_{L^p(B_{r})}\\
=2^{p-1}\left( \frac {c}{c-1}  \right)^{Q+sp}\Bigg[ \omega_Q^{(p-1)/p}\frac {1+r^{Q+sp}}{r^{Q/p+sp}}  \norm{u}_{\mathcal{L}_s^{p-1}(B_R\setminus B_{r})}^{p-1} \  + \ \frac{\sigma_Q}{sp} r^{-sp}\norm{u}^{p-1}_{L^p(B_{r})}\Bigg]\norm{\phi}_{L^p(B_{r})},
\end{multline}
using \eqref{e:polar coord} and that for
$|g|\geq r$ we have
\begin{equation}\label{e:switch to trace norm}
\frac {1}{|g|^{Q+sp}}\leq \frac {1+r^{Q+sp}}{r^{Q+sp}}\frac {1}{1+|g|^{Q+sp}}.
\end{equation}
Thus, from \eqref{e:J2alternative} and
 \eqref{e:local forms finite 1}, the local form $D_R$ can be estimated as follows,
\begin{multline}\label{e:local form estimate alternative}
|D_R(u,\phi)|  \leq [u]_{s,p,B_{cr}}^{p-1}\, [\phi]_{s,p,B_{cr}} + 2^{p}\left( \frac {c}{c-1}  \right)^{Q+sp} \left[C_1  \norm{u}_{\mathcal{L}_s^{p-1}(B_R\setminus B_r)}^{p-1} +C_2\norm{u}^{p-1}_{L^p(B_{r})}\right]\norm{\phi}_{L^p(B_{r})},\\
\text{where}\quad C_1=\omega_Q^{(p-1)/p}\frac {1+r^{Q+sp}}{r^{Q/p+sp}}\quad\text{and} \quad C_2 = \frac{\sigma_Q}{sp}r^{-sp}.\hfill.
\end{multline}

The estimate for $I_R(u,\phi)$ in \eqref{e:global form finite 1} is very similar to the estimate of the integral $J_2$ above. It can be obtained from \eqref{e:J2alternative} by replacing $B_R$ with $\bG$ and $B_{r}$ with $B_R$. In other words we can let $R\rightarrow\infty$ in \eqref{e:J2} and then replace $r$ with $R$ to obtain
from \eqref{e:J2alternative} the inequality
\begin{multline}\label{e:I_R estimate alternative}
I_R(u,\phi)    \leq 2^{p-1}\left( \frac {c}{c-1}  \right)^{Q+sp}\Bigg[ \omega_Q^{(p-1)/p}\frac {1+R^{Q+sp}}{R^{Q/p+sp}}  \norm{u}_{\mathcal{L}_s^{p-1}(\bG\setminus B_{R})}^{p-1}\\
+ \frac{\sigma_Q}{sp}R^{-sp}\norm{u}^{p-1}_{L^p(B_{R})}\Bigg]\norm{\phi}_{L^p(B_{R})}.
\end{multline}

Finally, $|D(u,\phi)|$ can be estimated with the help of \eqref{e:global form finite 1} and the above estimates for $D_R(u,\phi)$ and $I_R(u,\phi)$, which give, recalling that $\phi\in C_{0}^{\infty}(B_{r})$ with  $0<r<cr<R$,
\begin{multline}\label{e:estimate in defs of frac operator}
|D(u,\phi)|\leq  [u]_{s,p,B_{cr}}^{p-1}\, [\phi]_{s,p,B_{cr}}\\
+ 2^{p}\left( \frac {c}{c-1}  \right)^{Q+sp} \left[\omega_Q^{(p-1)/p}\frac {1+r^{Q+sp}}{r^{Q/p+sp}}\norm{u}_{\mathcal{L}_s^{p-1}(B_R)}^{p-1} + \frac{\sigma_Q}{sp}r^{-sp} \norm{u}^{p-1}_{L^p(B_{r})}\right]\norm{\phi}_{L^p(B_{r})}\\
+ 2^{p}\left( \frac {c}{c-1}  \right)^{Q+sp}\Bigg[ \omega_Q^{(p-1)/p}\frac {1+R^{Q+sp}}{R^{Q/p+sp}}  \norm{u}_{\mathcal{L}_s^{p-1}(\bG\setminus B_{R})}^{p-1}
+ \frac{\sigma_Q}{sp}R^{-sp}\norm{u}^{p-1}_{L^p(B_{R})}\Bigg]\norm{\phi}_{L^p(B_{R})}.
\end{multline}

The case of a general domain $\Om$ can be obtained from the model case of a ball using a partition of unity subordinate to a finite cover of $K=\supp \phi \Subset\Om$ by balls contained in the bounded domain $$K_R=\{g\in K : d(g,\partial\Om)<R\},$$ where $R=d(K,\bG\setminus \Om)/2$ taking into account that $|D(u,\phi)| $ is sub-linear in $\phi$, i.e.,
\[
|D(u,\phi_1+\phi_2)|\leq |D(u,\phi_1)|+ |D(u,\phi_2)|.
\]
Thus we obtain \eqref{e:estimate in defs of frac operator general case}, which completes the proof.
\end{proof}

\begin{corollary}\label{c:C1 in Wspsloc}
Let $1< p<\infty$ and $0<s<1$. Let $\Om$ be an open set in $\bG$, not necessarily bounded.
  If  $u\in \tailspace \cap \mathcal{C}^1(\Om)$ and  $\phi\in \mathcal{C}_{0}^{1}(\Om)$, then $D(u,\phi)$ is finite.
\end{corollary}

\begin{proof} From proposition  \ref{p:C1 in Wspsloc},  if  $u\in  \mathcal{C}^1(\Om)$ then $u\in \WsplocOm$. Then, proposition \ref{p:defs of frac operator} shows that if $u\in\WsplocOm \cap \tailspace$ and $\phi\in C_{0}^{1}(\Om)$ then $D(u,\phi)$ is finite.
\end{proof}

\subsubsection{Homogeneous Sobolev spaces}
In view of the above facts, the following definitions determine non-trivial spaces.

\begin{definition}\label{d: dsp-Om} Let $1< p<\infty$ and $0<s<1$.
For an open set $\Om \subseteq \bG$ we define the $L^p$ Dirichlet space
${\Dsp}(\Om)$ as the completion of $C_{0}^{\infty}(\Om)$ with respect to the \textit{global} norm

\begin{equation}\label{d: dsp-Br}
    \norm{u}_{\Dsp(\Om)}=D(u,u)^{1/p}= \left( \int_{\G} \int_{\G} \frac{\abs{u(g) - u(h)}^{p}}{|h^{-1}g|^{Q+ps}} dg dh \right)^{\frac{1}{p}}.
\end{equation}
\end{definition}
Clearly, ${\Dsp}(\Om)\hookrightarrow {\Dsp}(\bG)$.  From H\"older's inequality we have, cf. \eqref{e:local forms},
\begin{equation}\label{e:Holder on Dirichlet space}
|D_\Om(u,v)|\leq  \norm{u}^{p-1}_{\Dsp(\Om)}  \norm{v}_{\Dsp(\Om)}, \quad u,v \in {\Dsp(\Om)}.
\end{equation}
When $sp<Q$, i.e., \eqref{e:range for p} holds true, the Sobolev type embedding theorem \eqref{t: Global Sobolev} and H\"{o}lder's inequality, noting that $p^*(s)>p$, show that  we have the continuous embedding
\begin{equation}\label{e:Dsp in Wsploc}
\Dsp(\Om)\hookrightarrow \, \WsplocOm.
\end{equation}
Furthermore, the space ${\Dsp}(\bG)$ can be identified with the space of functions
\begin{equation}\label{e:alternative Dsp def)}
 \Dsp(\G)=\{ u\in L^{p^*(s)}(\bG): [u]_{s,p,\bG}<\infty \}.
\end{equation}
The proof of this claim relies on truncation and approximations of the identity.  Details in the case $p=2$ can be found in \cite{Gar-Loi-Vas-2024-OP}, see also  \cite[Theorem 2.2]{del2018spectrum} or \cite[Theorem 3.1]{BG-CV21}, where it is shown that in the Euclidean setting a function in $L^{p^*(s)}(\R^n)$ with a finite $[\, .\, ]_{s,p}$ semi-norm can be approximated in ${\Dsp}(\bG)$  by a sequence in $\COG$.  The generalization to the case of general $p$ considered here is straightforward taking into account the regularization via convolution results in sub-section \ref{ss:convolutions}. Furthermore, similar arguments show that the space $\COG$ is dense in $\WspG$. In particular, we have the continuous embedding
\begin{equation}\label{e:Wsp in Dsp}
    \WspG \hookrightarrow \DspG.
\end{equation}

The global homogeneous norm can be put in the form
\begin{equation}\label{e:def frac norm}
\norm{u}_{\mathcal{D}^{s,p}(\bG)}=\left(\int_{\bG}\int_{\bG}\frac {|u(g)-u(h)|^p}{| h^{-1}\, g|^{Q+ps}}\, {dg}dh\right)^{1/p}= \left(\int_{\bG}\int_{\bG}\frac {|u(gh)-u(g)|^p}{|h|^{Q+ps}}\, {dg}dh\right)^{1/p}
\end{equation}
using the change of variables $h'= h^{-1}\, g$ and the fact that $|h^{-1}|=|h|$. In particular, the $\mathcal{D}^{s,p}(\bG)$ norm is invariant under left-translations since it involves the left-translation invariant distance $d(g,h)=|h^{-1}g|$ and the left-invariant Haar measure, which we recall is  bi-invariant in our case. It is worth pointing out that unlike the Abelian setting, in general, only the first formula in \eqref{e:MVT} is relevant for the integrals in \eqref{e:def frac norm}.

\subsection{Approximations with smooth functions}\label{ss:convolutions}
We will consider as in \cite[p. 15]{FS82} the ``right convolution'' on $\bG$ defined by the formula
\begin{equation}\label{e:left convolution}
u*v (g)=\int_{\bG} u(gh^{-1}) v(h)dh= \int_{\bG} u(h) v(h^{-1}g)dh.
\end{equation}
The use of this convolution is motivated by the choice of the left-translation invariant distance, cf. \eqref{e:right inv dusatnce}. Notice that, cf. \cite[p. 22]{FS82}, we have
\begin{equation}\label{e:derivative and covol}
    X_j(u*v)=u*(X_jv), \ Y_j (u*v)=(Y_ju)*v, \quad u*(Y_jv)=(X_ju)*v.
\end{equation}
The convolution satisfies Young's inequalities, \cite[Proposition (1.18) and (1.19)]{FS82}, 
\begin{align}
& \norm {u*v}_{L^{r}(\bG)} \leq C_{p,q} \norm{u}_{L^{p}(\bG)} \, \norm{v}_{L^{q}(\bG)}, \quad 1\leq  p,\, q,\, r\leq \infty, \ p^{-1}+q^{-1}=r^{-1}+1,\\
 & \norm {u*v}_{L^{q,\infty}(\bG)} \leq C_{q} \norm{u}_{L^{1}(\bG)} \, \norm{v}_{L^{q,\infty}(\bG)}, \quad 1<  q <\infty,\\
 & \norm {u*v}_{L^{r}(\bG)} \leq C_{p,q} \norm{u}_{L^{p}(\bG)} \, \norm{v}_{L^{q,\infty}(\bG)}, \quad 1< p,\, q,\, r<\infty, \ p^{-1}+q^{-1}=r^{-1}+1.
\end{align}

For  $\eps>0$ and $\rho\in C^\infty_0(\bG)$ such that $\int_{\bG} \rho dh =1$, let
\begin{equation}\label{e:delta aproxim}
 \rho_\eps(g)=\eps^{-Q}\rho(\delta_{1/\eps} g)\quad \text{ and }\quad u_\eps =\rho_\eps * u.
\end{equation}
Notice that \eqref{e:derivative and covol} shows that $u_\eps\in C^\infty(\bG)$ whenever finite. This follows similarly to the arguments in \cite{Gar-Loi-Vas-2024-OP} about the Schwartz class of functions and its interaction with the fractional operator $\frlap$. In fact, taking into account \cite[proposition (1.26)]{FS82}, any derivative with respect to  the exponential coordinates $x_j$, $j=1,\dots,n$, can be expressed as a linear combination with polynomial coefficients of the left- (or right-) invariant vector fields. Conversely, any  left- (or right-) invariant vector field can be written as a linear combination with polynomial coefficients of the  derivative $\frac {\partial}{\partial_{x_j}}$,  $j=1,\dots,n$.

By \cite[proposition (1.20)]{FS82} we have the following facts.
\begin{prop}\label{p:conv prop FS} The next limits hold true as  $\eps\rightarrow 0$.
    \begin{enumerate}[(i)]
    \item If $1\leq p<\infty$, then for $u\in L^p(\bG)$, $u_\eps\rightarrow u$ in $L^p(\bG)$.
    \item If $u$ is bounded and right uniformly continuous, then $u_\eps\rightarrow u$ in $L^\infty(\bG)$.
    \item  If $u$ is bounded on $\bG$ and continuous on an open set $\Om\subset \bG$, then $u_\eps\rightarrow u$ uniformly on compact subsets of $\Om$.
\end{enumerate}
\end{prop}
The following facts were proven in \cite{Gar-Loi-Vas-2024-OP} in the case $p=2$. We use again $u_\eps=\rho_\eps * u$ as in \eqref{e:delta aproxim} and  $R_h$ for the right translation with $(R_hu)(g)=u(gh)$.
\begin{prop}\label{l:conv and Sobolev norm} Let $1\le p<\infty$.
\begin{enumerate}[a)]
\item  If  $u\in L^p_{loc}(\bG)$  and  $ [u]_{s,p,\bG}<\infty$, then we have
\begin{equation}\label{e:Lp norm of R_hu -u}
    R_h u -u \in {L^p(\bG)} \ \text{ and }\ \norm {R_h u -u}_{L^p(\bG)}\leq C |h|^{s}[{u}]_{{{s,p},\bG}}.
\end{equation}
\item The mollification \eqref{e:delta aproxim} does not increase the $L^p$ norm and the fractional semi-norm, assuming   $u\in L^p_{loc}(\bG)$  and  $ [u]_{s,p,\bG}<\infty$.  In the latter case:
\begin{equation}\label{e:conv and Sobolev norm}
\norm{u_\eps}_{L^p(\bG)}\leq \norm{u}_{L^p(\bG)}, \quad \quad[{u_\eps}]_{{s,p},\bG}\leq [{u}]_{{{s,p},\bG}}.
\end{equation}
\item If  $u\in L^p_{loc}(\bG)$  and  $ [u]_{s,p,\bG}<\infty$, then $[u_\eps - u]_{s,p,\bG}\rightarrow 0$  as $\eps\rightarrow 0$. In particular, if $u\in\mathcal{D}^{s,p}(\bG) $ then $u_\eps \rightarrow u$ in $\mathcal{D}^{s,p}(\bG) $ while if $u\in\WspG $ then $u_\eps \rightarrow u$ in $\WspG$ as $\eps\rightarrow 0$.
\end{enumerate}
\end{prop}

\begin{proof}

a) Let $r=|h|$ and $y\in B_r(g)=\{ y\in\bG: \, |y^{-1}g| < r\}$, cf. \eqref{e:right-invarinat balls}. By the triangle inequality  we have
\[
|R_hu (g)-u(g)|=|u(gh)-u(g)|\leq |u(gh)-u(y)| + |u(y)-u(g)|.
\]
Averaging over the ball $B_r(g)$ we have
\[
|u(gh)-u(g)|\leq \fint_{B_r(g)}|u(gh)-u(y)|dy + \fint_{B_r(g)}|u(y)-u(g)|dy.
\]
Using $(a+b)^p\leq 2^{p-1}(a^p+b^p)$ for $p\geq 1$, see \eqref{e: triangle with p}, and Holder's inequality we have
\begin{equation}\label{e:triangle ineq for Lhu-u}
\frac 1{2^{p-1}} |u(gh)-u(g)|^p\leq \fint_{B_r(g)}|u(gh)-u(y)|^p \, dy + \fint_{B_r(g)}|u(y)-u(g)|^p\, dy.
\end{equation}
We set
$$A(g) \overset{def}{=} \fint_{B_r(g)}|u(gh)-u(y)|^p \, dy\mbox{ and } B(g) \overset{def}{=} \fint_{B_r(g)}|u(y)-u(g)|^p\, dy.$$
Next, we estimate separately the integrals of $A$ and $B$. First, we compute
\begin{multline}
A(g)=\frac{1}{\omega_{Q}} r^{-Q}\int_{|y^{-1}g|<r} |u(gh)-u(y)|^p\, dy = \frac{1}{\omega_{Q}} r^{-Q}\int_{|hz|<r} |u(gh)-u(ghz)|^p\, dz\\
\leq  \frac{1}{\omega_{Q}} r^{-Q}\int_{|z|<2r} |u(gh)-u(ghz)|^p\, dz
\end{multline}
by first letting $y=ghz$, hence $y^{-1}g=(hz)^{-1}$, and then noting that for $|y^{-1}g|<r$ we have from the triangle inequality
$
|z|=|h^{-1}hz|\leq |hz| +|h^{-1}|=|hz| +|h|< 2r.
$
Integrating over $\bG$, applying Fubini's theorem and using $|z|<2r$ we come to
\begin{multline}\label{e:A term}
\int_{\bG} A(g)\volg \leq \frac{1}{\omega_{Q}} r^{-Q}\int_{\bG}\int_{|z|<2r} |u(gh)-u(ghz)|^p\, dz\volg\\
\leq \frac{1}{\omega_{Q}} 2^{Q+ps} r^{ps}\int_{|z|<2r} \int_{\bG}\frac {|u(gh)-u(ghz)|^p}{|z|^{Q+ps}}\volg dz
\leq  \frac{1}{\omega_{Q}} 2^{Q+ps} r^{ps}\int_{\bG} \int_{\bG}\frac {|u(g)-u(gz)|^p}{|z|^{Q+ps}}\volg dz\\
\leq \frac{1}{\omega_{Q}} 2^{Q+ps} r^{ps}[u]_{s,p} 
\end{multline}
using the change of variable $g'=gh$ in the next to last inequality and, finally,  the definition of the $[u]_{s,p} $ semi-norm from \eqref{e:sp seminorm def}.

Turning to the estimate of the integral $B(g)$ defined in \eqref{e:triangle ineq for Lhu-u}, an integration over $\bG$ gives
\begin{multline}\label{e:B term}
\int_{\bG} B(g)\volg\leq \frac{1}{\omega_{Q}} r^{-Q}\int_{\bG}\int_{|y^{-1}g|<r} |u(g)-u(y)|^p\, dy\volg \\
\leq \frac{1}{\omega_{Q}} r^{ps}\int_{\bG}\int_{|y^{-1}g|<r} \frac {|u(g)-u(y)|^p}{r^{Q+ps}}dy\volg
\leq \frac{1}{\omega_{Q}} r^{ps}\int_{\bG}\int_{|y^{-1}g|<r} \frac {|u(g)-u(y)|^p}{|y^{-1}g|^{Q+ps}}dy\volg\\
\leq \frac{1}{\omega_{Q}} r^{ps}\int_{\bG}\int_{\bG} \frac {|u(g)-u(y)|^p}{|y^{-1}g|^{Q+ps}}dy\volg= \frac{1}{\omega_{Q}} r^{ps}[u]_{s,p} . 
\end{multline}
The claim of part a) follows from \eqref{e:triangle ineq for Lhu-u}, \eqref{e:A term} and \eqref{e:B term}.

b) Minkowski's inequality and the invariance of the measure by (left-)translations gives
\begin{multline}\label{e:right translation and mollification}
    \norm{u_\eps}_{L^p(\bG)}
    = \left[\int_\bG  \left\vert \int_\bG \rho_\eps(h) u(h^{-1}g)\, dh\right\vert^p \,dg\right]^{1/p}
    \leq \int_\bG   \left[ \int_\bG \left\vert u(h^{-1}g)\right\vert^p\,dg\right]^{1/p}\rho_\eps(h)\, dh \\
    =\int_\bG \norm{u}_{L^p(\bG)} \rho_\eps(h)\, dh =  \norm{ u}_{L^p(\bG)}.
\end{multline}
To prove the second inequality in \eqref{e:conv and Sobolev norm}, as in the Euclidean case, following \cite[(15.8)]{Po16}, we rely on Jenssen's inequality. First, we recall that the dilations $\delta_\eps$ are a one-parameter group of  automorphisms, hence,
\begin{equation}\label{e:inverses of dilation}
\delta_{1/\eps}^{-1}h=\delta_{\eps}h\quad \text{and}\quad (\delta_\eps \, h)^{-1}=\delta_\eps \, h^{-1}.
\end{equation}
Therefore,  we have
\[
u_\eps\,  (g)= \eps^{-Q}\int_\bG  \rho(\delta_{1/\eps} h) u(h^{-1}g)\, dh = \int_\bG  \rho( h) u((\delta_\eps h)^{-1}g)dh
=\int_\bG  \rho( h) u((\delta_\eps h^{-1})g)\, dh.
\]
Thus, with the help of Jenssen's inequality, the invariance by left-translations of the distance and of  the Haar measure,  the norm of the convolution can be estimated as follows,
\begin{multline*}
[u_\eps]^p_{{s,p}}\leq \int_\bG\int_\bG \rho(z)\int_\bG \frac {|u((\delta_\eps h^{-1})g) - u((\delta_\eps h^{-1})y)|^p}{|y^{-1}g|^{Q+ps}}dgdy dh\\
=\int_\bG\int_\bG \rho(h)\int_\bG \frac {|u(g) - u(y)|^p}{|y^{-1}g|^{Q+ps}}dgdy dh=\int_\bG\int_\bG \frac {|u(g) - u(y)|^p}{|y^{-1}g|^{Q+ps}}dgdy \int_\bG \rho(h) \, dh =[{u}]_{{s,p}}\ ,
\end{multline*}
which proves \eqref{e:conv and Sobolev norm}.

c) Following \cite[Theorem 6.62]{Leoni23}, let $\triangle_h u$ be the function $$\triangle_h u (g) = (R_hu-u)(g)=u(gh)-u(g).$$ From part a) for $ u\in \Dsp(\bG)$ we have $\triangle_h u\in L^p(\bG)$.  From \eqref{e:def frac norm} we have
\begin{equation}\label{e:Rh and sp norm}
||u||_{\mathcal{D}^{s,p}(\bG)} 
= \left(\int_{\bG}\int_{\bG}\frac {|\triangle_h u (g)|^p}{|h|^{Q+ps}}\, {dg}dh\right)^{1/p}= \left(\int_{\bG}\frac {\norm{\triangle_h u}_{L^p(\bG)}^p}{|h|^{Q+ps}}\, dh\right)^{1/p}<\infty. 
\end{equation}
Note that $R_hu_\eps=(R_hu)_\eps$ since
\[
    R_hu_\eps\,  (g)=u_\eps\,  (gh)= \eps^{-Q}\int_\bG  \rho(\delta_{1/\eps} z) u(z^{-1}g)dz = \eps^{-Q}\int_\bG  \rho(\delta_{1/\eps} z) (R_hu)(z^{-1}g)dz=(R_hu)_\eps\,  (g),
\]
thus we have
\begin{equation}\label{e:[triangle,Rh]}
  (\triangle_h u)_\eps=\triangle_h u_\eps.
\end{equation}
 Invoking \eqref{e:right translation and mollification} we obtain

\begin{equation}\label{e:Rh and mollification}
   \norm{\triangle_h u_\eps}_{L^p(\bG)}= \norm{(\triangle_h u)_\eps}_{L^p(\bG)}
    \leq\norm{\triangle_h u}_{L^p(\bG)}.
\end{equation}
 From \eqref{e:Rh and sp norm} we have
$$[u_\eps - u]_{s,p,\bG}=\left(\int_{\bG}\frac {\norm{\triangle_h (u_\eps - u)}_{L^p(\bG)}^p}{|h|^{Q+ps}}\, dh\right)^{1/p}. $$
By the properties of the mollification we have  $\norm{(\triangle_h u)_\eps- \triangle_h u}_{L^p(\bG)}\rightarrow 0$. 
The claim that $[u_\eps - u]_{s,p}\rightarrow 0$  as $\eps\rightarrow 0$ follows from Lebesgue's dominated convergence theorem taking into account that $\norm{\triangle_h (u_\eps - u)}_{L^p(\bG)}\leq 2 \norm{\triangle_h  u}_{L^p(\bG)}\in L^1({{|h|^{-Q-ps}}\, dh})$ by \eqref{e:Rh and sp norm} since $u\in\Dsp(\bG)$.

Finally,  \eqref{e:alternative Dsp def)} and proposition \ref{p:conv prop FS} imply the rest of the claims.
\end{proof}

\subsection{The chain rule}\label{ss:chain rule}
The usual chain rule is valid in the considered setting. Thus, if $f:(a,b)\rightarrow \R$ is Lipschitz continuous, $\Om\subset\bG$ is an open set, and $u\in{\Dsp}(\Om)$ takes values in $(a,b)$, then $f\circ u \in {\Dsp}(\Om)$. In fact, if $|f(t_1)-f(t_2)|\leq L |t_1-t_2|$, $t_1,\, t_2\in (a,b)$, then we have:
\begin{equation}\label{e:chain rule}
[{f\circ u}]_{s,p,\Om}\leq L [{u}]_{s,p,\Om}\quad \text{and}\quad
  \norm{f\circ u}_{\mathcal{D}^{s,p}(\Om)}\leq L \norm{u}_{\mathcal{D}^{s,p}(\Om)}.
\end{equation}
Furthermore, if $u\in W^{s,p}_{loc}(\Om)$ then $f\circ u\in W^{s,p}_{loc}(\Om)$; in particular, $ (u -l)_{\pm}\in W^{s,p}_{loc}(\Om)$ for any constant $l\in\R$.

\subsection{Localizations}\label{ss:localizations}
The proof of the next lemma is contained in \cite[(3.5)]{Gar-Loi-Vas-2024-OP}
\begin{lemma}[Existence of Test Functions]\label{l: existence of test functions} Suppose that $0 < r < R$.  Let $g_{0} \in \G$ and let $B_{R} = B(g_{0},R)$, $B_{r} = B(g_{0}, r)$. Then there is a function $\phi \in C_{0}^{\infty}(B_{R})$ and a constant $c$ such that $\phi \equiv 1$ on $B_{r}$ and $$|\phi(g) - \phi(h)| \leq \frac{4}{R-r} |h^{-1}g|.$$
\end{lemma}

For $\phi\in  \mathcal{C}_{0}^{1}(\Om)$ we consider the mapping properties of the multiplication operator $u\mapsto M_\phi(u)=\phi u$. Together with the chain rule, these will be used subsequently freely when testing the equation with various test functions from $\WsplocOm$. From proposition \ref{p:C1 in Wspsloc} and H\"{o}lder's inequality we have
\[
|D_\Om(u,\phi)|\leq [u]_{s,p,\Om}^{p-1}\, [\phi]_{s,p,\Om}.
\]
The next proposition implies, in particular, that the local Sobolev space can be characterized in the following manner,
$$\WsplocOm=\{u\in L^p_{loc}(\Om) : \text{ for every } \phi \in \mathcal{C}_{0}^{\infty}(\Om)\ \text{ it holds }\ \phi u\in {W^{s,p}(\Om)}  \}.$$
\begin{prop}\label{p:localization} Let  $\Om$ be an open set and $\phi\in  \mathcal{C}_{0}^{1}(\Om)$. Let $M_\phi$ denote the multiplication by $\phi$ operator.
\begin{enumerate}[a)]
    \item We have $[\phi u]_{s,p,\bG}  \leq [u]_{s,p,\bG}$ and   the map $M_\phi:\Dsp (\bG)\rightarrow \Dsp (\Om)$  is continuous.
    \item The maps $M_\phi:\WsplocOm\rightarrow \WspG$ and $M_\phi:\WsplocOm\rightarrow \DspG$ are continuous.
\end{enumerate}
In particular, $M_\phi:\Dsp (\Om)\rightarrow \WsplocOm$,   $M_\phi:\Dsp (\Om) \rightarrow \Dsp(\bG)$,  and $M_\phi:\WsplocOm\rightarrow \WspOm$ are continuous embeddings.
\end{prop}

\begin{proof} Let $K=\supp \, \phi \Subset \Om$ and $v=\phi u$ 
with   $\supp\, v\subset K $. Let us fix an additional open set $U$ such that $K\Subset U\Subset \Om$.

a) 
  We will estimate the semi-norm $[{v }]_{s,p,\bG}$. 
  Decomposing the integral on $\bG\times\bG$ as the sum of the integrals on $U\times U$ and $(\bG\setminus U) \times U$, and invoking the polar coordinates to estimate  the integral on $(\bG\setminus U)\times U$, cf. lemma \ref{l:Om complement kernel estimate}, we obtain
\begin{multline}\label{e:Mphi to Dsp}
[{v }]^p_{s,p}  \leq 
 [\phi u]^p_{s,p,U} + 2\int_K |u\phi(h)|^p\int_{\bG\setminus U} \frac {1}{|h^{-1}g|^{Q+ps}}\, dhdg 
\leq [v]^p_{s,p,U} + C_K\frac {2^Q\sigma_Q}{sp\delta^{sp}} \norm{u}_{L^p(K)} \norm{\phi}_{L^\infty}^{p},
\end{multline}
 where  $\delta=d(K,\bG\setminus U)$  using the gauge distance. The semi-norm $ [\phi u]_{s,p,U} $ in \eqref{e:Mphi to Dsp} can be estimated as follows,
 \begin{multline}\label{e:uphi localization 1}
   [v]_{s,p,U} \leq 2^{p-1}\int_{U}  |u(h)|^{p} \int_{L_{h^{-1}}U}  \frac{|\phi(hg) - \phi(h)|^{p} }{|g|^{Q+ps}}\, dg dh +2^{p-1} \int_{U} \int_{U}  \frac{|u(g) - u(h)|^{p} |\phi(h)|^{p} }{|h^{-1}g|^{Q+ps}} \,dg dh\\
   \leq 2^{p-1}C \norm{u}_{L^p(U)}   + 2^{p-1}\norm{\phi}_{L^\infty}[ u]_{s,p,U},
 \end{multline}
using \eqref{e:MVT cpct supp}, with a constant $C$ depending on $\phi$, and lemma \ref{l: phi^p-kernel-lemma}.
Inequalities \eqref{e:Mphi to Dsp} and \eqref{e:uphi localization 1} imply $[{v }]_{s,p}<\infty$.

From the Sobolev embedding \eqref{t: Global Sobolev} we have $u\in L^{p^*(s)}(\bG)$, hence $v\in L^{p^*(s)}(\bG)$. Thus, from \eqref{e:alternative Dsp def)} we obtain $v\in\Dsp(\bG)$ with bounds on the norms. From proposition \ref{l:conv and Sobolev norm} it follows that $v_\eps=\rho_\eps *v \rightarrow v$ in $\Dsp(\bG)$, and by the properties of the convolution $v_\eps\in \mathcal{C}^\infty_0(\Om)$ for all $\eps$ sufficiently small. Thus, $v\in \Dsp(\Om)$. The above estimates show that $M_\phi: \DspG\rightarrow \DspOm$ is continuous. The last claim  follows from \eqref{e:Wsp in Dsp}.

b) 
For $u\in\WsplocOm$ since $\supp \, v \subset K$ we have   $v\in L^{p}(\bG)$. 
From \eqref{e:Mphi to Dsp} and \eqref{e:uphi localization 1} it follows
$v\in\WspG$ and the claimed continuity of $M_\phi$.

The additional claims follow from the continuity of the embeddings
$$\Dsp(\Om)\hookrightarrow \WsplocOm, \quad \Dsp(\Om)\hookrightarrow \DspG,\quad \WspG\hookrightarrow\WspOm ,$$ taking into account the Sobolev embedding \eqref{t: Global Sobolev} and $p<p^*(s)$ for the first inclusion and parts a) and b) for the last two.
\end{proof}

\subsection{Solutions, sub-solutions and super-solutions}\label{ss:weak solutions defs}
Let $1<p<\infty$ and $p'=p/(p-1)$ be its Sobolev conjugate. Let   $\Om$ be an open set and $ F\in L^{p^*(s)'}_{loc}(\Om)$. As noted in  \eqref{e: Main equation sub-sup-solution} a super-solution (resp. sub-solution)  of the equation
$\frlap u  =F$
will mean  a  function $u\in \WsplocOm \cap \tailspace$ such that  for all $0 \leq \phi \in \mathcal{C}_{0}^{\infty}(\Om)$ we have
\begin{equation}\label{e:weak solution}
  \int\limits_{\G}F(g) \phi(g)\, dg  \leq  D(u,\phi)\qquad  \left(\, \text{resp.}\quad \int\limits_{\G}F(g) \phi(g)\, dg  \geq D(u,\phi)\ \right).
\end{equation}
  As usual, a solution will mean a function which is both a sub-solution and a super-solution.  If $F\in L^{p^*(s)'}_{loc}(\Om)$, then \eqref{e:weak solution} will hold for any  $\phi\in \DspOm$, assuming when necessary the non-negativity of $\phi$. This follows from  \eqref{e:Holder for p-functional}
  and the Sobolev embedding theorem \eqref{t: Global Sobolev}. In particular, for a solution  $u\in \WsplocOm \cap \tailspace$ and  $\phi\in\mathcal{C}_{0}^{\infty}(\Om)$ the functions $\phi\, u$ and $\phi(u -l)_{\pm}$ are  legitimate test functions, cf.  proposition \ref{p:localization} and section \ref{ss:chain rule}. 
  
We have the following representation of the fractional operator.
\begin{theorem}\label{t:symm of d form and hypersing int}
For $u,\, \phi\in C_0^\infty(\bG)$
we have
\begin{equation}\label{e:pv integral 2}
{D}(u,\phi)=\int_{\bG}\phi(g)\, \frlap u (g) dg,
\end{equation}
where 
\begin{equation}\label{e:frac L via hypersing app app}
\frlap u(g)= 2\lim_{\eps\rightarrow
0}\int_{\bG\setminus B(g,\eps)} \frac {\tpwr{u(g)-u(h)}}{|
h^{-1}\, g|^{Q+ps}}\,dh.
\end{equation} 
\end{theorem}

\begin{proof}
The proof follows by an adaption of the arguments in the $p=2$ case noting that here $D(u,\phi)$ is not symmetric, see \cite[Theorem 1.1]{Gar-Loi-Vas-2024-FR}.

By a change of variables, the fact that $|h|=|h|^{-1}$, and that  $t^\fep$ is an odd function, we have for any $\eps>0$ the identity
\begin{equation}
 \int_{\bG\setminus B(g,\eps)} \frac {\tpwr{u(g)-u(h)}}{|
h^{-1}\, g|^{Q+ps}}\,dh = \frac12 \int_{\bG\setminus B(g,\eps)} \frac {\tpwr{u(g)-u(gh)}- {\tpwr{u(gh^{-1})-u(g)}}}{|
h|^{Q+ps}}dh 
\end{equation}
For $u\in C^\infty_0(B_R)$, from \eqref{e:precise Lip}, \eqref{e:MVT cpct supp} and  \eqref{e:2nd order difference remainder} we have the next bound for some $d>2$,
\[
\left\vert \frac {\tpwr{u(g)-u(gh)}- {\tpwr{u(gh^{-1})-u(g)}}}{|
h|^{Q+ps}}\right\vert \leq 2 T_2 M_2\,   (R^{d-1}+1)|h|^{d+(p-2)-Q-sp}.
\] 
Now, by \eqref{e:polar coord} we have $|h|^{d+(p-2)-Q-sp}\in L^1_{loc}(\bG)$ exactly when $d+(1-s)p>2$, i.e., $d>2$, taking into account \eqref{e:range for p}. On the other hand, on $\bG\setminus B_R$ as in \cite{Gar-Loi-Vas-2024-FR} we rely on the condition  $u\in \tailspace$. Thus, we have
\[
\frac{\tpwr{u(g)-u(gh)}- {\tpwr{u(gh^{-1})-u(g)}}}{|
h|^{Q+ps}}\in L^1(\bG),
\]
hence
\[
 \lim_{\eps\rightarrow 0}\int_{\bG\setminus B(g,\eps)} \frac {\tpwr{u(g)-u(h)}}{|
h^{-1}\, g|^{Q+2s}}\,dh = \frac12 \int_{\bG} \frac {\tpwr{u(g)-u(gh)}- {\tpwr{u(gh^{-1})-u(g)}}}{|
h|^{Q+ps}}dh. 
\]

By proposition \ref{p:defs of frac operator} we have that 
\[
\frac {\tpwr{u(g)-u(h)} (\phi(g)-\phi(h))}{|
h^{-1}\, g|^{Q+2s}}\in L^1(\bG\times\bG),
\]
hence Fubini's theorem gives
\[
{D}(u,\phi)= \int_{\bG}\frac {D(h)}{|h|^{Q+ps}}dh,\quad D(h)= \int_{\bG}
 {\tpwr{u(g)-u(gh)} (\phi(g)-\phi(gh))}\,
{dg}.
\]
A change of variables  gives
\begin{multline*}
D(h)= \int_{\bG}
 {\tpwr{u(g)-u(gh)}} \phi(g)\,
{dg} - \int_{\bG}
 {\tpwr{u(g)-u(gh)} \phi(gh)}\,
{dg}\\
=\int_{\bG} \phi(g) \left[{\tpwr{u(g)-u(gh)}- {\tpwr{u(gh^{-1})-u(g)}}}\right]\, dg,
\end{multline*}
which allows us to finish the proof of the claim.

\end{proof}

\section{Proof of the local boundedness, theorem \ref{t: Local Boundedness Interpolation Inequality}}\label{s:proof of loc bdd}

The proof of theorem \ref{t: Local Boundedness Interpolation Inequality} is based on De Giorgi's method, see \cite[Theorem 1.1]{DiC-Kuu-Pal-2016} and \cite[Theorem 1.1]{Man-Pal-Pic-2023} for related results in the case of bounded (or vanishing) $F$. We note that our inequality \eqref{e: Local Boundedness Interpolation Inequality} contains a term  $||F_{+}||_{L^{q}(B_{R})}$ on the right hand, while the corresponding result in \cite{Man-Pal-Pic-2023} does not show explicitly the dependence on  $||F_{+}||_{L^{\infty}(B_{R})}$. However, the constant $\textbf{c}$ therein depends on $||F||_{L^{\infty}(B_{R})}$. In fact, the reader should compare \eqref{e: local boundedness I3}, \eqref{e: reg and bound tildek greater eqn} and \eqref{e: reg and bound tildek def}  with equations  \cite[(5.6), (5.7) and (5.8)]{Man-Pal-Pic-2023}, which include $\norm{f}_{L^\infty(B_r)}$ rather than the correct $\norm{f}_{L^\infty(B_r)}/{\tilde k ^{p-1} }$.  In the linear case, $p=2$, an $L^\infty$ result for solutions of a Schr\"odinger type equation and the Yamabe equation on homogeneous groups were proven in \cite[Theorem 1.1 and Theorem 1.2]{Gar-Loi-Vas-2024-OP}. 
\begin{proof}[Proof of  theorem \ref{t: Local Boundedness Interpolation Inequality}]
First, we will prove the claimed inequality in the case $g_0=e$ and $R=1$. The general case follows by dilation and translation.

Let us observe that  
\begin{equation}\label{e:r' def}
1-\frac{p}{{p^*(s)}}=\frac{sp}{Q}>\frac{1}{q}.
\end{equation}
Thus,  the assumption $q>Q/(sp)$ implies that  $F\in L^{p^*(s)'}_{\loc}(\bG)$. Indeed, by H\"older's and Sobolev's inequalities we have
\[
\frac{1}{p^*(s)'} - \frac 1q=1-\frac{1}{{p^*(s)}}- \frac 1q=\frac{sp}{Q}- \frac 1q+ \frac{p-1}{{p^*(s)}} >0,
\]
hence $q>{p^*(s)'}$.

For the rest of the proof it will be convenient to introduce the sets
\begin{equation}\label{e:Akr}
A(k,r)=B_r\cap \{u\geq k \}.    
\end{equation}
Notice that for $\tilde r<r$ and $k<\tilde k$ we have $A(\tilde k,\tilde r) \subset  A(k,r)$ and
\begin{equation}\label{e:prop of Akr}
  |A(\tilde k,\tilde r) | = |B_{\tilde r}\cap \{u-k\geq \tilde k - k \}|\leq \frac{1}{(\tilde k - k)^p}\int_{B_{\tilde r}} (u-k)_{+} ^p\,dg .
\end{equation}

For $j=0,1, \dots$ let
$$r_{j} = \frac{1}{2}(1 + 2^{-j})R,\quad \tilde{r}_{j} = \frac{r_{j} + r_{j+1}}{2},\quad
B_{j} = B_{r_{j}}(g_{0}),\quad \tilde{B}_{j} = B_{\tilde{r}_{j}}(g_{0}).$$
Clearly, we have
\begin{equation}\label{e:proprts of rj}
\frac 12 < {r_{j+1}}<{\tilde{r}_{j}}<{r_{j}}<r_0=1,  \qquad \frac 12<\frac {r_j}{r_{j+1}}<2 \quad\text{and }\quad  r_{j} - {\tilde{r}_{j}} = {\tilde{r}_{j}}-r_{j+1}  =\frac {1}{2^{j+3}} \geq \frac{ r_j}{2^{j+3}}.
\end{equation}
Let $\phi_{j} \in C_{0}^{\infty}(\tilde{B}_{j})$ be a bump function constructed using  lemma \ref{l: existence of test functions} and \eqref{e:proprts of rj},
\begin{equation}\label{e:cut-off for loc bd}
0 \leq \phi_{j} \leq 1, \quad   \phi_{j} \equiv 1 \ \text{on}\  B_{j+1}, \qquad \text{and} \quad |\phi_{j}(g) - \phi_{j}(h)| \leq \frac{2^{j+5}}{r_{j}} |h^{-1}g|.
\end{equation}
In addition,  for $\tilde{k} > 0$  define the increasing sequences
$$k_{j} =  (1-2^{-j}) \tilde{k}\quad \text{ and } \quad \tilde{k}_{j} = \frac{k_{j+1} + k_{j}}{2},$$
hence we have
\begin{equation}\label{e: props of kj}
k_j<\tilde k_j<k_{j+1}<\tilde{k} \qquad\mbox{and}\qquad k_{j+1} - \tilde{k}_{j} =   \tilde{k}_{j} - k_{j} = \frac{\tilde{k}}{2^{j+2}}.
\end{equation}
Finally, let
\begin{equation}\label{e:T1.2 A_j}
\begin{aligned}
\tilde{w}_{j} = (u-\tilde{k}_{j})_{+}, & \qquad w_{j} = (u - k_{j})_{+}, \quad \mbox{and} \quad M_{j}= \left(\int_{B_{j}} w_{j}^{p}\, dg \right)^{\frac{1}{p}},\\
& \text{thus}\quad{w}_{j+1}<\tilde{w}_{j}<{w}_{j}.
\end{aligned}
\end{equation}
By the chain rule in section \ref{ss:chain rule} we have that $\tilde{w}_{j},\ {w}_{j}\in\Wsploc $. Our goal is to eventually apply  lemma \ref{l: iteration lemma} to the averages $M_j$ with a suitable choice of $\tilde k$, see \eqref{e: reg and bound tildek greater eqn}.

From \eqref{e: props of kj} and \eqref{e:T1.2 A_j} we have for $w_{j+1} > 0$, resp. $\tilde{w}_{j} > 0$, we have the inequalities
\begin{equation}\label{e: local boundedness wj-tilde}
\frac 14\cdot \frac{\tilde{k}}{2^{j}}<{k}_{j+1}-\tilde{k}_{j}< \tilde{w}_j,\quad
 \text{resp. }\   \frac 14\cdot\frac{\tilde{k}}{2^{j}} <\tilde {k}_{j}-{k}_{j}<{w}_{j}.
 \end{equation}
Thus, we have
\begin{equation}\label{e:relation between Akr and Mr}
\begin{gathered}
\{\tilde w_j\phi_j\not=0\}\subset A(\tilde{k}_j,\tilde {r}_j)\subset A({k}_j, {r}_j),\\
|\{\tilde w_j\phi_j\not=0\}|\leq |A(\tilde{k}_j,\tilde {r}_j)| 
\leq \frac{1}{(\tilde k_j - k_j)^p}\int_{A(\tilde{k}_j,\tilde {r}_j)} \tilde w_j^p \,dg\leq \frac{1}{(\tilde k_j - k_j)^p}\int_{ B_j}  w_j^p \,dg\leq 4^p\frac{2^{jp}}{\tilde k^p} M_j^p.
\end{gathered}
\end{equation}
Taking into account \eqref{e:cut-off for loc bd}, \eqref{e:T1.2 A_j}, \eqref{e:relation between Akr and Mr} and H\"older's inequality we have
\begin{equation}\label{e: local boundedness LHS}
M_{j+1}^p\leq  \int_{A(\tilde{k}_j,\tilde {r}_j)} |\tilde w_{j}\phi_j|^{p} dg
\leq \left(\int_{\bG} |\tilde{w}_{j}\phi_j|^{p^*(s)}dg\right)^{\frac{p}{p^*(s)}} |A(\tilde{k}_j,\tilde {r}_j)|^{1-\frac{p}{{p^*(s)}}}.
\end{equation}

Proposition \ref{p:localization}, the assumption $sp<Q$, the inequalities \eqref{e:proprts of rj}, and   the localized Sobolev inequality theorem \ref{t: Local Sobolev embedding} 
 imply  the next inequality with  constant $c=c(p,Q,s)$, the latter depending on the best constant in the Sobolev type embedding theorem \eqref{t: Global Sobolev},
\begin{multline}\label{e: local boundedness initial 0}
 \left(\int_{\tilde B_j} |\tilde{w}_{j}\phi_j|^{p^*(s)}dg\right)^{\frac{p}{p^*(s)}} |A(\tilde{k}_j,\tilde {r}_j)|^{1-\frac{p}{{p^*(s)}}}\\
\leq c  \left(
\left(\frac{r_{j}}{{\tilde{r}_{j}}}\right)^{Q}
\left(\frac{\tilde{r}_{j}}{r_{j} -{\tilde{r}_{j}}}\right)
\left(\frac{r_{j}}{r_{j} - {\tilde{r}_{j}}}\right)^{sp} + 1 \right) \,|A(\tilde{k}_j,\tilde {r}_j)|^{1-\frac{p}{{p^*(s)}}}\, D_{B_j}(\tilde{w}_{j}\phi_j,\tilde{w}_{j}\phi_j)
\\
\leq
c {2^{j(p+1)} }\,|A(\tilde{k}_j,\tilde {r}_j)|^{1-\frac{p}{{p^*(s)}}}\,  D_{B_j}(\tilde{w}_{j}\phi_j,\tilde{w}_{j}\phi_j).
\end{multline}
Since $u$ is a sub-solution, an application of Caccioppoli's inequality theorem \ref{t: Caccioppoli} to the right hand side of \eqref{e: local boundedness initial 0} gives the next bound with  a constant $a=a(p,Q,s)$,
\begin{equation}
\label{e: local boundedness initial}
 \left(\int_{\tilde B_j} |\tilde{w}_{j}\phi_j|^{p^*(s)}dg\right)^{\frac{p}{p^*(s)}} |A(\tilde{k}_j,\tilde {r}_j)|^{1-\frac{p}{{p^*(s)}}}
\leq a\, 2^{j(p+1)}  |A(\tilde{k}_j,\tilde {r}_j)|^{1-\frac{p}{{p^*(s)}}}\left( J_{1} + J_{2} + J_{3}\right),
\end{equation}
where $J_{1}, J_{2},$ and $J_{3}$ denote integrals obtained from Caccioppoli's inequality and the properties of $\phi_j$, cf. \eqref{e:cut-off for loc bd},
\begin{gather}\label{e: I's in local boundedness initial}
\begin{aligned}
 & J_1=  \int_{B_{j}} \int_{B_{j}}  \tilde{w}_{j}^{p}(h)\abs{\phi_{j}(g) - \phi_{j}(h) }^{p} K(g,h) \,dg dh, 
\\
& J_2= \left( \sup_{g \in \tilde{B}_{j} } \int_{\G \setminus B_{j}} \tilde{w}_{j}^{p-1}(h) K(g,h) \,dh \right)\, \int_{\tilde B_{j}} \tilde{w}_{j} \phi_j^p\,dg,
\ \text{and}\
 J_3= \int_{\tilde B_{j}} F \tilde{w}_{j}\phi_j^p \,dg .
\end{aligned}
\end{gather}
To estimate $J_{3}$ in  \eqref{e: I's in local boundedness initial} we use H\"older's inequality, and then Young's inequality with $p$ and its H\"older conjugate $p'=p/(p-1)$ which give
\begin{multline}\label{e: local boundedness I3}
J_3 \leq   \int_{\bG} F_{+}(g)\tilde{w}_{j}\phi_j^p\,  dg
\leq \norm{F_{+}}_{L^q(\tilde B_j)} \left(\int_{\tilde B_j} |\tilde{w}_{j}\phi_j|^{p^*(s)}dg\right)^{\frac{1}{p^*(s)}}  |A(\tilde{k}_j,\tilde {r}_j)|^{1-\frac 1q -\frac{1}{{p^*(s)}}}\\
\leq \delta_j  \left(\int_{\tilde B_j} |\tilde{w}_{j}\phi_j|^{p^*(s)}dg\right)^{\frac{p}{p^*(s)}} + J'_3, \quad J'_3=\frac{C_p}{\delta_j^{1/(p-1)}} \norm{F_{+}}^{p'}_{L^q(\tilde B_j)}  |A(\tilde{k}_j,\tilde {r}_j)|^{p'\left(1-\frac 1q -\frac{1}{{p^*(s)}}\right)}.
\end{multline}
Taking into account \eqref{e: local boundedness initial} we let
\begin{equation}\label{e:choice of dj}
\delta_j=\frac{ 2^{-j(p+1)} }{2a}, 
\end{equation}
which allows us to infer from \eqref{e: local boundedness LHS} and \eqref{e: local boundedness initial} the next bound with a constant $c=c(s,p,Q)$,
\begin{equation}\label{e: local boundedness initial improved 1}
     M_{j+1}^p
\leq c\,2^{j(p+1)}  |A(\tilde{k}_j,\tilde {r}_j)|^{1-\frac{p}{{p^*(s)}}}\left( J_{1} + J_{2} + J'_{3}\right).
\end{equation}
 Furthermore, the assumption $q>Q/(sp)$ is equivalent to
\begin{equation}\label{e:q condition}
   1-\frac{p}{{p^*(s)}}+ {p'\left(1-\frac 1q -\frac{1}{{p^*(s)}}\right)}=1+ p'\left(\frac{sp}{Q} -\frac{1}{q} \right)>1.
\end{equation}
Since
\[
|A(\tilde{k}_j,\tilde {r}_j)|\leq |A( {\tilde{k}}/{2}, R)|\leq  \left( \frac{2}{\tilde k} \right)^p\int_{B_R}u^p,
\] 
we have $|A(\tilde{k}_j,\tilde {r}_j)|<1$ provided $\tilde k$ is sufficiently large, 
\begin{equation}\label{e:1st cond on k}
  \tilde k> 2\norm{u}_{L^p(B_1)}.
\end{equation}
Assuming this condition for the rest of the proof, it follows from \eqref{e: local boundedness initial improved 1} that for $\gamma$ defined by
\begin{equation}\label{e:def of gamma}
\gamma =\min \left\{ \frac{sp}{Q},\, p'\left(\frac{sp}{Q} -\frac{1}{q} \right)\right\}=\left\{
                  \begin{array}{ll}
                   p'\left(\frac{sp}{Q} -\frac{1}{q} \right), & {{Q}/{(ps)}<q\leq  {Q}/{s} ;}\\
  {(sp)}/{Q},  & {q> {Q}/{s},}
                  \end{array}
                \right.,
\end{equation}
noting that $$\frac{sp}{Q} - p'\left( \frac{sp}{Q} -\frac{1}{q}\right)= p'\left( \frac 1q  - \frac sQ\right),$$
we have  the next inequality  with a constant $c=c(p,Q,s)$,
\begin{multline}\label{e: local boundedness initial improved 2}
     M_{j+1}^p
\leq c\,2^{j(p+1)}\left[\left( J_{1} + J_{2}\right)|A(\tilde{k}_j,\tilde {r}_j)|^{\gamma} + {2^{j\frac{p+1}{p-1}}}\,|A(\tilde{k}_j,\tilde {r}_j)|^{1+\gamma}\, \norm{F_{+}}^{p'}_{L^q(\tilde B_j)}\right]\\
\leq c\,2^{j(p+1)}\left[ \left( J_{1} + J_{2}\right)\left(\frac{2^{jp}}{\tilde k^p} M_j^{p}\right)^\gamma  +  {2^{j\frac{p+1}{p-1}}} \norm{F_{+}}^{p'}_{L^q(\tilde B_j)}\ \left(\frac{2^{jp}}{\tilde k^p} M_j^{p}\right)^{1+\gamma}\right],
\end{multline}
invoking \eqref{e:relation between Akr and Mr} in the last inequality.  Next, we turn to the estimates of $J_1$ and $J_2$.

Using \eqref{e:cut-off for loc bd}, \eqref{e:proprts of rj} and \eqref{e:T1.2 A_j}  we obtain, with a constant $c=c(p,Q,s)$, the inequality
\begin{equation}\label{e: local boundedness I1}
J_{1} \leq  \frac {2^{jp}}{r_{j}^{p} } \int_{B_{j}} \int_{B_{j}} |h^{-1}g|^{-Q-sp+p} \tilde{w}_{j}^{p}(h) \,dgdh\leq  c\frac{2^{jp}{r_{j}^{p(1-s)}}}{r_{j}^p}\int_{B_{j}} w_{j}^{p}(h) \,dh =c  \frac {2^{jp}}{r_{j}^{sp} }\, M_j^p,
\end{equation}
taking into account that  for  $h \in B_{j}= B_{r_{j}}(g_{0})$ we have $ B_{r_{j}}(g_{0}) \subset B_{2r_{j}}(h)$  and  lemma \ref{l: phi^p-kernel-lemma}.

To estimate the term $J_{2}$ in  \eqref{e: I's in local boundedness initial},  note that for $g \in  \tilde{B}_{j}$ and $h \in \G \setminus B_{j}$  the triangle inequality and \eqref{e:proprts of rj} give
$$ |g^{-1}h| \geq |g_{0}^{-1} h| -|g^{-1}g_{0}|\geq r_j -\tilde {r}_j\geq \frac {r_j}{2^{j+3}} \geq \frac {|h^{-1}g_{0}|}{2^{j+3}}.$$
Thus, we have
\begin{equation}\label{e: local boundedness I2 a}
J_{2} 
\leq 2^{(j+3)(Q+sp)} \left( \int_{\G \setminus B_{j}} \frac{\tilde{w}_{j}^{p-1}(h)}{|h^{-1}g_{0}|^{Q+sp}} dh\right) \left(\frac {2^{j+2}}{\tilde{k}}\right)^{p-1} M_j^p,
\end{equation}
taking into account  \eqref{e:T1.2 A_j} and the second inequality of \eqref{e: local boundedness wj-tilde}.
Now, from \eqref{e:proprts of rj}, \eqref{e:T1.2 A_j}, and the definition of the tail \eqref{e: Tail} we obtain with a constant $c=c(p,Q,s)$  the bound
\begin{equation}\label{e: local boundedness I2}
J_{2}
\leq c \frac {2^{j(p-1+Q+sp)}}{r_j^{sp}}\left(\frac{\Tail(w_{0}; g_{0}, 1/2)}{\tilde{k}}\right)^{p-1}\, M_j^p.
\end{equation}
From \eqref{e: local boundedness initial improved 2}, \eqref{e: local boundedness I1}, \eqref{e: local boundedness I2}, and \eqref{e:proprts of rj} we obtain
\begin{multline}\label{e: local boundedness initial improved 3}
     M_{j+1}^p
\leq c\,2^{j(p+1)}\left[ \left( \frac {2^{jp}}{r_{j}^{sp} } + \frac {2^{j(p-1+Q+sp)}}{r_j^{sp}}\left(\frac{\Tail(w_{0}; g_{0}, 1/2)}{\tilde{k}}\right)^{p-1}\right) M_j^p\left(\frac{2^{jp}}{\tilde k^p} M_j^{p}\right)^\gamma \right.\\ 
\left.+  {2^{j(2+p(1+\gamma))}} \norm{F_{+}}^{p'}_{L^q(\tilde B_j)}\ \left( \frac{M_j}{\tilde k} \right)^{p(1+\gamma)} \right]\\
 \leq c\, \frac{2^{j(3p+2+p\gamma +Q+sp)}}{r_{j}^{sp}}\left[ 1 + \frac{\left({\Tail(w_{0}; g_{0}, 1/2)}\right)^{p-1}}{\tilde{k}^{p-1}}  +  \frac{\norm{F_{+}}^{p'}_{L^q(\tilde B_j)}}{\tilde k^{p}}\ \right]    \frac{{M_j^{p(1+\gamma)}}}{\tilde k^{p\gamma}}.
\end{multline}
Thus, with constants $b=2^{3+\frac {2+Q}{p} +\frac{sp}{Q}+s}>2^{7}$ and  $C_0=C_0(p,Q,s)=c^{1/p}$, where $c$ is the constant in \eqref{e: local boundedness initial improved 3}, 
 we have the next inequality  
\begin{equation}\label{e: local boundedness recursion initial}
\frac{M_{j+1}}{\tilde{k}} 
\leq C_0 {b^j} \left[ 1 + \frac{\left({\Tail(w_{0}; g_{0}, 1/2)}\right)^{p-1}}{\tilde{k}^{p-1}}  +  \frac{\norm{F_{+}}^{p'}_{L^q(\tilde B_j)}}{\tilde k^{p}}\ \right]^{\frac{1}{p}} \left(\frac{M_j}{\tilde k}\right)^{1+\gamma}.
\end{equation}
Let us note that since $u  \in  \tailspace$ we have $\Tail \left(w_{0}; g_{0}, {1}/{2}\right)<\infty$.
Therefore, for any $\delta \in (0,1]$ if we assume, in addition to \eqref{e:1st cond on k}, the inequality
\begin{equation}\label{e: reg and bound tildek greater eqn}
\tilde{k} \geq \delta \left[\Tail\left(w_{0}; g_{0}, {1}/{2}\right) + \norm{F_{+}}_{L^q( B_1)}^{\frac{1}{p-1}}\right],
\end{equation}
then \eqref{e: local boundedness recursion initial} implies the next inequality with  $\gamma\in (0, sp/Q)$ defined in \eqref{e:def of gamma} and $b=b(Q,s,p)>2^7$,
\begin{equation}\label{e: local boundedness recursion}
\frac{M_{j+1}}{\tilde{k}}
\leq  \frac {C_0 }{\delta}\, b^{j} \,\left(\frac{M_{j}}{\tilde{k}}\right)^{1+\gamma}.
\end{equation}
To apply lemma \ref{l: iteration lemma} to the sequence ${M_{j}}/{\tilde{k}}$ with $c_{0} = \frac{C_0}{\delta } $  we  define
\begin{equation}
\label{e: reg and bound tildek def}
\tilde{k} \overset{def}{=} \eps +  \delta \left[\Tail\left(w_{0}; g_{0}, {1}/{2}\right) + \norm{F_{+}}_{L^q( B1)}^{\frac{1}{p-1}}\right] +  {c_0}^{\frac{1}{\gamma}} b^{\frac{1}{\gamma^{2}}} M_{0}>0, \quad c_{0} = {C_0}/\delta,
\end{equation}
noting that with the made choices we have also \eqref{e:1st cond on k} 
taking into account that $1<1/\delta$ and we can assume  $C_0\geq 1$.

From Lemma \ref{l: iteration lemma} it follows $M_{j} \rightarrow 0$, thus with a constant $C=C(Q,p,s)$ we have the inequality
\begin{equation}\label{e: Local Boundedness eps}
\sup_{B_{1/2}} u_{+} \leq  \tilde{k}
\leq  \eps +  \delta \Tail \left(u_{+}; g_{0}, {1}/{2}\right)
+ \delta  ||F_{+}||_{L^{q}(B_{1})}^{\frac{1}{p-1}}
+  C {\delta^{-\frac{1}{\gamma}} } \left( \fint_{B_{1}} u^{p}_{+} \right)^{\frac{1}{p}}.
\end{equation}
Letting $\eps\rightarrow 0$ we obtain \eqref{e: Local Boundedness Interpolation Inequality} in the case $R=1$  for the positive part of a sub-solution. The general case follows by applying the case $R=1$ to the function $v(g)=u(g_0\delta_\lambda g)$ noting that
\begin{equation}
\begin{gathered}
\frlap v (g)= F_\lambda(g)=\lambda^{sp}F(g_0\delta_\lambda g), \ T(v;e,r)= T(v;g_0,\lambda r), \\ \norm{F_\lambda}_{B_r(e)} = \lambda^{sp-Q/q}\norm{F}_{L^{q}(B_{\lambda r}(g_0))}.
\end{gathered}
\end{equation}
\end{proof}

For a quick reference we end this section with the following well known lemma, which was used in \eqref{e: Local Boundedness eps} in the last proof.
\begin{lemma}
\cite[Lemma 7.1]{Giu-2003}
\label{l: iteration lemma}
Let $\beta >0$ and $A_{j}$ be a sequence of real positive numbers such that
$$A_{j+1} \leq c_{0}b^{j}A_{j}^{1+ \beta}$$
with $c_{0} >0$ and $b > 1$.

If $A_{0} \leq c_{0}^{-\frac{1}{\beta}} b^{-\frac{1}{\beta^{2}}}$ then we have
$$A_{j} \leq b^{-\frac{j}{\beta}}A_{0},$$
which in particular yields $A_{j} \rightarrow 0$.
\end{lemma} 

Let us observe that  theorem \ref{t: Local Boundedness Interpolation Inequality} implies the following corollary.
\begin{cor}\label{c:Local Boundedness Interpolation Inequality}
If $u $ is a solution to $\D_{s,p}u= F$ with $F\in  L^{q}_{\loc}(\bG)$ and $q>Q/(sp)$, then $u\in L^\infty_{\loc}(\bG)$. Furthermore there is a constant  $C=C(Q,p,s)$ such that for $g_0\in\bG$ and $\delta \in (0,1]$  we have the inequality
\begin{equation}\label{e: Local Boundedness Interpolation Inequality for solutions}
\sup_{B_{R/2}} |u| \leq   C {\delta^{-\frac{1}{\gamma}} }\left( \fint_{B_{R}} |u|^{p} \right)^{\frac{1}{p}}+
 \delta \Tail \left(u; g_{0}, {R}/{2}\right)
+ \delta R^{\frac {1}{(p-1)} \left(sp-\frac Qq\right) }
||F||_{L^{q}(B_{R})}^{\frac{1}{p-1}}
,
\end{equation}
where $\gamma$ is defined in \eqref{e:def of gamma}.
\end{cor}

\begin{proof}
The claim  follows by applying \ref{t: Local Boundedness Interpolation Inequality}  to $u_{+}$ and $u_{-}$ taking into account the chain rule \ref{ss:chain rule} and the definitions of sub- and super- solutions \ref{ss:weak solutions defs}.

\end{proof}

\section{An estimate of the tail of a  sub-solutions}\label{ss:tail est}
The following estimate of the tail was proven for solutions on bounded domains with Dirichlet boundary data in the Euclidean case with $F=0$ in \cite[Lemma 4.2]{DiC-Kuu-Pal-2014} and in  the Heisenberg group for bounded  $F$ \cite[Lemma 4.4]{Pal-Pic-2022}, assuming that $u$ is non-negative on a some ball and $p \in (1,2n/(1-s))$. 

\begin{lemma}
\label{l: bound on tail lemma} Let  $q>Q/(sp)$, $F\in L^{p^*(s)'}_{\loc}(\bG)$ and  $F_{-} \in L^{q}_{\loc}(\bG)$. Let $u$ be a sub-solution to $\D_{s,p}u= F$.  Let $0<r<R$, $g_{0} \in \G$ and $B_{r} = B_r(g_{0})$. If  $u$ is non-negative in $B_{R}$, then there is a constant $c = c(Q,s,p)$ such that
$$\Tail(u_{+};g_{0},r)\
\leq \ c\sup_{B_{r}} u \ + \ cr^{\frac{1}{p-1} \left({{sp}-\frac Qq}\right)}\, \norm{F_{-}}^{1/(p-1)}_{L^{q}(B_{r})} \ + \ \left(\frac{r}{R}\right)^{\frac{sp}{p-1}}\Tail(u_{-};g_{0}, R).$$
\end{lemma}

\begin{proof}
First note that if $\sup_{B_{r}} u = \infty$ the claim is trivial. Thus, we suppose that $u$ is essentially bounded on $B_{r}$ and let $k = \sup_{B_{r}} u.$ Using lemma \ref{l: existence of test functions}, let $\phi \in C_{0}^{\infty}(B_{3r/4})$, $ \phi \equiv 1$ on $B_{r/2}$,  $|\phi(g) - \phi(h)| \leq 16/r$ and define
\begin{equation}\label{e:eta bump}
    \eta \overset{def}{=} (u-2k) \phi^{p}.
\end{equation}
Notice that $\eta$ is supported in $B_{r}$ and  $-2k\leq \eta \leq 0$. Since $u$ is a weak sub-solution  we have
\begin{equation}
\label{e: initial decomp for bound on tail}
\int\limits_{B_{r}} F(g)\eta(g) dg\geq D_{r}(u,\eta) + 2I_{r}(u,\eta)
\end{equation}
where $D_{r}(u,\eta)$ and $I_{r}(u,\eta)$ are defined in \eqref{e:local forms balls}.

To estimate $I_{r}(u,\eta)$ we use lemma \ref{l:test fn level set} which gives the next inequality
\begin{multline}\label{e: I2 + I3 for bound on tail prelim}
 I_{r}(u,\eta) \geq  \frac {k}{2^{p-1}}\int\limits_{\G\setminus B_{r}}u_{+}^{p-1}(h) \int_{B_{r}} \frac{\phi^{p}(g)}{|h^{-1}g|^{Q+sp}}\, dgdh  - 2^pk\int_{\bG \setminus B_{r}}  u(h)_{-}^{p-1} \int_{B_{r}}\frac{\phi^{p}(g)}{|h^{-1}g|^{Q+sp}}\, dgdh \\
- (1+2^p)k^p\int_{\G\setminus B_{r}} \int_{B_{r}} \frac{\phi^{p}(g) }{{|h^{-1}g|^{Q+sp}}}\, dgdh.
\end{multline}
Next, we use the properties of $\phi$ and the triangle inequality to obtain the following estimates for  $h\in \bG\setminus B_r$,
\[
\int_{B_{r}} \frac{\phi^{p}(g)}{|h^{-1}g|^{Q+sp}} dg \geq \int_{B_{r/2}} \frac{dg}{|h^{-1}g|^{Q+sp}}  \geq \frac{1}{2^{Q+sp}} \int_{B_{r/2}} \frac{dg}{|h^{-1}g_0|^{Q+sp}}  = \frac {|B_r|}{2^{2Q+sp}|h^{-1}g_0|^{Q+sp}}
\]
and, with a constant $C=C(s,p,Q)$,
\[
 \int_{B_{r}}\frac{\phi^{p}(g)\, dg}{|h^{-1}g|^{Q+sp}} \leq \int_{B_{3r/4}}\frac{dg}{|h^{-1}g|^{Q+sp}} \leq    \frac {C}{|h^{-1}g_0|^{Q+sp}},
\]
using also corollary \ref{c: ker_est}  for the last inequality.
Thus, with a constant $C=C(s,p,Q)$ we have the inequality
\begin{multline}\label{e: I2 + I3 for bound on tail}
 I_{r}(u,\eta) 
\geq  \frac {k}{2^{2Q+sp+p-1}}  |B_{r}| r^{-ps} \left(\Tail(u_{+}; g_{0}, r)\right)^{p-1}  - Ck|B_{r}| R^{-sp}(T(u_{-}; g_{0}, R))^{p-1} - C k^{p} r^{-ps} |B_{r}|,
 \end{multline}
recalling that $u_{-}=0$ on $B_R$, thus the domain of the first integral in the last term in the first line  of \eqref{e: I2 + I3 for bound on tail prelim} can be replaced with $\bG\setminus B_R$.

To estimate $D_{r}(u,\eta)$ we use lemma \ref{l: p-coerc_lemma} for  $w = u-2k$, $w_0=w(g)$ and $w_1=w(h)$, which gives
\begin{equation}
\label{e: I1_est for bound on tail }
D_{r}(u,\eta) \geq  -ck^{p} \int\limits_{B_{r}} \int\limits_{B_{r}} \frac{|\phi(g)- \phi(h)|^{p}}{|h^{-1}g|^{Q+ps}} dg dh \geq -ck^{p}r^{-ps}|B_{r}|,
\end{equation}
taking also into account lemma \ref{l: phi^p-kernel-lemma}. 

Finally,  since $\eta < 0$ and  $|u-2k| \leq 2k$ on $B_{r}$ we have
\begin{equation}\label{e:F estimate}
    \int_{B_{r}} F(g) \eta dg = \int_{B_{r}} F_{-}(g) |\eta| dg - \int_{B_{r}} F_{+}(g) |\eta| dg
     \leq 2k||F_{-}||_{L^{q}(B_{r})}|B_{r}|^{\frac {1}{q'}}.
\end{equation}

Combining \eqref{e: initial decomp for bound on tail}, \eqref{e: I2 + I3 for bound on tail}, \eqref{e: I1_est for bound on tail } and \eqref{e:F estimate}   the claimed inequality follows taking into account \eqref{e: triangle with p} and $|B_{r}|^{\frac {1}{q'}-1}= \omega_Q^{\frac {1}{q'}-1} r^{-Q/q}$.
\end{proof}

\section{Expansion of positivity}\label{s:expansion of positi} 
The results in this section borrow elements from the proofs of \cite[Lemma 1.3]{DiC-Kuu-Pal-2016}  and  \cite[Lemma 3.1, Lemma 3.2]{DiC-Kuu-Pal-2014} which dealt with  homogeneous equations in bounded domains in the Euclidean setting and super-solutions $u\in W^{s,p}(\Rn)$.  We consider weak super-solutions of non-homogeneous equations on unbounded domain in a homogeneous group. Furthermore, we allow solutions in the local Sobolev space $u\in\Wsploc$  so that the result can be used in  arguments involving truncations or translations of solutions, which preserve the local space $\Wsploc$ but not the global space $\WspG$.  This invariance is essential for the oscillation arguments in Section \ref{s:holder cont} leading to the H\"older's continuity of weak solutions.

\subsection{Expansion of positivity lemma}\label{ss:expansion of positi}
The following lemma is inspired from \cite[Lemma 3.1]{DiC-Kuu-Pal-2014}.   Here we note a slight discrepancy with the statement found therein concerning the strict inequality in $|B_{r} \cap \{ u < 2\delta k \}|$ which we use below in \eqref{e:first bound on level sets}.  Notice that if  $k = 0$ and  $u \equiv 0$ on $B_{r}$ then
$$|B_{r} \cap \{ u \leq 2\delta k \}| = |B_{r} \cap \{ u \leq 0 \}| = |B_{r}| > 0,$$ but, regardless of the choice of $\bar C$, we have
$$\lim_{\delta \rightarrow 0} \frac{\bar C}{\sigma \ln  (1/2\delta) }|B_{r}| = 0.$$  

The proof of lemma \ref{l: first bound on level sets} also uses elements from the proof of the logarithmic lemma \cite[Lemma 1.3]{DiC-Kuu-Pal-2016} where \eqref{e: level sets combined integral estimate} below was shown  for $F\equiv 0$   and $u\in W^{s,p}(\R^N)$ a weak super-solution of $\D_{s,p} u=0$ in a bounded open set $\Omega\subset \Rn$ with $u=g$ in $\R^n\setminus\Omega$. Later \cite[Lemma 1.4]{Man-Pal-Pic-2023} gave a version of our \eqref{e: level sets combined integral estimate} for weak solutions  $u\in W^{s,p}(\Hn)$ on a bounded domain $\Omega$ in the Heisenberg group, $\D_{s,p} u=f\in L^{\infty}(\Om)$ in $\Omega\subset \R^n$ with $u=g$ in $\Hn\setminus\Omega$, and $u\geq 0$ on a ball $B_R(x_0)$, but the estimate contained $\norm{f}_{L^\infty(B_r)}$ in the right-hand side.  We mention also \cite[Lemma 5.1 \& Lemma 5.2]{FangZhang24} for related results on the Heisenberg group with $F=0$.
\begin{lemma}
\label{l: first bound on level sets} Let  $q>Q/(sp)$, $F\in L^{p^*(s)'}_{\loc}(\bG)$ and  $F_{-} \in L^{q}_{\loc}(\bG)$. 
Let $u$ be super-solution of $\frlap u =F$ such that $u$ is non-negative on $B_{R}=B_R(g_0)$. 
There is a constant $\bar C=\bar C(p,Q,s)$ 
such that if for some $k \geq 0 $, $\sigma \in (0,1]$ and $r$, $0 < 16 r < R$  the following density inequality holds true,
\begin{equation}\label{e:density at k}
|B_{r} \cap \{ u \geq k\}| \geq \sigma |B_{r}|, \quad B_{r}=B_r(g_0),
\end{equation}
then  for all $a \in (0,1/2]$ we have
\begin{equation}\label{e:first bound on level sets}
\left|B_{6r} \cap \left\{ u <  k a - \frac{d}{2} \right\} \right| \leq \frac{\bar C}{\sigma \ln  (1/a) } |B_{6r}|,
\end{equation}
where $d$ is any non-negative number satisfying
\begin{equation}\label{e:def of d}
d \geq d_{R,r}(g_0)=\left(\frac{r}{R}\right)^{\frac{sp}{p-1}}\Tail(u_{-};g_{0}, R)\ + \  \norm{F_{-}}_{L^q(B_{R})}^{\frac{1}{p-1}}  r^{\frac{Q}{p-1}\left({\frac{sp}{Q}-\frac 1q }\right)}.
\end{equation}
\end{lemma}

\begin{proof}  We will first consider the case   $d>0$.
The first step of the proof is to prove inequality \eqref{e: level sets combined integral estimate final}. Let $w = u + d$, $d>0$. Using lemma \ref{l: existence of test functions}, let $\phi \in C_{0}^{\infty}(B_{7r})$ with $\phi \equiv 1$ on $B_{6r}$  such that  $$|\phi(g) - \phi(h)| \leq \frac{8}{r}|h^{-1}g|.$$
 Since  $u$ is a super-solution to \eqref{e: Main equation} we have, recalling \eqref{e:local forms balls} and \eqref{e:weak solution}, 
\begin{equation}\label{e: decomp estimate}
\int_{B_{7r}} F_{+} w^{1-p} \phi^p\,  dg - \int_{B_{7r}} F_{-} w^{1-p} \phi^p\,  dg =\int\limits_{\G}F w^{1-p} \phi^{p} dg \leq D_{7r}(u,w^{1-p}\phi^{p}) + 2I_{7r}(u,w^{1-p}\phi^{p}).
\end{equation}
From lemma \ref{l: p-coerc_lemma log}, noting that $w>0$ on $\supp \phi$,  there is a constant  $c = c(p,Q,s)$ such that we have
\begin{multline}
\label{e: level sets I1}
D_{7r}(u,w^{1-p}\phi^{p}) \leq -\frac{1}{c}\int_{B_{r}} \int_{B_{r}} \frac{\min \{ \phi(g), \phi(h)\}^p }{|h^{-1}g|^{Q+ps}}\left| \ln  \frac{ |w(g)|}{|w(h)|} \right|^{p} \, dgdh + c D_r(\phi,\phi) \\
\leq 
-\frac{1}{c}\int_{B_{7r}} \int_{B_{7r}}\frac{\min \{ \phi(g), \phi(h)\}^p }{|h^{-1}g|^{Q+ps}}\left| \ln  \frac{ w(g)}{w(h)} \right|^{p} \, dgdh + cr^{Q-sp}\\
 \leq -\frac{1}{c}\int_{B_{6r}} \int_{B_{6r}}\left| \ln  \frac{ w(g)}{w(h)} \right|^{p} K(g,h)\,dgdh + cr^{Q-sp},
\end{multline}
taking into account  \ref{l: ker_est} and $\phi\equiv 1$ on $B_{6r}$.

To estimate the term $I_{7r}(u,w^{1-p}\phi^{p})$ in \eqref{e: decomp estimate}  we use the trivial inequality $a-b\leq a_{+}+b_{-}$
and then the monotonicity of   the function $ t^\fep$, see {section \ref{a:A}}.  Thus, we have the inequality, cf. \eqref{e: triangle with p},
\begin{equation}\label{e:french power estimate 1}
\tpwr{w(g)-w(h)}\leq \paren{w_{+}(g) +w_{-}(h)}^{p-1}
\leq 2^{p-1}\paren{ w_{-}^{p-1}(h)+w_{+}^{p-1}(g) }, \quad g,\, h \in \bG.
\end{equation}
Notice that  $w_{+}(g)=w(g)\geq d > 0$ for $g\in B_R$ and $w_{-}(h)=0$ for  $h\in B_R$.
Thus, we have
\begin{multline}\label{e: level sets I22}
I_{7r}(u,w^{1-p}\phi^{p})
\leq 2^{p-1}\int_{\bG\setminus B_{r}} \int_{B_{7r}}w_{-}^{p-1}(h) w^{1-p}(g) \phi^{p}(g) K(g,h)\, dgdh
\\
+  2^{p-1}\int_{\bG\setminus B_{r}} \int_{B_{7r}} w^{p-1}(g)w^{1-p}(g) \phi^{p}(g) K(g,h)\, dgdh
\leq  cr^{Q} d^{1-p} R^{-sp}  [\Tail(u_{-}; g_{0}, R)]^{p-1} + cr^{Q-sp} ,
\end{multline}
with a constant $c = c(p,Q,s)$, taking into account Lemma  \ref{l: ker_est} for the estimate of the second  integral where $w$ cancels.

Combining the estimates in  \eqref{e: decomp estimate}, 
\eqref{e: level sets I1}  and \eqref{e: level sets I22}, and  the properties of $\phi$ we obtain with a constant $c = c(Q,p,s)$ the next inequalities
\begin{multline}
\label{e: level sets combined integral estimate}
\int\limits_{B_{6r}} \int\limits_{B_{6r}}\left| \ln \frac{w(g)}{w(h)} \right|^{p} K(g,h)dgdh +  \int_{\bG} F_{+} w^{1-p} \phi^p\,  dg \\
    \leq cr^{Q-sp}
    + \int_{B_{7r}} F_{-} w^{1-p}\phi^p\, dg
    +  c \frac {r^{Q}} {d^{p-1}  R^{sp}}  [\Tail(u_{-}; g_{0}, R)]^{p-1}\\
     \leq  cr^{Q-sp} +  \frac {\norm{F_{-}}_{L^q(B_{7r})} |B_{r}|^{\frac {1}{q'}}}{d^{p-1}} + c \frac {r^{Q}} {d^{p-1}  R^{sp}}  [\Tail(u_{-}; g_{0}, R)]^{p-1}\\
    \leq cr^{Q-sp} + c\frac {\norm{F_{-}}_{L^q(B_{7r})}r^{sp-Q/q} }{d^{p-1}}r^{Q-sp} + c \frac {[\Tail(u_{-}; g_{0}, R)]^{p-1} }{d^{p-1}}  \left(\frac{r}{R}\right)^{sp}r^{Q-sp}
    \leq c   r^{Q-sp},
\end{multline}
using \eqref{e:def of d} and noting that $Q/q'-(Q-sp)=sp-Q/q>0$ and $u\geq 0$ on the ball $B_R$.
From \eqref{e:def of d} 
and \eqref{e: level sets combined integral estimate} we come to the inequality
\begin{equation}\label{e: level sets combined integral estimate final}
\int_{B_{6r}} \int_{B_{6r}}\left| \ln \frac{w(g)}{w(h)} \right|^{p} K(g,h)dgdh  \leq Cr^{Q-sp}
\end{equation}
for some constant $C=C\left(Q,p,s\right)$.

For  $a\in (0,1/2)$  define
\begin{equation}\label{e:def of v}
v = \min \left\{ \ln \frac{1}{a}, \left(\ln\frac{k+d}{w} \right)_{+} \right\}.
\end{equation}
Clearly, $0 \leq v \leq \ln \frac{1}{a}$  and $v$ vanishes precisely when $w \geq k + d$.  Thus, \eqref{e:density at k} implies
\begin{equation}\label{e:density at k bis}
|B_{6r} \cap \{ v = 0 \}|\geq |B_{r} \cap \{ v = 0 \}| \geq \sigma  |B_{r}|= \frac {\sigma}{6^Q}|B_{6r}|.
\end{equation}
With these observations in hand, letting $v_{6r}$ denote the average of $v$ on the ball $B_{6r}$, it follows
\begin{multline}\label{e: log_1_over_a_initial}
\ln \frac{1}{a} = \frac{1}{|B_{6r} \cap \{v = 0 \}|}
\int_{B_{6r} \cap \{v = 0 \}} \left( \ln \frac{1}{a} - v(g) \right) \volg
\\
\leq  \frac{6^Q}{\sigma|B_{6r}|} \int\limits_{B_{6r}} \left( \ln \frac{1}{a} - v(g) \right) dg
= \frac{6^Q}{\sigma} \left( \ln \frac{1}{a} - v_{6r} \right).
\end{multline}
An integration of  inequality \eqref{e: log_1_over_a_initial} over the set $B_{6r} \cap \{v = \ln {a^{-1}} \}$ yields
\begin{equation}
\label{e:4.1_before_holder}
 \left\vert B_{6r} \cap \left\{v = \ln\frac{1}{a} \right\}\right\vert \, \ln \frac{1}{a}
\leq \frac{6^Q}{\sigma} \int_{B_{6r} \cap \{v = \ln {a^{-1}}\}} v(g) - v_{6r} \volg
\leq \frac{6^Q}{\sigma} \int_{B_{6r}} |v(g) - v_{6r}| \volg.
\end{equation}
From H\"older's inequality and  Poincar\'e's inequality, theorem \ref{t:Poincare}, taking into account $p>1$,  we have the inequality
\begin{equation}\label{e: level sets app of Holder}
\frac{6^Q}{\sigma} \int\limits_{B_{6r}} |v(g) - v_{6r}| \volg
\leq  \frac{C_{p,Q}}{\sigma} |B_{6r}|^{1/p'} \left(  r^{ps} \int\limits_{B_{6r}} \int\limits_{B_{6r}} \frac{|v(g) - v(h)|^{p}}{|h^{-1}g|^{Q + ps}} \volg dh \right)^{1/p}.
\end{equation}
Recalling the definition \eqref{e:def of v} of $v$  and using that the $\min$ function is a Lipschitz function, we have
\begin{equation}\label{e: v-log estimate} |v(g)-v(h)| \leq \left|\ln \frac{k+d}{w(g)} - \ln \frac{k+d}{w(h)}\right| = \left| \ln\frac{w(g)}{w(h)} \right|.
\end{equation}
From \eqref{e:4.1_before_holder}, \eqref{e: level sets app of Holder} and \eqref{e: level sets combined integral estimate final} 
 we obtain with a constant $\bar C=\bar C\left(Q,p,s\right)$ the inequality,
 \begin{equation}\label{e:sub-level estimate} \left\vert B_{6r} \cap \left\{v = \ln\frac{1}{a} \right\}\right\vert\, \ln  \frac{1}{a} \leq \frac{\bar C}{\sigma} |B_{6r}|.
 \end{equation}
At this point we note that for   $0<a\leq \frac 12$ we have $(1-a)d \geq d/2$, hence $ka-(1-a)d \leq ka-d/2$ and
\begin{equation}\label{e:sub-level est L5.1}
    \left\{v = \ln\frac{1}{a} \right\} = \left\{ w \leq a(k + d) \right\} =\left\{u\leq ka-(1-a)d \right\}  \subset \left\{u \leq ka-\frac{d}{2} \right\}.
\end{equation}
The latter together with the inequality \eqref{e:sub-level estimate} implies
\begin{equation}\label{e:final for L5.1}
\left\vert B_{6r} \cap \left\{u \leq ka- \frac{d}{2} \right\}\right\vert\leq \frac{\bar C}{\sigma \ln \frac{1}{a} } |B_{6r}|,
\end{equation}
which gives \eqref{e:first bound on level sets}  when $d>0$.

If  $d_{R,r}=0$ in \eqref{e:def of d} then $\Tail(u_{-}; g_{0}, R) = 0$. Thus, by taking $d\rightarrow 0$ in \eqref{e:final for L5.1} we obtain  \eqref{e:first bound on level sets}   in the case $d=0$, taking into account that $\bar C$ is independent of $d$.

\end{proof}

The proof of the following lemma relies on \cite[Lemma 3.2]{DiC-Kuu-Pal-2014} where the case  $F\equiv 0$  with  $u\in W^{s,p}(\R^N)$ a weak super-solution of $\D_{s,p} u=0$ in a bounded open set $\Omega\subset \Rn$ with $u=g$ in $\R^n\setminus\Omega$ was considered.  Later, \cite[Lemma 4.2]{Pal-Pic-2022} contains a somewhat similar lemma assuming $F\in L^\infty(\Om)$, but \eqref{e: intro to inf conclusion} below is shown to hold on a ball contained in $B_r$, which is not adequate for our application in \eqref{e:use of exp of positivity}.

\begin{lemma}\label{l: intro of inf}
Let  $q>Q/(sp)$, $F\in L^{p^*(s)'}_{\loc}(\bG)$ and  $F_{-} \in L^{q}_{\loc}(\bG)$.
Let $u$ be super-solution of $\frlap u =F$ such that $u$ is non-negative on $B_{R}=B_R(g_0)$. 
If $\sigma \in (0,1]$,  $0 < 16r < R$ and $k \in \R$ are such that   we have
\begin{equation}\label{e:density condition}
|\{g \in B_{r} : u(g) \geq k \} | \geq \sigma |B_{r}|, \quad B_{r}=B_{r}(g_0),
\end{equation}
then there exists a $\bar\delta=\bar\delta\left(\sigma, p,Q,s\right)\in (0,1/4)$ for which
\begin{equation} \label{e: intro to inf conclusion}
    \inf_{B_{4r}} u \geq \delta k-d
\end{equation}
 for all $0<\delta\leq \bar\delta$ and 
 $d$ is as in Lemma \ref{l: first bound on level sets}, i.e, $
d \geq d_{R,r}(g_0)=\left(\frac{r}{R}\right)^{\frac{sp}{p-1}}\Tail(u_{-};g_{0}, R)\ + \  \norm{F_{-}}_{L^q(B_{R})}^{\frac{1}{p-1}} \, r^{\frac{Q}{p-1}\left({\frac{sp}{Q}-\frac 1q }\right)}$. 
\end{lemma}

\begin{proof} Let $\bar C=\bar C\left(p,Q,s\right)$ be the constant in  lemma \ref{l: first bound on level sets} and $\bar\delta$ be defined as follows,
\begin{equation} \label{e:def of delta}
 {\ln \frac{1}{2\bar\delta}}= \frac{\overline{C}}{\sigma} c_{0}^{\frac{1}{\beta}} b^{\frac{1}{\beta^{2}}},
\end{equation}
where for some constant $c=c(p)$ (as in \eqref{e: estimate of w- on balls}) we have
\begin{equation}\label{e:def of c0 and b}
c_{0} =   c^{1+\beta}, 
\quad b =  2^{Q+2p+sp},\quad \beta =\frac{p^{*}(s)}{p\,q'} - 1=\frac{Q}{Q-sp}\left(\frac{sp}{Q}-\frac{1}{q} \right)>0.
\end{equation}
We will show that for any $\delta\in (0,\min\{\bar\delta, \frac 14\})$ we have \eqref{e: intro to inf conclusion}.
Let us observe that since $u$ is nonnegative on $B_{R}$,  inequality \eqref{e: intro to inf conclusion} is trivial when  $\delta k -   d \leq 0$. Thus, since $d \geq 0$, we only need consider the case  
\begin{equation} \label{e:condition on deltak}
    \delta k >   d \geq  d_{R,r}=d \geq d_{R,r}(g_0).
\end{equation}
Next we consider balls $B_{j}$ decreasing to $B_{4r}(g_0)$ and level sets defined by $l_{j}$ decreasing to $\delta k$ such that lemma \ref{l: iteration lemma} will give
$$A_{j} = \frac{ |B_{j} \cap \{ u < l_{j} \}| }{|B_{j}|} \rightarrow  0,$$  which would imply that  $|B_{4r} \cap \{ u < \delta k \}|=0$, hence $\inf_{B_{4r}} u \geq \delta k$.  For $j=0,1,\dots$ and $r>0$ define balls  $B_{j} = B_{\rho_{j}}(g_{0})$ and  $\tilde B_{j} = B_{\tilde\rho_{j}}(g_{0})$ letting
$$
\rho_{j} = (4 + \frac {1}{2^{j-1}})r\in (4r,6r) \quad  \mbox{  and  }\quad \tilde{\rho_{j}} = \frac{\rho_{j+1} + \rho_{j}}{2}.
$$
Thus, we have
\begin{equation}\label{e:rhoj in expansion}
4r<\rho_{j+1}<\tilde{\rho_{j}}< \rho_{j}<6r,\quad \frac 12< \frac {\tilde{\rho_{j}} }{\rho_{j}}<2,  \mbox{  and  } \rho_{j} - \tilde{\rho_{j}}=\frac {r}{2^{j+1}} >\frac {\rho_{j} }{2^{j+4}}.
\end{equation}
Let $\phi_{j} \in C_{0}^{\infty}(\tilde B_{{j}})$ with $0 \leq \phi_{j} \leq 1$, $\phi_{j} \equiv 1$ on $B_{j+1}$ and
\begin{equation}\label{e:bump in 2nd expa posit lemma}
    |\phi_{j}(g) - \phi_{j}(h)| \leq \frac{4}{\rho_{j} - \tilde{\rho}_{j}}|h^{-1}g| = \frac{2^{j+3}}{r}|h^{-1}g|.
\end{equation}
Define, in addition, a decreasing sequence of levels, $l_{j} \rightarrow \delta k$, such that, cf.  \eqref{e:condition on deltak},
\begin{equation}\label{e:relation of l_0 and d}
  l_{0} = \frac{3}{2} \delta k < 2 \delta k - \frac{1}{2}  d,\quad l_{j} = (1 +  \frac {1}{2^{j+1}})\delta k <2\delta k, \quad  l_{j}  - l_{j+1} =  \frac {\delta k}{2^{j+2}}\geq \frac {l_{j}}{ 2^{j+3}}, 
\end{equation}
and a decreasing sequence of functions 
$w_{j} = ( u-l_{j})_{-}$  such that $w_{j} \rightarrow (u-\delta k )_{-}$. Thus, letting $\chi_{\{u<l_j\}}$ be the characteristic function of the set $ {\{u<l_j\}}$, we have
\begin{equation}\label{e: w_j_def}
\begin{gathered}
\{u < l_{0} \} \subset \{u < 2 \delta k - \frac{1}{2}  d \},\quad w_j\leq l_j \chi_{\{u<l_j\}}, \\
w_{j} = (l_{j} - u)_{+}
\geq l_{j} - l_{j+1}
\geq \frac{ l_{j}}{2^{j+3}}\quad  \text{on}\ \{u < l_{j+1}\}.
\end{gathered}
\end{equation}
From  \eqref{e: w_j_def}, $\phi_{j} \equiv 1$ on $B_{j+1}$, \eqref{e:rhoj in expansion}, and proposition \ref{p:localization}, using the  localized Sobolev inequality theorem \ref{t: Local Sobolev embedding}  as in \eqref{e: local boundedness initial 0}, we obtain with a constant $c=c(Q,p,s)$  the inequality
\begin{multline}\label{e: positivity lemma 2 - part 2 summary}
\frac{ l_{j}^pA_{j+1}^{\frac {p}{p^*(s)}}}{2^{(j+3)p}}\leq (l_{j} - l_{j+1})^{p} A_{j+1}^{\frac {p}{p^*(s)}}
\leq c\left(\fint_{B_{j}} |w_{j} \phi_{j}|^{p^{*}(s)}\, dg\right)^{\frac {p}{p^*(s)}}
\leq c \frac {2^{j(p+1)} \rho_{j}^{sp}}{|B_j|}  D_{B_j}({w}_{j}\phi_j,{w}_{j}\phi_j).
\end{multline}
We pause for a moment to note that \cite[(3.9)]{DiC-Kuu-Pal-2014} should include one more term in the right-hand side coming from the integral  $(\R^n\setminus B_j)\times B_j$ when the Sobolev inequality is used, compare with \eqref{e: positivity lemma 2 - part 2 summary} where the localized Sobolev inequality theorem \ref{t: Local Sobolev embedding} is used.
To estimate the last integral $D_{B_j}({w}_{j}\phi_j,{w}_{j}\phi_j)$ we use first that $u$ is a super-solution,  
thus we have
\begin{equation} \label{e: int_decomp}
- \int\limits_{\G}F_{-}(g)w_{j}(g)\phi_j^{p}(g) dg\leq \int\limits_{\G} F(g)\, w_{j}(g) \phi_j^{p}(g) \, dg\leq  D_{B_j}(u,w_{j} \phi_j^{p}) + 2I_{B_j}(u,w_{j} \phi_j^{p}).
\end{equation}
The integrals $D_{B_j}(u,w_{j} \phi_j^{p})$ and $I_{B_j}(u,w_{j} \phi_j^{p})$ in \eqref{e: int_decomp} can be estimated exactly as in \cite[p. 1822-1823]{DiC-Kuu-Pal-2014}. In fact, since  $p > 1$, using Caccioppoli's   lemma \ref{l: from_caccioppoli_lemma} for $w_{-}$, and \eqref{e:bump in 2nd expa posit lemma}  we have the inequalities
\begin{multline} \label{e: J1 from Caccioppoli}
D_{B_j}(u,w_{j} \phi_j^{p})\leq  - \frac{1}{c}  D_{B_j}(w_{j} \phi_j,w_{j} \phi_j)
 +   c  \int_{B_{_j}} \int_{B_{_j}} w_{j}^p(g) \,|\phi_j(g) - \phi_j(h)|^{p} K(g,h)\,  dg dh\\
 \leq - \frac{1}{c}  D_{B_j}(w_{j} \phi_j,w_{j} \phi_j) +c  2^{jp}l^{p}_{j} r^{-sp} | B_{j} \cap \{u < l_{j} \}|,
\end{multline}
where $c$ is a constant  depending only on $p$. 
On the other hand, using that  $u$ is nonnegative in $B_{R}$,  some obvious estimates, and  lemma \ref{l:Om complement kernel estimate} we obtain
\begin{multline}\label{e: J2+J3 bounded by measure of ball}
   I_{B_j}(u,w_{j} \phi_j^{p})\leq cl_j |B_{_j} \cap \{u < l_j \}| \sup_{g \in \tilde {B}_j} \int_{\G \setminus B_{\rho_j}}  (l_j + u_{-}(h))^{p-1}  K(g,h) dh\\
\leq c l_{j} 2^{j(Q+sp)}  r^{-sp} \abs{B_{j} \cap \{ u < l_{j} \}} \left(l_{j}^{p-1} + \left(\frac{r}{R}\right)^{sp} (\Tail(u_{-};g_{0},R))^{p-1} \right)
\\
\leq c l_{j}^{p} 2^{j(Q+sp)}  r^{-sp} \abs{B_{j} \cap \{ u < l_{j} \}},
\end{multline}
taking into account  \eqref{e:relation of l_0 and d}   and the assumption \eqref{e:condition on deltak} which implies $\left(\frac{r}{R}\right)^{\frac{sp}{p-1}} \Tail(u_{-},g_{0}, R) \leq \delta k <  l_{j}$. 
Finally, from \eqref{e: w_j_def} and $\phi_j^p\leq \chi_{B_j}$  
we have
\begin{equation}\label{e:F- estimate}
\int_{\bG} F_{-}(g)w_{j}(g)\phi_j^{p}(g) \volg\leq \norm{F_{-}}_{L^q(B_{j})}\,l_j \abs{B_{j} \cap \{ u < l_{j} \}}^{1/q'} = \norm{F_{-}}_{L^q(B_{j})}\, l_j\abs{B_{j}}^{1/q'}A_j^{1/q'}.
\end{equation}

Combining the estimates \eqref{e: positivity lemma 2 - part 2 summary}, \eqref{e: int_decomp}, 
\eqref{e: J1 from Caccioppoli},   \eqref{e: J2+J3 bounded by measure of ball} and \eqref{e:F- estimate}, noting that $A_j\leq 1$, we 
we obtain with a constant $c=c(p)$  the inequality 
\begin{equation}\label{e: estimate of w- on balls}
 A_{j+1}^{\frac {p}{p^*(s)}} \leq c 2^{j(Q+2p+sp)}  \left(1 + \frac {\norm{F_{-}}_{L^q(B_{j})} r^{sp-Q/q}} {l_j^{p-1}}\right) A_j^{1/q'}\leq  c 2^{j(Q+2p+sp)}  A_j^{1/q'},
\end{equation}
using \eqref{e:condition on deltak} and \eqref{e:relation of l_0 and d} for the last inequality. 
Recalling  
\eqref{e:def of c0 and b}, we can put \eqref{e: estimate of w- on balls} in the form
$$A_{j+1} \leq c_{0}b^{j}A_{j}^{1+ \beta}.$$ In order to satisfy the condition $A_{0} \leq c_{0}^{\frac{-1}{\beta}} b^{\frac{-1}{\beta^{2}}}$ of lemma \ref{l: iteration lemma} we use \eqref{e:relation of l_0 and d} and  \eqref{e:first bound on level sets} that imply 
\begin{equation}
    \label{e:def of overline c}
    A_{0} = \frac{|B_{6r} \cap \{ u < \frac{3\delta k}{2} \}|}{|B_{6r}|} \leq \frac{|B_{6r} \cap \{ u < 2\delta k - \frac{1}{2} d \}|}{|B_{6r}|} \leq \frac{\overline{C}}{\sigma}\frac{1}{\ln \frac{1}{2\delta}},
\end{equation}
where $\bar C=\bar C(p,Q,s)$ is the constant in  lemma \ref{l: first bound on level sets}. Since $ \ln \frac{1}{2\delta}$ is a continuous decreasing function of $\delta$ we have that for all $0<\delta\leq \bar\delta$, where $\bar\delta$ is defined in \eqref{e:def of delta}, we have the inequality
\begin{equation}
\frac{\overline{C}}{\sigma} c_{0}^{\frac{1}{\beta}} b^{\frac{1}{\beta^{2}}}={\ln \frac{1}{2\bar\delta}}\leq {\ln \frac{1}{2\delta}} .
\end{equation}
Lemma \ref{l: iteration lemma} gives
$$0 = \lim_{j \rightarrow \infty} A_{j} = \lim_{j \rightarrow \infty} \frac{|B_{j} \cap \{ u < l_{j} \}|}{|B_{j}|}.$$ By construction of $B_{j}$ and $l_{j}$ this implies that $|\{ u < l = \delta k \}| = 0$ on $ B_{4r}= \lim_{j \rightarrow \infty} B_{j} $, i.e.,
$\inf_{B_{4r}} u \geq \delta k,$
which completes the proof.

\end{proof}

In the proof of lemma \ref{l: mean bounded by inf lemma} we will need a stronger version of lemma \ref{l: intro of inf} allowing balls that are possibly not centered at $g_0.$ This is the content of  the following corollary. Notice that $\bar\delta$ is the same for all of the balls $B_r(h_0)$ and is independent of $k$.
\begin{cor}\label{c:  intro of inf}
With all of the assumption  of lemma \ref{l: intro of inf}, the conclusion \eqref{e: intro to inf conclusion}   
 holds true for  $0<32r<R$ and any ball $B_r=B_r(h_0)$, not necessarily centered at $g_0$, as long as $B_r\subset B_{R/2}(g_0)$ and the following stronger than \eqref{e:def of d} condition for $d$  is assumed
\begin{equation}\label{e:def of d better}
d \geq 2^{\frac {Q+sp}{p-1}}d_{R,r}(g_0)= 2^{\frac {Q+sp}{p-1}}\left[
\left(\frac{r}{R}\right)^{\frac{sp}{p-1}}\Tail(u_{-};g_{0}, R)\ + \  \norm{F_{-}}_{L^q(B_{R})}^{\frac{1}{p-1}}  r^{\frac{1}{p-1} \left({{sp}-\frac Qq}\right)}\right].  
\end{equation}
\end{cor}

\begin{proof}
We begin by noting that
 \begin{equation}\label{e:compare d's}
  \left( \frac {2}{3}\right)^{\frac {Q+sp}{p-1}} d_{R,r}(g_0) \leq d_{\frac {R}{2},r}(h_0) \leq  {2}^{\frac {Q+sp}{p-1}} d_{R,r}(g_0).
 \end{equation}
 This follows from  the non-negativity of $u$ in $B_R(g_0)$ which gives
 \[
 \int_{\bG\setminus B_r(h_0) } \frac{u_{-}^{p-1}(g)}{|g^{-1}h_0|^{Q+sp}}\, dg =  \int_{\bG\setminus B_R (g_0)} \frac{u_{-}^{p-1}(g)}{|g^{-1}h_0|^{Q+sp}}\, dg
 \]
 and  the triangle inequality which implies
 \[
 \frac {1}{2}|g^{-1}g_0|\leq |g^{-1}h| \leq \frac {3}{2}|g^{-1}g_0|, \quad h\in B_r(h_0)\subset B_{R/2}(g_0), \ g\in\bG\setminus B_R (g_0).
 \]
From lemma \ref{l: intro of inf}  applied to $B_r(h_0)\subset B_{R/2}(h_0)$ we obtain \eqref{e: intro to inf conclusion} with $B_r=B_r(h_0) $ and any $d\geq d_{\frac R2, r}(h_0)$. Taking $d$ larger, as in \eqref{e:def of d better}, from \eqref{e:compare d's} it follows that \eqref{e: intro to inf conclusion}  holds true for any ball $B_r$ as claimed.
\end{proof}

\section{Proof of the  Harnack inequality, theorem \ref{t: Harnack}}\label{s:harnack ineq}
\subsection{A weak Harnack inequality}
Lemma \ref{l: mean bounded by inf lemma} below  will use the following  result of Krylov \& Safonov, which we recall for ease of reference.

\begin{theorem}\cite[Theorem 7.2]{Kin-Sha-2001}
\label{l:krylov-safonov}
Let $g \in \bG$ and $B_{r}=B_r(g_{0}) $. For a measurable set $A \subset B_{r} $ and  $0 < \gamma < 1$, let
$$[A]_{\gamma} = \cup_{\rho >0} \Big\{  B_{3 \rho}(h) \cap B_{r} :  h \in B_{r} \mbox{ and } |A \cap B_{3 \rho}(h)| > \gamma |B_{ \rho}| \  \Big\}.$$
Then, either $[A]_{\gamma} = B_{r}$ or $|[A]_{\gamma}| \geq \frac {1}{\kappa_0\gamma} |A|$, where $\kappa_0 =2^Q$ is the
doubling constant of the Haar measure.
\end{theorem}
Note that for $\rho$ which contributes a non-empty set to $[A]_{\gamma}$ we have trivially from $A \subset B_{r} $ the inequality
\begin{equation}\label{e:KS rmrk}
  \rho< \frac{1}{\gamma^{1/Q} }r.
\end{equation}
\begin{lemma}\label{l: mean bounded by inf lemma}
Let  $q>Q/(sp)$, $F\in L^{p^*(s)'}_{\loc}(\bG)$ and  $F_{-} \in L^{q}_{\loc}(\bG)$. Let  $u$ be a super-solution of $\frlap u =F$ such that $u$ is non-negative on $B_{R}=B_R(g_0)$ and $2^9r< {R}$.  Then, there are constants $\beta>0$ and $c_3 \geq 1$ both depending  on $p$, $Q$, $s$, cf. Lemma \ref{l: first bound on level sets},  such that for any ball $B_r=B_r(g_0)$ we have
\begin{equation}\label{e: mean bounded by inf lemma}
\left( \fint_{B_{r}} u^{\nu} \right)^{1/\nu} \leq c_3 \inf_{B_{r}} u +c_3\left[\left(\frac{r}{R}\right)^{\frac{sp}{p-1}}\Tail(u_{-};g_{0}, R)\ + \  \norm{F_{-}}_{L^q(B_{R})}^{\frac{1}{p-1}}  \,r^{\frac{1}{p-1} \left({{sp}-\frac Qq}\right)}\right], \quad 0<\nu<\beta/2.
\end{equation}
In particular, if $u$ is non-negative on $\bG$ and $F\geq 0$ on $B_R$, then we  have
\begin{equation}\label{e: mean bounded by inf lemma non-negative}
\left( \fint_{B_{r}} u^{\nu} \right)^{1/\nu} \leq C_{p,Q,s} \inf_{B_{r}} u.
\end{equation}
\end{lemma}

\begin{proof} Let $\delta$ 
is as in lemma \ref{l: intro of inf} corresponding to  $\sigma = \frac{\gamma}{3^{Q}}$  with $\gamma = {1}/{(2\kappa_0)}=2^{-(Q+1)}<1$ with $\kappa_0=2^Q$  the doubling constant of the Haar measure.  
Note that $\delta=\delta(Q,p,s)$.

We will show the existence of $\beta$, $\beta>0$, depending only on  $Q,p,s$ such that
for  $d$  given in \eqref{e:def of d better} and $t>0$ we have 
\begin{equation}\label{e:key for layered cake}
    \frac{|B_{r} \cap \{ u > t - \frac{d}{1-\delta} \} |}{|B_{r}|}  
\leq  \delta^{-\beta} t^{-\beta} \left({m} + \frac{d}{1-\delta}\right)^{\beta}, \quad m={\essinf_{B_{r}} u}.
\end{equation}
The claimed inequality \eqref{e: mean bounded by inf lemma} follows then, as usual, from   the layered cake formula. Indeed, 
any $a>{m} + \frac{d}{1-\delta}\geq 0$ we have 
\begin{multline*}
\fint_{B_{r}} u^{\frac{\beta}{2}} dg
= \frac{\beta}{2} \int_{0}^{\infty} t^{ \frac{\beta}{2} - 1} \frac{|B_{r} \cap \{u > t \} |}{|B_{r}|} dt
= \frac{\beta}{2} \int_{0}^{a} t^{\frac{\beta}{2}} \frac{|B_{r} \cap \{u > t \} |}{|B_{r}|} dt +  \frac{\beta}{2}\int_{a}^{\infty} t^{ \frac{\beta}{2} - 1} \frac{|B_{r} \cap \{u > t \} |}{|B_{r}|} dt 
\\
\leq  \frac{\beta}{2} \int_{0}^{a} t^{ \frac{\beta}{2} - 1}  dt +
 \frac{\beta}{2} \int_{a}^{\infty} t^{ \frac{\beta}{2} - 1} \frac{|B_{r} \cap \{u > t - \frac{d}{1-\delta}\} |}{|B_{r}|} dt
\\
\leq a^ {\frac{\beta}{2}} +  \frac{\beta}{2} \delta^{-\beta} \left({m} + \frac{d}{1-\delta}\right)^{\beta} \int_{a}^{\infty} t^{ \frac{\beta}{2} - 1- \beta} dt 
\leq  a^{\frac{\beta}{2}} - \frac{\frac{\beta}{2}}{\frac{\beta}{2}- \beta} \delta^{-\beta}\left({m} + \frac{d}{1-\delta}\right)^{\beta} a^{\frac{\beta}{2}-\beta} 
\\
\leq a^{\frac{\beta}{2}} + \frac 12 \delta^{-\beta} a^{\frac{\beta}{2}}=\left( 1+ \frac 12 \delta^{-\beta}\right)a^{\frac{\beta}{2}}.
\end{multline*}
Since $a$ can be any number such that $a>{m} + \frac{d}{1-\delta}\geq 0$, inequality \eqref{e: mean bounded by inf lemma} follows from H\"{o}lder's inequality for any $\nu\leq \beta/2$.

We turn to the proof of \eqref{e:key for layered cake}.
For $t > 0$ let $k_i= t \delta^{i} - \frac{d}{1-\delta}$, $i=0,1,\dots$, where $d$ is fixed as in \eqref{e:def of d better}. A  computation shows that $\delta k_{i-1} - d=k_i$, which will be used below. Define further the sets
$$A_{t}^{i} = \left\{ g \in B_{r} : u(g) \geq k_i\right\}. $$
In particular, we have ${|A_{t}^{0}|}={|B_{r} \cap \{ u > t - \frac{d}{1-\delta} \} |},$ which is one of the sets appearing in  \eqref{e:key for layered cake}.
It is enough to consider the case ${|A_{t}^{0}|}>0$ since otherwise \eqref{e:key for layered cake} holds trivially.  Furthermore,  we have the inclusion, cf. also the notation in theorem \ref{l:krylov-safonov},
\begin{equation}\label{e: At containment for k.s}
[A_{t}^{i-1}]_{\gamma} \overset{def}{=} \cup_{\rho >0} \Big\{  B_{3\rho}(h) \cap B_{r} :  h \in B_{r} \mbox{ and } | A_{t}^{i-1}\cap B_{3\rho}(h)| > \gamma |B_{\rho}| \Big\} \subseteq A_{t}^{i},
\end{equation}
where we have let $\gamma = {1}/{(2\kappa_0)}<1$. 
Indeed, for each set $B_{3\rho}(h)$ contributing in a non-trivial way to the union that defines $[A_{t}^{i-1}]_{\gamma}$ we have
\begin{equation*}
|\{g \in B_{3\rho}(h) : u(g) \geq k_{i-1} \}| \geq   
|A_{t}^{i-1} \cap B_{3\rho}(h)| > \gamma |B_{\rho}(h)| = \frac{\gamma}{ 3^{Q}}|B_{3\rho}(h))|.
\end{equation*}
In addition, by \eqref{e:KS rmrk}
we have $ \rho< \frac{1}{\gamma^{1/Q} }r= 2^{(Q+1)/Q}r<4r$. Since $h\in B_r$, from the triangle inequality and $2^9r<R$ we have  $B_{3\rho}(h)\subset B_{16r}\subset B_{R/2}$ and  $32(3\rho)<2^9r<R$. Thus, we can apply  corollary \ref{c: intro of inf} with $\sigma = \frac{\gamma}{ 3^{Q}}$ to each ball $B_{3\rho}(h)$, which gives
\begin{equation}\label{e:use of exp of positivity}
\inf_{B_{12\rho}(h)} u \geq \delta k_{i-1} - d = k_{i},
\end{equation}
hence $ B_{3\rho}(h) \subset A_{t}^{i}  $, i.e., \eqref{e: At containment for k.s} holds true.

Next, we will show that there is a natural number $j$ such that $A_{t}^{j} = B_{r}$. 
By the Krylov \& Safonov's theorem \ref{l:krylov-safonov}, in which we let  $\gamma=1/(2\kappa_0)$, either $[A_{t}^{i-1}]_{\gamma} = B_{r}$  or for every $i$ we have
\begin{equation}\label{e:KS 2nd case}
   |A_{t}^{i}|\geq  |[A_{t}^{i-1}]_{\gamma}|\geq \frac {1}{\gamma \kappa_0} |A_{t}^{i-1}|=2|A_{t}^{i-1}|,
\end{equation} 
taking into account \eqref{e: At containment for k.s}.

If the first case holds for some $j$, from \eqref{e: At containment for k.s} it follows $A_{t}^{j} = B_{r}$.
If the second case holds for all  $i$, recalling that $|A_{t}^{0}| > 0$, let $j$  be the smallest non-negative integer such that 
$|A_{t}^{0}| > 2^{-j}|B_{r}|$.
Then,  from  \eqref{e:KS 2nd case} we have
\[
|A_{t}^{j}|\geq 2 |A_{t}^{j-1}| \geq \dots \geq 2^j|A_{t}^{0}| > |B_{r}|.
\]
Since $A_{t}^{j} \subset |B_{r}|$ we have  $A_{t}^{j} = B_{r}$. Thus, in any case, letting
\[
k=1+ \frac {1}{\ln (\gamma\kappa_0) }\ln \frac{|A_{t}^{0}|}{|B_{r}|}=1+  \frac {1}{\ln 2}\, \ln \frac {|B_{r}|}{|A_{t}^{0}|}
\]
we have  $u > t\delta^{k} - \frac{d}{1-\delta}$ on $B_{r}$. A small calculation shows \eqref{e:key for layered cake} with $\beta {=} \frac{\ln (\gamma\kappa_0)}{\ln  \delta}=\frac{\ln 2}{\ln  (1/\delta)}>0$ since $0<\delta < 1$, which completes the proof of \eqref{e:key for layered cake}, thus the proof of the lemma.

\end{proof}

\subsection{Proof of Theorem \ref{t: Harnack} - Harnack's Inequality}\label{s:harnack}

The proof presented below is based on \cite[Theorem 1.1]{Pal-Pic-2022} and \cite[Theorem 1.1]{DiC-Kuu-Pal-2014}.   We note that \cite{DiC-Kuu-Pal-2014} is stated  for $p \in (1,\infty)$ while \cite{Pal-Pic-2022} allows for
$p \in (1,2n/(1-s))$. Here, we have $p \in (1,Q/s)$ due to  the global Sobolev embedding, theorem \eqref{t: Global Sobolev}, which was used in the proofs of certain preliminary results. A different approach  in the homogeneous case in the Euclidean space can be found in  \cite{Kassmann09}, \cite{Cozzi17}, \cite{DydaKassmann20} and \cite{KuusiMingSire15}.

\begin{proof}[Proof of theorem \ref{t: Harnack}]
In the proof we will continue to use $B_{r}$ for $B_{r}(g_{0})$. Let $\rho \leq r$. Theorem \ref{t: Local Boundedness Interpolation Inequality} implies the bound 
$$ \sup_{B_{\rho/2}} u_{+} \leq  
 \delta \Tail \left(u_{+}; g_{0}, {\rho}/{2}\right)
+ \delta \rho^{\frac {1}{(p-1)} \left(sp-\frac Qq\right) }
||F_{+}||_{L^{q}(B_{\rho})}^{\frac{1}{p-1}}
+  C {\delta^{-\frac{1}{\gamma}} }\left( \fint_{B_{\rho}} u_{+}^{p} \right)^{\frac{1}{p}},$$
where $c=c(Q,p,s)$ and $0<\delta\leq 1$. Applying lemma \ref{l: bound on tail lemma} to the tail term, recalling the assumptions $p > 1$, we arrive at the inequality
\begin{multline}
    \sup_{B_{\rho/2}} u
    \leq c \delta^{-\gamma} \left( \fint_{B_{\rho}} u^{p} dh \right)^{1/p}
    \\
    + c \delta \left(\sup_{B_{\rho/2}} u + c\rho^{\frac{{sp-\frac Qq}}{p-1}} ||F_{-}||_{L^{q}(B_{\rho/2})}^{1/(p-1)}  + \left(\frac{\rho}{2R}\right)^{\frac{sp}{p-1}} \Tail (u_{-}; g_{0}, R)\right) + c\delta \rho^{\frac{{sp-\frac Qq}}{p-1}}||F_{+}||_{L^{q}(B_{\rho})}^{1/(p-1)}
    \\
    \leq c \delta^{-\gamma} \left( \fint_{B_{\rho}} u^{p} dh \right)^{1/p}
    + c \delta \left(\sup_{B_{\rho/2}} u  + \left(\frac{\rho}{2R}\right)^{\frac{sp}{p-1}} \Tail (u_{-}; g_{0}, R)\right) +c\delta\rho^{\frac{{sp-\frac Qq}}{p-1}} ||F||_{L^{q}(B_{\rho})}^{1/(p-1)}.
\end{multline}
Thus, we obtain
$$
(1-c\delta)\sup_{B_{\rho/2}} u
\leq c  \delta^{-\gamma} \left(\fint_{B_{\rho}} u^{p} dh \right)^{1/p} + c\delta\left(\frac{r}{R}\right)^{\frac{sp}{p-1}} \Tail(u_{-}; g_{0}, R)  + +c\delta r^{\frac{{sp-\frac Qq}}{p-1}} ||F||_{L^{q}(B_{r})}^{1/(p-1)}.
$$
Let $1/2 \leq \sigma' < \sigma \leq 3/2$. A covering argument gives the inequality
\begin{multline}
\label{e: sup_sigma'r}
(1-c\delta)\sup_{B_{\sigma' r}} u
\leq  \frac{c\delta^{-\gamma} }{(\sigma - \sigma')^{Q/p}} \left( \fint_{B_{\sigma r}} u^{p} dh \right)^{1/p}
+ c\delta\left(\frac{r}{R}\right)^{\frac{sp}{p-1}} \Tail(u_{-};g_{0}, R)
+ c\delta r^{\frac{{sp-\frac Qq}}{p-1}} ||F||_{L^{q}(B_{r})}^{1/(p-1)}.
\end{multline}
For any $\nu>0$ we have
  $  u^{p} = u^{p-\nu + \nu} \leq (\sup u)^{p-\nu} u^{\nu},$ 
which  allows to infer from \eqref{e: sup_sigma'r} the inequality
\begin{multline*}
(1-c\delta) \sup_{B_{\sigma' r}} u  \leq  \frac{c\delta^{-\gamma} }{(\sigma - \sigma')^{Q/p}}(\sup_{B_{\sigma r}} u)^{(p-\nu)/p}  \left( \fint_{B_{\sigma r}} u^{\nu} dh \right)^{1/p}
+ c\delta \left(\frac{r}{R}\right)^{\frac{sp}{p-1}} \Tail(u_{-};g_{0}, R)\\
+ c\delta r^{\frac{{sp-\frac Qq}}{p-1}} ||F||_{L^{q}(B_{r})}^{1/(p-1)}.
\end{multline*}
Note that increasing the value of $c$ makes the right side bigger and the left side smaller, thus we will assume that $c \geq 1/4$. With $\delta = \frac{1}{4c} \in (0,1]$ we obtain
\begin{multline*}
\frac{3}{4} \sup_{B_{\sigma' r}} u
\leq \frac{c}{(\sigma - \sigma')^{Q/p}}(\sup_{B_{\sigma r}} u)^{(p-\nu)/p}  \left( \fint_{B_{\sigma r}} u^{\nu} dh \right)^{1/p}
+ \frac{1}{4} \left(\frac{r}{R}\right)^{\frac{sp}{p-1}} \Tail(u_{-};g_{0}, R)
+ \frac{1}{4}  r^{\frac{{sp-\frac Qq}}{p-1}} ||F||_{L^{q}(B_{r})}^{1/(p-1)}.
\end{multline*}
From Young's inequality it follows
\begin{multline*}
\sup_{B_{\sigma' r}} u
\leq \frac{2}{3}\sup_{B_{\sigma r}} u + \frac{c_\nu}{(\sigma - \sigma')^{Q/\nu}} \left( \fint_{B_{ r}} u^{\nu} dh \right)^{1/\nu}
+ \frac{1}{4} \left(\frac{r}{R}\right)^{\frac{sp}{p-1}} \Tail(u_{-};g_{0}, R)
+ \frac{1}{4}  r^{\frac{{sp-\frac Qq}}{p-1}} ||F||_{L^{q}(B_{r})}^{1/(p-1)},
\end{multline*}
where
\[
c_\nu=2^{Q/\nu}\left(\frac{3}{2}\right)^{\frac{p-\nu}{\nu}}\frac{\nu}{p}\left({\frac{p-\nu}{\nu}} \right)^{{\frac{p-\nu}{\nu}}}c^{p/\nu}
\]
is a constant depending on $p,Q,s$ and $ \nu$, but  $\nu$ will be fixed below.  An application of lemma \ref{l: iteration lemma 2} for $f(t)=\sup_{B_{t r}}u $ gives with a constant $C=C(p,Q,s,\nu)$ the inequality
$$\sup_{B_{r}} u
\leq  C\left( \fint_{B_{ r}} u ^{\nu} dh \right)^{1/\nu}
+ C\left(\frac{r}{R}\right)^{\frac{sp}{p-1}} \Tail(u_{-};g_{0}, R)
+ C  r^{\frac{{sp-\frac Qq}}{p-1}} ||F||_{L^{q}(B_{r})}^{1/(p-1)}.$$ 
Since $p \in (1,Q/s)$ we can apply lemma \ref{l: mean bounded by inf lemma} with a fixed  $ \nu $  to bound the first term and get
\begin{equation*}
 \sup_{B_{r}}  u
\leq C c_3  \inf_{B_{r}} u
+ C c_3 \left[\left(\frac{r}{R}\right)^{\frac{sp}{p-1}}\Tail(u_{-};g_{0}, R)\ + \  \norm{F}_{L^q(B_{R})}^{\frac{1}{p-1}}  \,r^{\frac{{sp-\frac Qq}}{p-1}}\right] .
\end{equation*}
Recalling that the constant $c_3$  from lemma \ref{l: mean bounded by inf lemma} depends only on $p,Q$ and $s$, the inequality in Theorem \ref{t: Harnack} follows.

\end{proof}

The following iteration lemma is well known and was used in the above proof.

\begin{lemma} \cite[lemma 1.1]{Gia-Giu-1982}\label{l: iteration lemma 2}
Let $f = f(t)$ be a non-negative bounded function defined for $0 \leq T_{0} \leq t \leq T_{1}$.  Suppose there are constants $d_{1}, d_{2}$ such that for $T_{0} \leq t \leq \tau \leq T_{1}$ we have
$$f(t) \leq d_{1}(\tau- t)^{-\theta} + d_{2} + \zeta f(\tau),$$
where $d_{1}$, $d_{2}$, $\theta, \zeta$ are nonnegative constants and $\zeta < 1$.  Then there is a constant $d_{3}$ depending only on $\theta$ and $\zeta$, such that for every $\rho, R$ and $ T_{0} \leq \rho < R \leq T_{1}$ we have
$$f(\rho) \leq d_{3}[d_{1}(R-\rho)^{-\theta} + d_{2}].$$
\end{lemma}

A direct proof of the following weak Harnack inequality in the case $p=2$ in a very general setting can be found in \cite[Theorem 1.6]{DydaKassmann20}. Here, it follows trivially from theorem \ref{t: Harnack}.

\begin{corollary}[Weak Harnack Inequality]\label{c: Weak Harnack}
With the assumptions of theorem \ref{t: Harnack}, 
we have for $\nu\geq 1$ the inequality
\[
\left(\fint_{B_{r}} u^{\nu}(g) \volg \right)^{1/\nu} \leq C
\left[ \inf_{B_{r}(g_{0})} u
+  \left(\frac{r}{R}\right)^{\frac{sp}{p-1}}\Tail(u_{-};g_{0}, R)\ 
+   r^{\frac{{sp-\frac Qq}}{p-1}} \norm{F}_{L^{q}(B_{r})}^{1/(p-1)} \right].
\]
\end{corollary}

\section{Proof of theorem \ref{t: Holder Continuity} - H\"older continuity}\label{s:holder cont}

The proof follows the arguments in the local case as in \cite[Theorem 8.22]{Gil-Tru-1983}. 
Let us observe that since $u$ is a weak solution with $F \in L^{q}_{\loc}(\bG)$ with $q>Q/(sp)$ we have $u \in L^{\infty}_{\loc}(\bG)$ by theorem \ref{t: Local Boundedness Interpolation Inequality} applied to $u$ and $-u$.  Let $M_{t} = \sup_{B_{rt}} u$, $m_{t} = \inf_{B_{rt}} u$, and $\omega$ be the oscillation of $u$ on the ball $B_r=B_r(g_0)$,  $$\omega (r)=\osc_{B_r} u =\sup_{B_{r}}  u - \inf_{B_{r}} u=M_1-m_1.$$
From the chain rule in section \ref{ss:chain rule}, the definition of a weak solution \ref{ss:weak solutions defs} and \eqref{e:tail of constant} it follows that for any constant $C$ the function $C\pm u$ is also a weak solution.   
From the weak Harnack inequality, corollary \ref{c: Weak Harnack}, we have with $C=C(p,Q,s)$ and $$A(r)=      r^{\frac{{sp-\frac Qq}}{p-1}} \norm{F}_{L^{q}(B_{r})}^{1/(p-1)}< {r}^{\frac{{sp-\frac Qq}}{p-1}} ||F||_{L^{q}(B_{R)})}^{1/(p-1)}
, \quad r<\frac{R_0}{2^9},$$
the next inequalities 
\begin{multline}
\label{e:Holder sup minus u}
\fint_{B_{2r}} (M_{4} - u)
\leq C\inf_{B_{2r}} (M_{4} - u) + C\left[\left(\frac{2r}{R_{0}}\right)^{\frac{sp}{p-1}}\Tail((M_{4} -u)_{-};g_{0}, R_{0})  +A(r)\right]
\\
\leq C\inf_{B_{r}} (M_{4} - u) + C\left[\left(\frac{r}{R_{0}}\right)^{\frac{sp}{p-1}}  \Tail((M_{4} -u)_{-};g_{0}, R_{0})  +A(r)\right] 
\\
= C(M_{4} - M_{1}) + C\left[\left(\frac{r}{R_{0}}\right)^{\frac{sp}{p-1}}  \Tail((M_{4} -u)_{-};g_{0}, R_{0})  +A(r)\right]
\\
\leq C(M_{4} - M_{1}) + C\left[\left(\frac{r}{R_{0}}\right)^{\frac{sp}{p-1}} \Tail((u_{+};g_{0}, R_{0})  +A(r)\right].
\end{multline}
In the last inequality we used that $M_4\geq 0$, thus $(M_{4} - u)_{-}\leq u_{+}$, which gives $$\Tail((M_{4} -u)_{-};g_{0}, R_{0})\leq \Tail(u_{+};g_{0}, R_{0}).$$

Similarly, we have the estimates
\begin{multline*}
\label{e: Holder u minus sup 0}
\fint_{B_{2r}} (u-m_{4})\, dg
\leq C\inf_{B_{2r}} (u - m_{4}) + C\left[\left(\frac{2r}{R_{0}}\right)^{\frac{sp}{p-1}} \Tail((u-m_{4})_{-};g_{0}, R_{0}) +A(r) \right] \\
\leq C\inf_{B_{r}} (u - m_{4}) + C\left[\left(\frac{r}{R_{0}}\right)^{\frac{sp}{p-1}} \Tail((u-m_{4})_{-};g_{0}, R_{0}) +A(r)\right]\\
= C(m_{1} - m_{4}) + C\left[\left(\frac{r}{R_{0}}\right)^{\frac{sp}{p-1}} \Tail((u-m_{4})_{-};g_{0}, R_{0}) +A(r)\right] .
\end{multline*}
To estimate the term $\Tail((u-m_{4})_{-};g_{0}, R_{0}) $ we note that $(u-m_{4})_{-}$ is non-zero only when $u < m_{4}$, and if $0 \leq u < m_{4}$ then $(u-m_{4})_{-} < m_{4}$ while if $u < 0 \leq m_{4}$ then $|u-m_{4}|^{p-1} \leq 2^p(|u|^{p-1} + m_{4}^{p-1})$. Therefore we have
\[
\Tail((u-m_{4})_{-};g_{0}, R_{0}) \leq  \Tail(u_{-};g_{0}, R_{0})  + \Tail(m_{4};g_{0}, R_{0}).
\]
Furthermore, from  \eqref{e:tail of constant} we have $\Tail(m_{4};g_{0}, R_{0})=m_4\leq u(g_0)$. Thus, we have the inequality
\begin{equation}\label{e: Holder u minus sup}
\fint_{B_{2r}} (u-m_{4})\, dg
\leq C(m_{1} - m_{4}) + C\left(\frac{r}{R_{0}}\right)^{\frac{sp}{p-1}} \left[ \Tail(u_{-};g_{0}, R_{0}) + u(g_0)\right]+A(r).
\end{equation}
Adding  inequalities \eqref{e:Holder sup minus u} and \eqref{e: Holder u minus sup} we obtain with a constant $ C= C(p,Q,s)$ the inequality
\begin{align}\notag
   M_{4} - m_{4} = & \fint_{B_{2r}} M_{4} - m_{4} \,dg
    \leq C\left( \left(M_{4} - m_{4}\right) -\left( M_1 - m_{1}\right)\right) + k_0 r^{\frac{sp-Q/q}{p-1}}, \\ \label{e:k0 defn}
     k_0= & \left(\frac{1}{R_{0}}\right)^{\frac{sp-\frac Qq}{p-1}} \left[ C\Tail(u;g_{0}, R_{0})
    + C u(g_{0}) +  {R_0}^{\frac{sp-\frac Qq}{p-1}}\norm{F}_{L^{q}(B_{R_0})}^{1/(p-1)}\right],
\end{align}
taking into account the inequality $\Tail(u_{+};g_{0}, R_{0}) + \Tail(u_{-};g_{0}, R_{0}) \leq \Tail(u;g_{0}, R_{0})$.
Thus, the oscillation of $u$ satisfies the inequality
\begin{equation}\label{e: Holder cont oscilation 1}
\omega (r) \leq \left(1-\frac 1C \right) \omega({4r})\, u +  k_0 r^{\frac{sp-Q/q}{p-1}}\leq \gamma \omega({4r})\, u +  k_0 r^{\frac{sp-Q/q}{p-1}}, \quad \gamma=1-\frac {1}{C+1}\in (0,1),
\end{equation}
which, after replacing $r$ with $r/4$ can be put in the form
\begin{equation}
  \omega(\tau r)  \leq \gamma \omega(r) + \sigma (r) ,
\end{equation}
where $\tau = 1/4$ and $\sigma(r)=2^{\frac{-2sp}{p-1}} k_0 \, r^{\frac{sp-Q/q}{p-1}}$. 

The proof of theorem \ref{t: Holder Continuity} follows from  lemma \ref{l: Gil-Tru-1983} stated below, see \cite[Theorem 8.22]{Gil-Tru-1983} for details.

\begin{lemma}\cite[Lemma 8.23]{Gil-Tru-1983}
\label{l: Gil-Tru-1983}
Let $\omega$ be a non-decreasing function on $(0,R_{0}]$ satisfying , for all $R \leq R_{0}$, the inequality
$$\omega(\tau R) \leq \gamma \omega(R) + \sigma(R)$$
where $\sigma$ is also non-decreasing and $0 < \gamma$, $\tau < 1$. Then for any $\mu \in (0,1)$ and $R \leq R_{0}$,
$$\omega(R) \leq C \left(\left(\frac{R}{R_{0}} \right)^{\alpha} \omega(R_{0}) + \sigma(R^{\mu} R_{0}^{1-\mu}) \right)$$
with $C = C(\gamma, \tau)$ and $\alpha = \alpha(\gamma, \tau, \mu)=(1-\mu)\ln  \gamma/ \ln  \tau$.
\end{lemma}

\section{Some useful inequalities}\label{a:A}
In this section we show point-wise inequalities which underline several integral estimates throughout the paper.
Before we turn to  our results we recall the following elementary inequalities that are used frequently in the paper,
\begin{equation}\label{e: triangle with p}
\frac 12(a^q+b^q) \leq  \min \{ 1, 2^{q-1} \}(a^q+b^q)\leq (a+b)^q\leq \max\{1,2^{q-1}\} (a^q+b^q) \leq 2^q(a^q+b^q),
\end{equation}
where   $ a,\, b, \, q \geq 0$. In addition, we mention  \cite[Section 2.2]{IannizzottoMosconiSquassina2016}, \cite[Appendix C]{BrascoLindgrenSch18} and \cite[Chapter 12]{Lindqvist2019} for  other useful inequalities relevant for the equations modeled on the local or non-local p-Laplacian.

\subsection{Inequalities for the signed power function}
We recall that for $q>0$ the signed power function $t^\feq$ is the function on the real line defined by 
 \begin{equation}\label{e:french power}
  t\mapsto t^\feq=|t|^{q-1}t = t_{+}^{q} - t_{-}^{q}.  
 \end{equation}
 From the homogeneity and monotonicity
of the map $t\mapsto t^\feq$ we have for $a,\, b\in\R$ the inequality
\[
 \left(a^{\langle q\rangle}-b^{\langle q\rangle}\right)(a-b)
 \geq c_q |a-b|^{q+2}.
\]
The function $t^\feq$ is continuously differentiable on $\R\setminus\{0\}$, 
\[
 \frac{d}{dt}t^{\langle q\rangle}
 =q|t|^{q-1}, \ t\not=0,
\] 
and $|t|^{q-1}\in L^1_{loc}(\R)$ since $q>0$. Furthermore, for $q\geq 1$ ( resp.$0<q<1$ ) the signed power $t^\feq$ is convex (resp. concave) on $(0,+\infty) $ and concave (resp. convex) on $(-\infty,0)$. 

The following lemma was used in the proof of theorem \ref{t:symm of d form and hypersing int} and can be found in \cite[Lemma 3.4]{Korevenpaa-Kuusi-Lindgre19}. 
\begin{lemma}\label{l:precise Lip}
For $p>1$ there is a constant $C_p$ such that  for all $a,\, b\in \R$ we have
\begin{equation}\label{e:precise Lip}
\abs{{a}^\fep-b^\fep}
\leq 
C_p\bigl(\abs{b}+\abs{a-b}\bigr)^{p-2}\abs{a-b}.
\end{equation}

\end{lemma}

\begin{proof} We shall prove the following inequality which is  equivalent to the claimed inequality
\begin{equation}\label{e:precise Lip '}
\abs{\tpwr{a+b}-a^\fep}
\leq 
C_p\bigl(\abs{a}+\abs{b}\bigr)^{p-2}\abs{b}.
\end{equation}

Letting $\phi(t)=\tpwr{a+tb}$ we have $\phi'(t)=(p-1)\,b\, |a+tb|^{p-2}$, hence
\[
\abs{\tpwr{a+b}-a^\fep}\leq (p-1)|b|\,\int_0^1 \left\vert a+tb \right\vert^{p-2}\, dt.
\]
Let us consider the following function $F$ on the punctured plane $\abs{a}+\abs{b}>0$,
\[
F(a,b)\overset{def}{=}\int_0^1 \left\vert\frac{a}{\abs{a}+\abs{b}}+t\frac{b}{\abs{a}+\abs{b}}\right\vert^{p-2}\, dt.
\]
Clearly $F$ is continuous and homogeneous of degree 0. Since it is bounded on the unit circle, the inequality \eqref{e:precise Lip '} follows.
\end{proof}

The following  lemma is useful when testing with positive powers of a solution as in the proof of Caccioppoli's inequality, see  the proof of the intermediate lemma \ref{l: from_caccioppoli_lemma}, or when testing the equation with functions involving level sets of a solution, see the proof of lemma \ref{l: bound on tail lemma}.
\begin{lemma}\label{l: p-coerc_lemma}
Let $p>1$, $w_j\geq 0$ and $0\leq \phi_j\leq 1$, $j=0,1$. For any $\eps>0$ we have
\begin{equation*}
\label{e: I1 integrand_est1 for bound on tail}
\tpwr{{w_0} - {w_1}} \left({w_0} {{\phi_0}^p} - {w_1}{\phi_1^p}\right) 
\geq  (1-\eps)|{w_0}-{w_1}|^{p} \max\{{\phi_0}^p, \phi_1^p\} - {(p-1)}^{p-1}{\eps}^{1-p} |{\phi_0}-\phi_1|^{p} \max\{w_0^p,w_1^p\}.
\end{equation*}
\end{lemma}

\begin{proof} From the homogeneity in the variables and symmetry in the indices of the inequality it is enough to consider the case $$0\leq \phi_1=a <{\phi_0}= 1.$$
From $\abs{1 - {a^p}}  \leq p\abs{1 - a} $ we obtain the first of the following inequalities,
\begin{multline}
\label{e: I1 integrand_est1 for bound on tail in proof}
\tpwr{{w_0} - {w_1}} \left({w_0}  - {w_1}{a^p}\right) =\tpwr{{w_0} - {w_1}}\left[ ({w_0} - {w_1})  +w_1(1-{a^p})\right] 
\\
\geq \abs{{w_0} - {w_1}}^{p}
- p\abs{{w_0} - {w_1}}^{p-1} \abs{{w_1}} \abs{1- a}
\geq  (1-\eps)|{w_0}-{w_1}|^{p} - \left(\frac{\eps}{(p-1)}\right)^{-(p-1)}\abs{1 - a}^{p} |{w_1}|^{p}.
\end{multline}
using Young's inequality $$xy\leq \frac{\delta^{p'} x^{p'}}{p'}+\frac{\delta^{-p}y^{p}}{p}$$ with $x=|{w_0}-{w_1}|^{p}$, $\ y=\abs{{w_1}} \abs{1- a}$ and $\eps=(p-1)\delta^{p'}$. The proof is complete.
\end{proof}

The next lemma is useful when testing with functions involving negative powers of a solution, see \eqref{e: level sets I1}. A different inequality with the same goal  in mind was proven in  \cite[Lemma 2.5]{Kassmann09}.
\begin{lemma}\label{l: p-coerc_lemma log}
Let $p>1$. There is a constant  $B_p$ depending only on $p$, such that,
\begin{equation}\label{e: p-coerc_lemma log}
\tpwr{w_0-w_1}(w_0^{1-p}\phi_0^p- w_1^{1-p}\phi_1^p) \leq - B_p\left\vert \ln  \frac {w_0}{w_1} \right\vert^p \min\{\phi_0^p,\phi_1^p\} + B_p \lvert \phi_0-\phi_1\rvert^p.
\end{equation}
 for any $w_j> 0$ and  $\phi_j\geq 0$, $j=0,1$.
\end{lemma}

\begin{proof} It is enough to show the inequality for $w_0\geq w_1>0$ taking into account the symmetry in the two indices of the claimed inequality.
If $w_0= w_1> 0$ the inequality holds trivially, thus henceforth we assume $w_0> w_1> 0$.

First we consider the case $\phi_0=\phi_1=1$. Let $t=\ln (w_0/w_1)>0$. The negative of the left-hand side of \eqref{e: p-coerc_lemma log} can be written in terms of $t$ as a function $f(t)$ defined for $t\geq 0$, such that
\[
f(t)=(e^t-1)^{p-1}-(1-e^{-t})^{p-1}, \quad f(0)=0, \quad f(t)>0\  \text{ for } \ t>0,\quad \lim_{t\rightarrow \infty} f(t)=+\infty.
\]
It is then enough to show that $A_p\overset{def}{=}\inf_{t> 0} {f(t)}\,{t^{-p}} >0.$  Using the Taylor expansion of the logarithmic function we have for $t\rightarrow 0+$ the inequalities
\[
e^t-1\geq t(1+\frac t2) \quad \text{and}\quad 1-e^{-t}\leq t\left(1-\frac t2+ \frac {t^2}{6}\right),
\]
hence letting $g(t)=(1+\frac t2)^{p-1}-(1-\frac t2+ \frac {t^2}{6})^{p-1}$ we have
\[
 \frac{f(t)}{t^p}\geq \frac {g(t)}{t}.
\]
Since $g'(0)=p-1$ we have $\lim_{t\rightarrow 0+} \frac{f(t)}{t^p} =p-1>0,$ which shows the claimed inequality in the case  $\phi_0=\phi_1=1$,i.e, 
\begin{equation}\label{e: p-coerc_lemma log equal weights}
(w_0-w_1)^{p-1}\left(w_0^{1-p}- w_1^{1-p}\right) \leq - A_p\left( \ln  \frac {w_0}{w_1} \right)^p .
\end{equation}

The general case is a perturbation of the proven inequality. If $\phi_1\geq \phi_0$, then using also \eqref{e: p-coerc_lemma log equal weights} we have
\begin{equation}\label{e: p-coerc_lemma log 1}
(w_0-w_1)^{p-1}(w_0^{1-p}\phi_0^p- w_1^{1-p}\phi_1^p)\leq  (w_0-w_1)^{p-1}(w_0^{1-p}\phi_1^p- w_1^{1-p}\phi_1^p) \leq - A_p\left( \ln  \frac {w_0}{w_1} \right)^p \phi_1^p ,
\end{equation}
which is stronger than the claimed inequality. 

The remaining case, which we consider next, is when $\phi_1< \phi_0$. We recall that if $1\leq p<\infty$ then for any $\eps\in [0,1]$ there exists a constant $c_p$ such that we have the inequality
\begin{equation}\label{l: triangle-ish_lemma}
   (1+\eps)^{p-1}\leq 1+c_p\eps, \quad 0\leq \eps\leq 1, 
\end{equation} 
which follows trivially from $$\lim_{\eps\rightarrow0+}\frac{(1+\eps)^{p-1}- 1}{\eps}=p-1.$$
Thus, using the convexity of $t^p$ we have with $C_p=1+c_p$ the inequality, recall $0\leq\eps\leq 1$,
\[
\phi_0^p\leq (1+\eps)^p\left[\frac{\phi_1^p}{ 1+\eps} + \frac{\eps |\phi_0-\phi_1|^p}{1+\eps}\right] \leq (1+c_p\eps)\phi_1^p +(1+ c_p\eps)\eps^{1-p} |\phi_0-\phi_1|^p\leq (1+c_p\eps)\phi_1^p + C_p\eps^{1-p} |\phi_0-\phi_1|^p.
\]
The above inequality and \eqref{e: p-coerc_lemma log equal weights} show that we have
\begin{multline}\label{e: p-coerc_lemma log 2}
(w_0-w_1)^{p-1}(w_0^{1-p}\phi_0^p- w_1^{1-p}\phi_1^p)\leq
\phi_1^p \left[ (w_0-w_1)^{p-1}\left(w_0^{1-p} - w_1^{1-p}\right)  + c_p\eps\left(\frac{w_0-w_1}{w_0}\right)^{p-1}\right]\\
+C_p\left(\frac{w_0-w_1}{w_0}\right)^{p-1}\eps^{1-p}\lvert \phi_0-\phi_1\rvert^p\\
\leq \phi_1^p \left[ - A_p\left( \ln  \frac {w_0}{w_1} \right)^p +c_p\eps\left(\frac{w_0-w_1}{w_0}\right)^{p-1} \right]
+C_p\left(\frac{w_0-w_1}{w_0}\right)^{p-1}\eps^{1-p}\lvert \phi_0-\phi_1\rvert^p.
\end{multline}
Recalling the assumption $w_0> w_1> 0$, let
\[
\eps=\delta\frac{w_0-w_1}{w_0}\in (0,1],
\]
where $\delta\in(0,1)$ will be determined shortly.  A substitution in \eqref{e: p-coerc_lemma log 2} gives
\begin{multline}\label{e: p-coerc_lemma log 3}
(w_0-w_1)^{p-1}(w_0^{1-p}\phi_0^p- w_1^{1-p})
\leq \phi_1^p \left[ - A_p\left( \ln  \frac {w_0}{w_1} \right)^p +c_p\delta\left(\frac{w_0-w_1}{w_0}\right)^{p} \right]
+C_p\lvert \phi_0-\phi_1\rvert^p\\
\leq  - \frac{A_p}{2}\left( \ln  \frac {w_0}{w_1} \right)^p \phi_1^p
+C_p\lvert \phi_0-\phi_1\rvert^p,
\end{multline}
after taking $\delta$ sufficiently small and taking into account the trivial inequality
\[
\frac{w_0-w_1}{w_0} \leq  \ln  \frac {w_0}{w_1},
\]
which follows, for example, from the convexity of $e^{-t}$, hence its graph lies above the tangent at $(0,1)$.
\end{proof}

The next  lemma is useful when testing the equation with functions involving level sets of a solution, see \eqref{e: I2 + I3 for bound on tail prelim}.

\begin{lemma}\label{l:test fn level set}
Suppose $p>1$, $k\geq 0$, $b\in \R$, and  $0\leq a\leq k$. Then we have
\begin{equation}\label{e:cut-off and sgned power}
  \tpwr{a-b}(a-2k) \geq \frac {k}{2^{p-1}}\,b_{+}^{p-1}-2^{p}k\,b_{-}^{p-1} -(1+2^{p})k^{p}.
\end{equation}
\end{lemma}

\begin{proof}
We begin by taking into account $t_{+}=(-t)_{-}$ which gives
\begin{multline*}
  \tpwr{a-b}(a-2k)
  =(b-a)_{+}^{p-1} (2k-a) -  (b-a)_{-}^{p-1}\,(2k-a) 
  \geq (b-k)_{+}^{p-1}\, k - (b-k)_{-}^{p-1}\,(2k-a) 
\end{multline*}
using that $b-a=b-k +(k-a)\geq b-k$,  $t_{+}$ is  increasing while  $t_{-}$ is  decreasing, hence
\[
(b-a)_{+}\geq  (b-k)_{+}  \qquad\text{and}\quad  (b-a)_{-}\leq (b-k)_{-}.
\]
Now we use the following inequalities, which follow by a direct verification,
\[
b_{+}\leq (b-k)_{+}+ k\quad\text{and}\quad (b-k)_{-}\leq b_{-}+ k.
\]
Thus, we have
\[
b_{+}^{p-1} \leq ( (b-k)_{+}+ k)^{p-1}\leq 2^{p-1}[ (b-k)_{+}^{p-1} + k^{p-1} ]\ \text{and}\  (b-k)_{-}^{p-1}\leq 2^{p-1}[ b_{-}^{p-1} + k^{p-1}],
\]
which completes the proof.
\end{proof}

\appendix

\section{Integral estimates involving balls}
The following lemma will be used frequently to estimate $[\phi]_{s,p,B_r}^p$ with $\phi$ as in lemma \eqref{l: existence of test functions} or $\phi \in C_0^1(\G)$ as in  \eqref{e:MVT cpct supp}. 
\begin{lemma}\label{l: phi^p-kernel-lemma}
Let $p > 0$ and $s \in (0,1)$. 

a) If $U$ is a measurable set and $R>0$ such that $U\subset B_R(h)$ for all $h\in U$ we have 
$$\sup_{h\in U}\int_{U} \frac{|h^{-1}g|^{p}}{|h^{-1}g|^{Q+ps}}\,dg \leq \frac {\sigma_Q}{(1-s)p}R^{p(1-s)}.$$

b) For any ball $B_R(g_0)$ we have $$\int_{B_R(g_0)} \int_{B_R(g_0)}
\frac{|h^{-1}g|^{p}}{|h^{-1}g|^{Q+ps}}\,dgdh \leq \frac {\sigma^2_Q}{Q(1-s)p}R^{Q+p(1-s)}.$$
\end{lemma}

\begin{proof}
The proof follows from  \eqref{e:right inv dusatnce} and the polar coordinates formula  \eqref{e:polar coord}. Indeed,
\begin{multline*}
\sup_{h\in U}\int_{U} \frac{|h^{-1}g|^{p}}{|h^{-1}g|^{Q+ps}}\,dg
\leq \sup_{h\in U} \int_{B_R(h)}\frac{|h^{-1}g|^{p}}{|h^{-1}g|^{Q+ps}}\, dg =  \sup_{h\in U}\int_{L_{h^{-1}}B_R(h)} \frac {|g|^p}{|g|^{Q+sp}}\, dg\\
= \int_{B_{R}(e)} \frac {|g|^p}{|g|^{Q+sp}}\, dg = \frac {\sigma_Q}{(1-s)p}R^{(1-s)p},
\end{multline*}
using \eqref{e:polar coord} to compute the last integral.

The estimate for the double integral follows by noting that for $h\in B_R(g_0)$ we have by the triangle inequality the inclusion $B_R(g_0)\subset B_{2R}(h)$ and the identity $\omega_{Q}=\sigma_Q/Q$ relating the measures of the unit gauge ball and the unit sphere, cf. \eqref{e:polar coord}.

\end{proof}

\begin{lemma}\label{l:Om complement kernel estimate}
Let  $sp>0$ and $\Om$ be an open set. If  $g\in \Om$ and $\delta=\delta_g>0$ such that $B_{\delta}(g)\subset\Om$ we have 
\[
\int_{\bG\setminus \Om} \frac{dh}{|h^{-1}g|^{Q+ps}} \leq \frac{\sigma_Q}{sp}\delta_g^{-sp}.
\]
\end{lemma}

\begin{proof}  From  $\bG\setminus \Om\subset \bG\setminus B_{\delta}(g)$ and a left translation we have
  \begin{equation}\label{e:kernel estimate outside ball}
  \int_{\bG\setminus \Om} \frac{dh}{|h^{-1}g|^{Q+ps}} \leq \int_{\bG\setminus  B_{\delta}(g)} \frac{dh}{|h^{-1}g|^{Q+ps}} = \int_{\bG\setminus  B_{\delta}(e)} \frac{dh}{|h|^{Q+ps}}=\frac{\sigma_Q}{sp}\delta^{-sp}
  \end{equation}
  taking into account \eqref{e:polar coord} in  the last equality.  
\end{proof}

\begin{cor}\label{c: ker_est}
Let $p > 0$, $0 < s < 1$ and $0 < a < 1$.  If $g_{0} \in \G$, $B_{R} = B_R(g_{0})$, $B_{r} = B(g_{0})$ with $r=aR$,  then we have

\begin{equation}\label{e:kernel Br and BR}
\sup_{g\in B_r(g_0)} \int_{\bG\setminus B_R(g_0)}  \frac{dh}{|h^{-1}g|^{Q+sp}}  \leq \frac{\sigma_Q}{sp(1-a)^{sp} R^{sp}}=\frac{\sigma_Q}{sp(R-r)^{sp} }.
\end{equation}
and
$$ \int_{B_{r}}  \int_{\G\setminus B_{R}}\frac{dh}{|h^{-1}g|^{Q+ps}} \, dg \leq \frac{\sigma_{Q}^{2}\,a^Q\, R^{Q-sp}}{sp Q (1-a)^{sp}}
=\frac{ \sigma_{Q}^{2} }{spQ}\frac{r^{Q}}{(R-r)^{sp}}.$$
\end{cor}

\begin{proof}  The proof of \eqref{e:kernel Br and BR} follows from lemma \ref{l:Om complement kernel estimate} noting that for $g\in B_{aR}$ we have $B_{\delta}(g)\subset B_R$ for $\delta=(1-a)R$. The second inequality follows by integrating the first,
\[
\int_{B_{aR}}\int_{\bG\setminus B_R}  \frac{dh}{|h^{-1}g|^{Q+sp}} \, dg \leq \frac{\sigma_Q}{sp}\delta^{-sp} |B_{aR}|= \frac{a^Q\sigma_{Q}^{2}}{Qsp(1-a)^{sp}}R^{Q-sp}.
\]
\end{proof}

\begin{lemma}
\label{l: ker_est}
Let $p > 0$, $0 < s < 1$ and $0 < a < 1$.  Let $g_{0} \in \G$, $B_{R} = B(g_{0},R)$ and $B_{aR} = B(g_{0}, aR)$. We have
\begin{equation}\label{l: ker_est 1}
\sup_{h\in \bG\setminus B_R} \int_{B_{aR}} \frac{dg}{|h^{-1}g|^{Q+sp}}\leq \frac{a^Q}{(1-a)^{Q+sp}}\frac{|B_R|}{|h^{-1}g_0|^{Q+sp}}\leq \frac{\sigma_Q a^{Q}}{Q(1-a)^{Q+sp}}\frac{1}{R^{sp}}.
\end{equation}
\end{lemma}

\begin{proof} By a left translation by $g_0$ we reduce to the case $g_0=e$. The triangle inequality gives
\[
|h^{-1}g|\geq |h|-|g|\geq |h|-a|h|=(1-a)|h|,
\]
which gives \eqref{l: ker_est 1} taking into account \eqref{e:polar coord}.  
\end{proof}

\section{Poincar\'e and Sobolev inequalities}\label{s: functional}
The next results essentially follow from known arguments. We state and prove them in our context for the sake of completeness.
\subsection{Poincar\'e Inequalities}

We  have the following  fractional Poincar\'e inequality, which can be proven as in  \cite[(4.2)]{Ming-2003}.
\begin{theorem}\label{t:Poincare} Let $1\leq p<\infty$ and $0<s<1$. Let $g_{0} \in \G$, $R  > 0$ and $B_{R}=B(g_0,R)$.  There exists a constant $C=C(p,Q,s)$ such that for  $u \in \Wsploc$ we have
\begin{equation}\label{e: Poincare}
 \int_{B_{R}} |u-u_{R}|^{p} dg
\leq  C R^{sp} D_R(u,u),\quad D_R(u,u)=\int_{B_{R}} \int_{B_{R}} \frac{|u(g) - u(h)|^{p}}{|h^{-1}g|^{Q+sp}} dg dh, \quad\text{cf. \eqref{e:local forms balls} .}
\end{equation}

\end{theorem}

The following result is adapted from \cite[Lemma 2.2]{Bra-Par-2015}.

\begin{lemma}[Localized Poincare Inequality]
\label{l: localized Poincare}
Let $0 < p < \infty$, $0 < s < 1$, $g_{0} \in \G$ and $0<r<R $. Let $B_{r}$ and $B_{R}$ be  the balls centered at $g_{0}$ of radii $r$ and $R$ respectively.  For each $u \in \Dsp(\bG)$ which is compactly supported in $B_{r}$ we have
\begin{equation}\label{e: localized Poincare}
\frac{Q \omega_{Q}}{2^{Q+p}} \left(\frac{r}{R} \right)^{Q} \left( \frac{R-r}{r} \right) R^{-sp} \int_{B_{r}} |u|^{p} dg
\leq  D_R(u,u).
\end{equation}
Where $\omega_{Q}$ is the measure of the unit gauge ball in $\bG$.
\end{lemma}

\begin{proof}
Since $u$ is supported in $B_{r}$ and $R > r$, for any $h \in B_{R}$ such that $u(h) = 0$ we have
$$ |u(g)|^{p}
= |u(g) - u(h)|^{p}
\leq   2^{Q+sp}R^{Q+sp} \frac{|u(g) - u(h)|^{p}}{|h^{-1}g|^{Q+sp}}, \quad g \in B_{R}, $$
since $g,h \in B_{R}$ implies $|h^{-1}g| \leq 2R$. Integration with respect to $h$ over $\{h \in B_{R} : u(h) = 0\}$ then gives
\begin{equation}\label{e: localized Poincare intermediate}
\left|\{g \in B_{R} : u(g) = 0 \}\right|\ |u(g)|^{p}
\leq \int_{B_{R}} 2^{Q+sp}R^{Q+sp} \frac{|u(g) - u(h)|^{p}}{|h^{-1}g|^{Q+sp}} dh.
\end{equation}
Since $\supp u \subset B_{r}$, $u$ vanishes on $B_{R} \setminus B_{r}$, thus  we have
\[
\left|\{g \in B_{R} : u(g)
= 0 \}\right| \geq |B_{R}| - |B_{r}| 
= \omega_{Q}(R^{Q} - r^{Q}) \geq \omega_{Q} Q(R-r)  r^{Q-1}.
\]
With this, \eqref{e: localized Poincare intermediate} becomes
$$ {\omega_{Q} }(R-r) Q r^{Q-1} |u(g)|^{p} \leq \int_{B_{R}} 2^{Q+sp}R^{Q+sp} \frac{|u(g) - u(h)|^{p}}{|h^{-1}g|^{Q+sp}} dh.$$
Dividing by the constant on the right and integrating with respect to $g$ then gives the result
$$\frac{Q \omega_{Q}}{2^{Q+sp}} \left(\frac{r}{R} \right)^{Q} \left( \frac{R-r}{r} \right) R^{-sp} \int_{B_{r}} |u|^{p} dg
\leq \int_{B_{R}} \int_{B_{R}} \frac{|u(g) - u(h)|^{p}}{|h^{-1}g|^{Q+sp}} dg dh=D_R(u,u).$$

\end{proof}

\subsection{Sobolev embedding}
The following localized Sobolev embedding is adapted from \cite[Lemma 2.3]{Bra-Par-2015}. The integrals $D_R$ and $I_R$ we use were defined in  \eqref{e:local forms balls}.

\begin{theorem}[Localized Sobolev Embedding]
\label{t: Local Sobolev embedding}
Let $0<\rho < R $,  $0 < s < 1$ and $1 < p < Q/s$.  There is a constant $c = c(p,Q,s)$ such that for $u \in \Dsp(\bG)$ which is compactly supported in a ball $B_{\rho}=B_{\rho}(g_0)$, $g_0\in\bG$, we have
\begin{equation*}\left( \int_{B_{\rho}} |u|^{p^{*}(s)} dg \right)^{\frac{p}{p^{*}(s)}}
\leq c \left( \left(\frac{R}{\rho} \right)^{Q} \left( \frac{\rho}{R-\rho} \right) \left(\frac{R}{R-\rho}\right)^{sp} + 1\right) D_{R}(u,u).
\end{equation*}

\end{theorem}

\begin{proof} By the translation invariance of the involved quantities, it is enough to consider the case $g_0=e$.
Since $u  \in {\Dsp}(\G)$  we may apply the global Sobolev inequality \eqref{t: Global Sobolev} to obtain $c$ such that
\begin{equation*}
\left( \int_{B_{\rho}} |u(g)|^{p^{*}(s)} dg \right)^{\frac{p}{p^{*}(s)}}
\leq \left( \int_{\G} |u(g)|^{p^{*}(s)} dg \right)^{\frac{p}{p^{*}(s)}} \\
\leq c \int_{\G} \int_{\G} \frac{\abs{u(g) - u(h)}^{p}}{|h^{-1}g|^{Q+ps}} dg dh
\end{equation*}
Decomposing the integral on the right and using that
$u$ is supported in the ball $B_{\rho} \subset B_{R}$ we obtain
\begin{equation}
\label{e: local sobolev initial}
\left( \int_{B_{\rho}} |u|^{p^{*}(s)} dg \right)^{\frac{p}{p^{*}(s)}}   \\
= 2cI_R(u,u)+ c D_R(u,u).
\end{equation}
Again using that $u$ is supported in $B_{\rho}$, the first integral on the right becomes
\begin{equation}
\label{e: local sobolev intermediate}
 I_R(u,u)= \int_{B_{\rho}} \int_{\G \setminus B_{R}}  \frac{|u(g)|^{p}}{|h^{-1}g|^{Q+ps}} dh dg  = \int_{B_{\rho}} |u(g)|^{p} \int_{\G \setminus B_{R}}  \frac{1}{|h^{-1}g|^{Q+ps}} dh dg.
\end{equation}
Now observe that $(\G \setminus B_{R}) \subset (\G \setminus B_{R-\rho}(g))$ for each $g \in B_{\rho}$. Thus, from corollary \ref{c: ker_est} we have
$$
 I_R(u,u)\leq \frac{\omega_{Q}}{ps(R-\rho)^{ps}} \int_{B_{\rho}} |u(g)|^{p}  dg.$$
The localized Poincare inequality, lemma \ref{l: localized Poincare}, and  the last inequality give
\[
 I_R(u,u)\leq \frac{\omega_{Q}}{ps(R-\rho)^{ps}} \frac{2^{Q+p}}{Q \omega_{Q}} \left(\frac{R}{\rho} \right)^{Q} \left( \frac{\rho}{R-\rho} \right) R^{sp} \,D_R(u,u).
\]
Taking into account \eqref{e: local sobolev initial}, we obtain
\begin{equation*}
\left( \int_{B_{\rho}} |u(g)|^{p^{*}(s)} dg \right)^{\frac{p}{p^{*}(s)}}
\leq c \left( \left(\frac{R}{\rho} \right)^{Q} \left( \frac{\rho}{R-\rho} \right) \left(\frac{R}{R-\rho}\right)^{sp} + 1\right)\,  D_R(u,u) 
\end{equation*}
with a constant $c=c(Q,p,s)$ noting that $\frac{R}{R-\rho} > 1$.  The proof is complete.
\end{proof}

\section{A Caccioppoli type inequality}\label{ss:Caccioppoli's inequality}

 The following lemma contains an estimate for one of the terms that appears in the proof of  lemma \ref{l: intro of inf}, cf. \eqref{e: J1 from Caccioppoli}, and also in the Caccioppoli's inequality theorem \ref{t: Caccioppoli}. We state and provide a proof of this intermediate result following the arguments in the proofs of \cite[Theorem 1.4]{DiC-Kuu-Pal-2016},  \cite[proposition 8.5]{Cozzi17}, and our lemma \ref{l: p-coerc_lemma}. As before, we use the notation from \eqref{e:local forms balls} for the form $D_{\rho}(.,.)= D_{B_\rho}(., .)$.
  For simplicity, we state the lemma for $w_{+}$, but a version of inequality \eqref{l: from_caccioppoli_lemma} for $w_{-}$ follows  by replacing $u$ with $-u$ and $l$ with $-l$, taking also into account that $D_{\rho}(u,\phi)$ is an odd function of $u$ and the relation  $w_-=(-w)_+$ since  $   w_{\pm} =\max\{\pm w,0 \}$.

\begin{lemma}
\label{l: from_caccioppoli_lemma}
Let $p > 1$, $s \in (0,1)$. Let $\phi\in\mathcal{C}^{\infty}_o(B_{\rho})$ be a non-negative function, where $B_\rho=B_{\rho}(g_0)$ with $g_0\in\bG$. For any constant $l\in\R$ and  $u\in \Wsploc$ we have with  $w= u-l$ the inequality
\begin{multline}\label{e: from_caccioppoli_estimate_lemma}
D_{\rho}(u,w_{+}\phi^p)
\ +  4\left(1 +  {(p-1)^{p-1}}{ 2^{p-1}}\right) \int_{B_{\rho}}\int_{B_{\rho}}\lvert \phi (g)-\phi (h)\rvert^p \, w_{+}^p(g) \, K(g,h) dg dh
\\
 \geq   \frac {1}{2^{p-2} }D_\rho(w_{+}\phi,w_{+}\phi)\ + \   \int_{B_{\rho}}\int_{B_{\rho}} {w_{-}}^{p-1}(h){w_{+}}(g)\phi^{p}(h)\, K(g,h) dg dh
\end{multline}
\end{lemma}

\begin{proof}[Proof of Lemma \ref{l: from_caccioppoli_lemma}]
We begin with the proof of some point-wise estimates.  If  $u(g) > u(h)$ we have the following identities
\begin{multline}\label{e: from_caccioppoli_cases}
|u(g) - u(h)|^{p-2} \paren{u(g) - u(h)}\paren{{w_{+}}(g) \phi^{p}(g) - {w_{+}}(h)\phi^{p}(h)}
\\
= \begin{cases}
\paren{{w_{+}}(g) - {w_{+}}(h)}^{p-1} \paren{{w_{+}}(g)\phi^{p}(g) - {w_{+}}(h)\phi^{p}(h)},  &  l< u(h) < u(g) ;
\\
\paren{ {w_{+}}(g) + {w_{-}}(h)  }^{p-1}{w_{+}}(g)\phi^{p}(g),  & u(h) \leq  l <u(g); \\
0, &  u(h)<u(g)<l.
\end{cases}
\end{multline}
In the first case of  \eqref{e: from_caccioppoli_cases},  where $ w_{+}(g) > w_{+}(h)$, lemma \ref{l: p-coerc_lemma} gives for  $\eps>0$ the estimate
\begin{multline}\label{e: from p-coerc lemma}
({w_{+}}(g)-{w_{+}}(h))^{p-1}({w_{+}}(g)\phi^{p}(g)- {w_{+}}(h)\phi^{p}(h)) \\
\geq (1-\eps)\left\vert{w_{+}}(g)-{w_{+}}(h)\right\vert^p\max \{\phi (g),\phi (h)\}^p
- C_p \eps^{1-p}\lvert \phi (g)-\phi (h)\rvert^p \max \{{w_{+}}(g),{w_{+}}(h)\}^p,
\end{multline}
where $C_p=(p-1)^{p-1}$ is the constant from lemma \ref{l: p-coerc_lemma}.  On the other hand, in the second case of  \eqref{e: from_caccioppoli_cases} where $u(h) < l \leq u(g)$ we use \eqref{e: triangle with p},  
thus, noting that ${w_{+}}(h)=0$ 
and $p>1$, we obtain
\begin{equation}\label{e: from simple ineq}
\paren{ {w_{+}}(g) + {w_{-}}(h)  }^{p-1}{w_{+}}(g)\geq  \frac 12 \left[  |{w_{+}}(g) - {w_{+}}(h)  |^{p} + {w_{-}}^{p-1}(h){w_{+}}(g)\right].
\end{equation}
From the symmetry of  the first integral in \eqref{e: from_caccioppoli_estimate_lemma} with respect to $g$ and $h$,  a decomposition of the domain of integration based on the cases in \eqref{e: from_caccioppoli_cases},    inequalities \eqref{e: from p-coerc lemma} and \eqref{e: from simple ineq},  noting that $w_{+}(g)=w_{+}(h)=0$ in the third case of \eqref{e: from_caccioppoli_cases},  $w_{-}(h)=0$ in the first  case, while  $w_{+}(g)=0$ in the third case,  we obtain the next  inequality,
\begin{multline}\label{e:key for caccioppoli 1}
\frac 12 D_{\rho}(u,w_{+}\phi^p)
\geq  (1-\eps) \int_{B_{\rho}}\int_{B_{\rho}}\left\vert{w_{+}}(g)-{w_{+}}(h)\right\vert^p\max \{\phi (g),\phi (h)\}^p\, K(g,h) dg dh\\
-C_p \eps^{1-p}\int_{B_{\rho}}\int_{B_{\rho}}\lvert \phi (g)-\phi (h)\rvert^p \max \{{w_{+}}(g),{w_{+}}(h)\}^p \, K(g,h) dg dh\\
+  \frac 12 \int_{B_{\rho}}\int_{B_{\rho}}\left[ |{w_{+}}(g) - {w_{+}}(h)  |^{p} + {w_{-}}^{p-1}(h){w_{+}}(g)\right]\phi^{p}(h)\, K(g,h) dg dh.
\end{multline}
Due to the assumption  $u\in \Wsploc$ the terms in  \eqref{e:key for caccioppoli 1} are finite. Thus, with $\eps=\frac 12$ in  the  last inequality we have
\begin{multline}\label{e:key for caccioppoli 2}
\frac 12 D_{\rho}(u,w_{+}\phi^p)
+  C_p 2^{p-1}\int_{B_{\rho}}\int_{B_{\rho}}\lvert \phi (g)-\phi (h)\rvert^p \max \{{w_{+}}(g),{w_{+}}(h)\}^p \, K(g,h) dg dh
\\
\geq \int_{B_{\rho}}\int_{B_{\rho}} |{w_{+}}(g) - {w_{+}}(h)  |^{p} \phi^{p}(h)\, K(g,h) dg dh
  +\frac 12\int_{B_{\rho}}\int_{B_{\rho}} {w_{-}}^{p-1}(h){w_{+}}(g)\phi^{p}(h)\, K(g,h) dg dh.
\end{multline}
 Using  the triangle inequality and  \eqref{e: triangle with p} we obtain
\begin{multline}\label{e:equiv caccioppoli}
\left\vert w_{+} \phi(g) - w_{+}\phi (h) \right\vert^p
\leq \Bigl( \left\vert  \phi(g) - \phi (h)  \right\vert w_{+}(g)
+  \left\vert  {w_{+}}(g) - {w_{+}}(h)  \right\vert \phi (h)\Bigr)^p \\
\leq  2^{p-1} \left\vert  \phi(g) - \phi (h)  \right\vert ^p  \max \{{w_{+}}(g),{w_{+}}(h)\} ^p  + 2^{p-1} \left\vert  {w_{+}}(g) - {w_{+}}(h)  \right\vert ^p \phi (h)^p.
\end{multline}
From \eqref{e:key for caccioppoli 2} and \eqref{e:equiv caccioppoli}, taking also in account the symmetry with respect to the variables in the integral containing $\max \{{w_{+}}(g),{w_{+}}(h)\}^p\leq {w_{+}}(g)+{w_{+}}(h) $, we obtain  the inequality
\begin{multline}\label{e:key for caccioppoli 3}
\frac 12 D_{\rho}(u,w_{+}\phi^p)
+ 2\left( 1 +C_p2^{p-1}\right)\int_{B_{\rho}}\int_{B_{\rho}}\lvert \phi (g)-\phi (h)\rvert^p \,{w_{+}}^p(g) \, K(g,h) dg dh
\\
\geq  \frac {1}{2^{p-1} }D_\rho(w_{+}\phi,w_{+}\phi)\ 
 + \ \frac 12 \int_{B_{\rho}}\int_{B_{\rho}} {w_{-}}^{p-1}(h){w_{+}}(g)\phi^{p}(h)\, K(g,h) dg dh,
\end{multline}
which implies inequality \eqref{e: from_caccioppoli_estimate_lemma}.
\end{proof}

\begin{lemma}
\label{l:I_r estimate}
Let $p > 1$, $s \in (0,1)$ and $u\in \Wsploc$. Let $\phi\in\mathcal{C}^{\infty}_o(B_{r})$ be a non-negative function, where $B_\rho=B_{r}(g_0)$ with $g_0\in\bG$. For any non-negative  constant $l\in\R$   and  $w= u-l$  we have the inequality
\begin{equation}\label{e:I_r estimate}
I_{r}(u,w_{+}\phi^p) \geq  -2^{p-1}\int_{\G\setminus B_{r}} \int_{B_{r}}  w_{+}^{p-1}(h)w_{+}(g)\phi^p(g) K(g,h)dgdh,
\end{equation}
where the form $I_{r}(.,.)= D_{B_r}(.,\, .)$ is defined in  \eqref{e:local forms balls}.
\end{lemma}

\begin{proof}[Proof of lemma \ref{l:I_r estimate}]
We begin the estimate of the  integral $I_{r}(.,.)$ by a use of the trivial inequality $$a^\fe{p-1}=|a|^{p-2}a \geq-(-a)_{+}^{p-1},$$   
which gives the point-wise estimate
\begin{multline}\label{e:non-local Cacc}
\tpwr{u(g)-u(h)}w_{+}(g)=\tpwr{w(g)-w(h)}w_{+}(g)\\
\geq -\left(w_{+}(h) -w_{-}(h)+w_{-}(g)   -w_{+}(g) \right)_{+}^{p-1}w_{+}(g)\\
\geq -\left[w_{+}(h) +w_{-}(g)  \right]_{+}^{p-1}w_{+}(g)\geq - 2^{p-1} w_{+}^{p-1}(h)w_{+}(g),
\end{multline}
taking into account $w_{-}(g)\, w_{+}(g) =0$ and \eqref{e: triangle with p}.
Integrating inequality \eqref{e:non-local Cacc} we obtain \eqref{e:I_r estimate}.
\end{proof}
Again, let us observe that a version of the above lemma holds for $w_{-}$ taking into account that $I_r(.\, , .)$ is odd in each of the variables.

In the next theorem we give a Caccioppoli type inequality.
Other general versions of Caccioppoli's inequality in the Euclidean case can be found in  \cite[proposition 3.5]{Bra-Par-2015} and \cite[Theorem 3.2]{KuusiMingSire15}.

\begin{theorem}[Caccioppoli's Inequality]
\label{t: Caccioppoli}
Let $p \in (1,\infty)$. Suppose $u\in \Wsploc\cap\tailspace$ is a sub-solution (resp. super-solution) to \eqref{e: Main equation}. For any ball  $B_{r} = B_{r}(g_{0})$, nonnegative function $\phi \in C_{0}^{\infty}(B_{r})$, constant $l\in\R$, and   $w= u-l$, we have the inequality
\begin{multline}\label{e:caccioppoli}
\pm\int_{\bG} F(g) w_{\pm}\phi^p(g) \,dg +  C_p \int_{B_{r}}\int_{B_{r}}\lvert \phi (g)-\phi (h)\rvert^p \, w_{\pm}^p(g) \, K(g,h) dg dh 
\\
\geq  \frac{1}{2^{p-2}} D_r(w_{\pm}\phi,w_{\pm}\phi)
  + \int_{B_{r}}\int_{B_{r}}{w_{\mp}}^{p-1}(g){w_{\pm}}(h)\phi^{p}(h)\, K(g,h) dg dh \\
   - \frac{1}{2^{p-2}} \int_{\G\setminus B_{r}} \int_{B_{r}}w_{\pm}^{p-1}(h)w_{\pm}(g)\phi^p(g) K(g,h)dgdh,
\end{multline}
where $C_p$ is a constant depending on $p$, and the upper signs are used in the case of a sub-solution, while the lower signs are used in the case of a super-solution.
\end{theorem}

\begin{proof}[Proof of Theorem \ref{t: Caccioppoli} ]
Taking into account the comments in the paragraphs preceding lemma \ref{e: from_caccioppoli_estimate_lemma}  it is enough to consider the case of a sub-solution and prove the case of the upper signs.  By the definition of a sub-solution, cf. \eqref{e:weak solution}, with $\eta=w_{+}\phi^p$ we obtain
\begin{equation}\label{e: caccioppoli initial decomp}
\int_{\bG} F(g) w_{+}\phi^p(g) \,dg \geq D_{r}(u,w_{+}\phi^p)+2I_{r}(u,w_{+}\phi^p),
\end{equation}
where $D_{r}(.,.)$ and $I_{r}(.,.)$ were defined in \eqref{e:local forms balls}.  The proof follows from \eqref{e: from_caccioppoli_estimate_lemma} and \eqref{e:I_r estimate}.
\end{proof}

\bibliographystyle{plain} 
\bibliography{MasterBib}

\end{document}